\documentclass{amsart}

\newtheorem{theorem}{Theorem}[section]
\newtheorem{lemma}[theorem]{Lemma}

\newtheorem{proposition}[theorem]{Proposition}
\newtheorem{corollary}[theorem]{Corollary}

\theoremstyle{definition}
\newtheorem{definition}[theorem]{Definition}

\theoremstyle{remark}
\newtheorem{remark}[theorem]{Remark}

\usepackage{graphicx}
\usepackage{color}
\usepackage{amsmath}
\usepackage{amsfonts}
\usepackage{amssymb}
\newcommand{\be}{\begin{equation}}
\newcommand{\ee}{\end{equation}}
\newcommand{\ben}{\begin{equation*}}
\newcommand{\een}{\end{equation*}}

\newcommand{\al}{\alpha}

\newcommand{\bet}{\beta}

\newcommand{\Om}{\Gamma}

\newcommand{\om}{\omega}

\newcommand{\we}{{\stackrel{\scriptscriptstyle{W}}{\Gamma}}\phantom{}}

\newcommand{\D}{\mathcal{D}}

\newcommand{\dz}{\wedge}

\newcommand{\ba}{\begin{array}}

\newcommand{\ea}{\end{array}}

\newcommand{\beq}{\begin{eqnarray}}

\newcommand{\eeq}{\end{eqnarray}}

\newtheorem{lm}{lemma}

\newtheorem{thee}{theorem}

\newtheorem{proo}{proposition}

\newtheorem{co}{corollary}

\newtheorem{rem}{remark}

\newtheorem{deff}{definition}

\newcommand{\bd}{\begin{deff}}

\newcommand{\ed}{\end{deff}}

\newcommand{\bl}{\begin{lm}}

\newcommand{\el}{\end{lm}}

\newcommand{\bp}{\begin{proo}}

\newcommand{\ep}{\end{proo}}

\newcommand{\bt}{\begin{thee}}

\newcommand{\et}{\end{thee}}

\newcommand{\bc}{\begin{co}}

\newcommand{\ec}{\end{co}}

\newcommand{\brm}{\begin{rem}}

\newcommand{\erm}{\end{rem}}

\newcommand{\der}{{\rm d}}

\newcommand{\sgn}{\mathrm{sgn}}

\usepackage{t1enc}
\def\frak{\mathfrak}

\newcommand{\newc}{\newcommand}

\renewcommand{\exp}{\operatorname{exp}}

\let\ccdot\cdot
\def\cdot{\hbox to 2.5pt{\hss$\ccdot$\hss}}

\newc{\aR}{\mbox{\boldmath{$ R$}}}
\newc{\aS}{\mbox{\boldmath{$ S$}}}
\newc{\aT}{\mbox{\boldmath{$ T$}}}
\newc{\aW}{\mbox{\boldmath{$ W$}}}

\newc{\aK}{\mbox{\boldmath{$ K$}}}
\newc{\aL}{\mbox{\boldmath{$ L$}}}

\usepackage{amssymb}
\usepackage{amscd}

\newcommand{\hook}{\raisebox{-0.35ex}{\makebox[0.6em][r]
{\scriptsize $-$}}\hspace{-0.15em}\raisebox{0.25ex}{\makebox[0.4em][l]{\tiny
 $|$}}}

\newtheorem*{thintro}{Theorem}
\newtheorem*{defintro}{Definition}
\newtheorem*{propintro}{Proposition}

\newc{\obstrn}[2]{B^{#1}_{#2}}

\newcommand{\rpl}                         
{\mbox{$
\begin{picture}(12.7,8)(-.5,-1)
\put(0,0.2){$+$}
\put(4.2,2.8){\oval(8,8)[r]}
\end{picture}$}}

\newcommand{\lpl}                         
{\mbox{$
\begin{picture}(12.7,8)(-.5,-1)
\put(2,0.2){$+$}
\put(6.2,2.8){\oval(8,8)[l]}
\end{picture}$}}

\usepackage{ifthen}

\newcommand{\bbR}{\mathbb{R}}
\newcommand{\bbN}{\mathbb{N}}

\newcommand{\sog}{\mathbf{SO}}
\newcommand{\cog}{\mathbf{CO}}

\newcommand{\slg}{\mathbf{SL}}
\newcommand{\glg}{\mathbf{GL}}
\newcommand{\og}{\mathbf{O}}
\newcommand{\coa}{\frak{co}}
\newcommand{\soa}{\frak{so}}
\newcommand{\sla}{\frak{sl}}
\newcommand{\gla}{\frak{gl}}

\newcommand{\spg}{\mathbf{Sp}}
\newcommand{\sug}{\mathbf{SU}}

\newc{\tensor}[1]{#1}
\newc{\Mvariable}[1]{\mbox{#1}}
\newc{\down}[1]{{}_{#1}}
\newc{\up}[1]{{}^{#1}}

\newc{\JulyStrut}{\rule{0mm}{6mm}}
\newc{\midtenPan}{\mbox{\sf S}}
\newc{\midten}{\mbox{\sf T}}
\newc{\midtenEi}{\mbox{\sf U}}
\newc{\ATen}{\mbox{\sf E}}
\newc{\BTen}{\mbox{\sf F}}
\newc{\CTen}{\mbox{\sf G}}

\def\sideremark#1{\ifvmode\leavevmode\fi\vadjust{\vbox to0pt{\vss
 \hbox to 0pt{\hskip\hsize\hskip1em
 \vbox{\hsize3cm\tiny\raggedright\pretolerance10000
 \noindent #1\hfill}\hss}\vbox to8pt{\vfil}\vss}}}%

\newcommand{\bgw}{{\textstyle \bigwedge}}
\newcommand{\bgs}{{\textstyle \bigodot}}
\newcommand{\bgt}{{\textstyle \bigotimes}}

\numberwithin{equation}{section}

\newcommand{\bma}{\begin{pmatrix}}
\newcommand{\ema}{\end{pmatrix}}

\newcounter{romenumi}
\newcommand{\labelromenumi}{(\roman{romenumi})}

\newcommand{\ten}{\Upsilon}

\renewcommand{\Pi}{\tilde{\Gamma}}

\newcommand{\Gl}{\glg(2,\bbR)}
\newcommand{\gl}{\gla(2,\bbR)}
\newcommand{\Rv}{R_{\rm v}}
\allowdisplaybreaks

\begin{document}
\title{${\bf GL}(2,\bbR)$ geometry of ODE's}

\author{Micha\l~ Godli\'nski} 
\address{Instytut Matematyczny Polskiej Akademii Nauk, ul. Sniadeckich 8, Warszawa, Poland}
\address{Instytut Fizyki Teoretycznej, Uniwersytet Warszawski, ul. Hoza 69, Warszawa, Poland}
\email{godlinski@impan.gov.pl}

\author{Pawe\l~ Nurowski} 
\address{Instytut Fizyki Teoretycznej,
Uniwersytet Warszawski, ul. Hoza 69, Warszawa, Poland}
\email{nurowski@fuw.edu.pl} 
\thanks{This research was supported by the KBN grant 1 P03B 07529}
\date{1 October 2007}

\begin{abstract}
We study five dimensional geometries associated with the 5-dim\-ensio\-nal
irreducible representation of $\Gl$. These are special Weyl
geometries in signature $(3,2)$ having the structure group reduced 
from $\cog(3,2)$ to $\Gl$. The reduction is obtained by means of 
a conformal class of totally symmetric 3-tensors.
Among all $\Gl$ geometries we distinguish a subclass which we term 
`nearly integrable $\Gl$ geometries'. These define a 
unique $\gl$ connection which has totally skew symmetric torsion. This
torsion splits onto the $\Gl$ irreducible components having respective
dimensions 3 and 7.  

We prove that on the solution space of every 5th order ODE satisfying
certain three nonlinear differential conditions there exists a 
nearly integrable $\Gl$ geometry such that the skew symmetric torsion
of its unique $\gl$ connection is very special. In contrast to an
arbitrary nearly integrable $\Gl$ geometry, it belongs to the 
3-dimensional irreducible representation of $\Gl$. 
The conditions for the existence of the structure are lower order
equivalents of the Doubrov-Wilczynski conditions found recently 
by Boris Doubrov \cite{dub}.

We provide nontrivial examples of 5th order ODEs satisfying the three
nonlinear differential conditions, which in turn provides examples of
inhomogeneous $\Gl$ geometries in dimension five, with torsion in
$\bbR^3$.

We also outline the theory and the basic properties of $\Gl$
geometries associated with $n$-dimensional irreducible representations
of $\Gl$ in $6 \leq n \leq 9$. In particular we give conditions for an
$n$th order ODE to possess this geometry on its solution space.
\vskip5pt\centerline{\small\textbf{MSC classification}: 53A40, 53B05,
  53C10, 34C30}\vskip15pt
\end{abstract}
\maketitle

\tableofcontents

\section{Introduction}

Let us start with an elementary algebraic geometry in $\bbR^3$. 
Every \emph{point} on a curve $(1,x,x^2)$ in $\bbR^3$ defines a 
\emph{straight line} passing
through the origin in the dual space $(\bbR^3)^*$ via the relation:
\begin{eqnarray}
\theta^0 +2\theta^1 x+\theta^2 x^2&=&0\label{uw1}\\
\theta^1+\theta^2 x&=&0.\nonumber
\end{eqnarray}
Here $(\theta^0,\theta^1,\theta^2)$ parametrize points of
$(\bbR^3)^*$. When moving along the curve $(1,x,x^2)$ in $\bbR^3$, the
corresponding lines in the dual space $(\bbR^3)^*$ sweep out a ruled
\emph{surface} there, which is the cone
\be
(\theta^1)^2-\theta^0\theta^2=0\label{pn1}\ee
with the tip in the origin. The points $(\theta^0,\theta^1,\theta^2)$
lying on this cone may be thought as those, and only those, which 
admit a common root $x$ for the pair of equations (\ref{uw1}). A
standard method for determining if two polynomials have a common root is
to equate to zero their \emph{resultant}. In the case of equations
(\ref{uw1}) the resultant is:
$$R_3=\det\bma \theta^0&2\theta^1&\theta^2&0&0\\
0&\theta^0&2\theta^1&\theta^2&0\\0&0&\theta^0&2\theta^1&\theta^2\\
\theta^1&\theta^2&0&0&0\\
0&\theta^1&\theta^2&0&0\\
\ema.$$
It vanishes if and only if condition (\ref{pn1}) holds.

Before passing to $\bbR^n$ with general $n\geq 3$, it is instructive to repeat the above
considerations in the cases of $n=4$ and $n=5$.  

A \emph{point} on a curve $(1,x,x^2,x^3)$ in $\bbR^4$ defines a
\emph{plane} passing through the origin in the dual space $(\bbR^4)^*$ 
via the relation:
\begin{eqnarray}
\theta^0 +3\theta^1 x+3\theta^2 x^2+\theta^3 x^3&=&0\label{uw2}\\
\theta^1+2\theta^2 x+\theta^3 x^2&=&0.\nonumber
\end{eqnarray}
Now $(\theta^0,\theta^1,\theta^2,\theta^3)$ parametrize points of
the dual $(\bbR^4)^*$ and when moving along the curve $(1,x,x^2,x^3)$
in $\bbR^4$, the
corresponding planes in $(\bbR^4)^*$ sweep out a ruled
\emph{hypersurface} there, which is defined by the vanishing of the resultant
of the two polynomials defined in (\ref{uw2}). This is given by 
\be
-3(\theta^1)^2 (\theta^2)^2 + 4\theta^0 (\theta^2)^3 + 4
(\theta^1)^3\theta^3 -6 \theta^0 \theta^1\theta^2 \theta^3 +
(\theta^0)^2 (\theta^3)^2=0,
\label{pn2}
\ee
as can be easily calculated.  

For $n=5$, a \emph{point} on a curve $(1,x,x^2,x^3,x^4)$ in
$\bbR^5$ defines a
\emph{3-plane} passing through the origin in the dual space $(\bbR^5)^*$ 
via the relation:
\begin{eqnarray}
\theta^0 +4\theta^1 x+6\theta^2 x^2+4\theta^3 x^3+\theta^4 x^4&=&0\label{uw3}\\
\theta^1+3\theta^2 x+3\theta^3 x^2+\theta^4 x^3&=&0,\nonumber
\end{eqnarray}
where $(\theta^0,\theta^1,\theta^2,\theta^3,\theta^4)$ parametrize points of
the dual $(\bbR^5)^*$ as before. And now, when moving along the 
curve $(1,x,x^2,x^3,x^4)$
in $\bbR^5$, the
corresponding 3-planes in $(\bbR^4)^*$ sweep out a ruled
\emph{hypersurface} there, which is again defined by the vanishing 
of the resultant of the two polynomials defined in (\ref{uw3}). The
algebraic expression for this hypersurface in terms of the $\theta$
coordinates is quite complicated:
\begin{eqnarray}
&&-36 (\theta^1)^2 (\theta^2)^2 (\theta^3)^2 + 54 \theta^0 (\theta^2)^3 (\theta^3)^2 + 
   64 (\theta^1)^3 (\theta^3)^3 - 108 \theta^0 \theta^1 \theta^2 (\theta^3)^3 +\nonumber \\&&27 (\theta^0)^2 (\theta^3)^4 + 
   54 (\theta^1)^2 (\theta^2)^3 \theta^4 - 81 \theta^0 (\theta^2)^4 \theta^4 - 
   108 (\theta^1)^3 \theta^2 \theta^3 \theta^4 +\label{pn3}\\&& 180 \theta^0 \theta^1 (\theta^2)^2 \theta^3 \theta^4 + 
   6 \theta^0 (\theta^1)^2 (\theta^3)^2 \theta^4 - 54 (\theta^0)^2 \theta^2 (\theta^3)^2 \theta^4 + 
   27 (\theta^1)^4 (\theta^4)^2\nonumber\\&& - 54 \theta^0 (\theta^1)^2 \theta^2 (\theta^4)^2 + 18 (\theta^0)^2 (\theta^2)^2 (\theta^4)^2 + 12 (\theta^0)^2 \theta^1 \theta^3 (\theta^4)^2 - 
   (\theta^0)^3 (\theta^4)^3=0,\nonumber\end{eqnarray}
but easily calculable.  

The beauty of the hypersurfaces (\ref{pn1}), (\ref{pn2}) and (\ref{pn3})
consists in this that they are given by means of homogeneous
equations, and thus they descend to the corresponding projective spaces. From
the point of view of the present paper, even more
important is the fact, that they are $\Gl$ \emph{invariant}. By this we
mean the following:

Consider a real polynomial of $(n-1)$-th degree  
\be
w(x)=\sum_{i=0}^{n-1} \binom{n-1}{i} \theta^ix^i\label{pol}
\ee
in the real variable $x$ with real coefficients 
$( \theta^0,\theta^1,...,\theta^{n-1})$. The $n$-dimensional 
vector space $(\bbR^{n})^*$ of such polynomials may be
identified with the space of their coefficients. Now,  
replacing the variable $x$ by a new variable $x'$ such that 
\be
x=\frac{\al x'+\bet}{\gamma x'+\delta},\quad\quad
\al\delta-\bet\gamma\neq 0, \label{tra}
\ee
we define a new covector $(\theta'^0,\theta'^1,...,\theta'^{n-1})$ which is related
to $(\theta^0,\theta^1,...,\theta^{n-1})$ of (\ref{pol}) via
$$\sum_{i=0}^{n-1} \binom{n-1}{i} \theta'^ix'{^i}=(\gamma
x'+\delta)^{n-1}w(x).$$
It is obvious that $\theta'=(\theta'^0,\theta'^1,...,\theta'^{n-1})$ is  
linearly expressible in terms of $\theta=(\theta^0,\theta^1,...,\theta^{n-1})$:
\be
\theta'~=~\theta~\cdot~\rho_n(a),\quad\quad\quad a=\bma \al&\beta\\\gamma&\delta\ema.\label{glac}\ee
Here $a$ corresponds to the $\glg(2,\bbR)$ transformation (\ref{tra}), and 
the map $$\rho_n:\glg(2,\bbR)\to\glg((\bbR^{n})^*)\cong\glg(n,\bbR)$$
defines the \emph{real $n$-dimensional irreducible representation} of $\glg(2,\bbR)$.
For example, if $n=2$, we have $w(x)=\theta^0+2\theta^1 x +\theta^2x^2$, 
and we easily get 
$$\bma \theta'^0&\theta'^1&\theta'^2\ema=\bma \theta^0& \theta^1& \theta^2\ema\bma
\delta^2&\gamma\delta&\gamma^2\\
2\beta\delta &\al\delta+\beta\gamma&2\al\gamma\\
\beta^2&\al\beta&\al^2\ema,$$
so that $\rho_2$ is given by
$$\rho_2\bma\al&\beta\\\gamma&\delta\ema=\bma
\delta^2&\gamma\delta&\gamma^2\\
2\beta\delta &\al\delta+\beta\gamma&2\al\gamma\\
\beta^2&\al\beta&\al^2\ema.$$

Now, let us define $g(\theta,\theta)$,
${^4I}(\theta,\theta,\theta,\theta)$ and
${^5I}(\theta,\theta,\theta,\theta,\theta,\theta)$ by
\begin{eqnarray}
g(\theta,\theta)={\rm the~ left~ hand~ side ~of}~ (\ref{pn1})\nonumber\\
{^4I}(\theta,\theta,\theta,\theta)={\rm the~ left~ hand~ side ~of}~ (\ref{pn2})\label{mnmnn}\\
{^5I}(\theta,\theta,\theta,\theta,\theta,\theta)={\rm the~ left~ hand~ 
side ~of}~ (\ref{pn3}). \nonumber
\end{eqnarray}
We will often abbreviate this notation to the respective: 
$g(\theta)$, ${^4I}(\theta)$ and ${^5I}(\theta)$.

To explain our
comment about the $\Gl$ invariance of the respective 
hypersurfaces $g(\theta)=0$,
${^4I}(\theta)=0$ and
${^5I}(\theta)=0$ we calculate $g(\theta')$,
${^4I}(\theta')$ and
${^5I}(\theta')$ with $\theta'$ as in (\ref{glac}). The result is
\begin{eqnarray}
g(\theta')&=&(\al\delta-\beta\gamma)^2~g(\theta)\nonumber\\
{^4I}(\theta')&=&(\al\delta-\beta\gamma)^4~{^4I}(\theta)\nonumber\\
{^5I}(\theta')&=&
(\al\delta-\beta\gamma)^6~{^5I}(\theta).\nonumber\end{eqnarray}
Thus the vanishing of the expressions $g(\theta)$,
${^4I}(\theta)$ and
${^5I}(\theta)$ is invariant under the
action (\ref{glac}) of the irreducible $\Gl$ on $(\bbR^n)^*$.

We are now ready to discuss the general case $n\geq 3$ of the
\emph{rational normal curve} $(1,x,x^2,...,x^{n-1})$ in
$\bbR^n$. Associated with this curve is a pair of polynomials, namely 
$w(x)$ as in (\ref{pol}), and its \emph{derivative} $\frac{\der w}{\der
  x}$. We consider the 
relation  
\be
w(x)=0\quad\quad\&\quad\quad\frac{\der w}{\der x}=0.\label{uwn}\ee
This gives a \emph{correspondence} between the \emph{points} on the curve
$(1,x,x^2,...,x^{n-1})$ in $\bbR^n$ and the $(n-2)$-\emph{planes} passing
through the origin in the dual space $(\bbR^n)^*$ parametrized by $(\theta^0,\theta^1,...,\theta^{n-1})$. When moving along the rational normal curve in $\bbR^n$, the
corresponding $(n-2)$-planes in $(\bbR^n)^*$ sweep out a ruled
\emph{hypersurface} there. This is defined by the vanishing 
of the resultant, $R(w(x),\frac{\der w}{\der x})$, of the two
polynomials in (\ref{uwn}). The
algebraic expression for this hypersurface is the vanishing of 
a homogeneous polynomial, let us call it $I(\theta)$, 
of order $2(n-2)$, in the coordinates
$(\theta^0,\theta^1,...,\theta^{n-1})$. The hypersurface
$I(\theta)=0$ in $(\bbR^n)^*$ is $\Gl$ invariant, since the property of the
two polynomials $w(x)$ and $\frac{\der w}{\der x}$ to have a common
root is independent of the choice (\ref{tra}) of the coordinate $x$. Thus $\Gl$ \emph{is
included} in the stabilizer $G_I$ of $I$ under the action
of the full $\glg(n,\bbR)$ group. This stabilizer, by definition, 
is a subgroup of $\glg(n,\bbR)$
with elements $b\in G_I\subset\glg(n,\bbR)$ such that
$I(\theta\cdot b)=(\det b)^{\frac{2(n-2)}{n}}I(\theta)$. 
Moroever, in $n= 4,5$, it turns out that $G_I$ is \emph{precisely} the group $\Gl$ in the 
corresponding irreducible representation $\rho_n$. Thus if $n=4,5$
one can characterize the irreducible $\Gl$ in $n$ dimensions as the
stabilizer of the polynomial $I(\theta)$.

Crucial for the present paper is an observation of Karl
W\"unschmann
that the algebraic geometry and the correspondences we were describing
above, naturally appear in the analysis of solutions of
the ODE $y^{(n)}=0$. Indeed, following W\"unschmann\footnote{We are
  very grateful to Niels Schuman, who found a copy of W\"unschmann
  thesis in the
\emph{city} library of Berlin and sent it to us. It was this copy, which after
translation from German by Denson Hill, led us to write this introduction.}  (see the
Introduction in his PhD thesis \cite{wun}, pp. 5-6), we note the
following:

Consider the third order ODE: $y'''=0$.
Its general solution is 
$y=c_0+2c_1x+c_2x^2,$
where $c_0,c_1,c_2$ are the integration constants. Thus, the solution space 
of the ODE $y'''=0$ is $\bbR^3$ with solutions identified with
points ${\bf c}=(c_0,c_1,c_2)\in\bbR^3$. The solutions to the ODE $y'''=0$ 
may be also identified with curves $y(x)=c_0+2c_1x+c_2x^2$, 
actually parabolas, in the plane $(x,y)$. Suppose now, that we
take two solutions of $y'''=0$ corresponding to two
neighbouring points ${\bf c}=(c_0,c_1,c_2)$ and ${\bf c}+\der{\bf c}=(c_0+\der c_0,c_1+\der
c_1,c_2+\der c_2)$ in $\bbR^3$. Among all pairs of 
neighbouring solutions we choose only those, which have the property
that their corresponding curves $y=y(x)$ and $y+\der y=y(x)+\der y(x)$ 
touch each other, at some point $(x_0,y_0)$ in the plane
$(x,y)$. If we do not require anything more about the properties of
this incidence of the two curves, we say that solutions ${\bf c}$ and
${\bf c}+\der{\bf c}$ 
have \emph{zero order contact} at $(x_0,y_0)$. 

In this `baby' example everything is very simple:

To get the criterion 
for the solutions to
have zero order contact we first write the curves 
$y=c_0+2c_1x+c_2x^2$ and $y+\der y=c_0+\der c_0+2(c_1+\der
c_1)x+(c_2+\der c_2)x^2$ corresponding to ${\bf c}$ and ${\bf
  c}+\der{\bf c}$. 
Thus the solutions have zero order contact at
$(x_0,y(x_0))$ provided that $\der y(x_0)=0$, i.e. if and only if $$\der c_0+2x_0\der
c_1+x_0^2\der c_2=0.$$
This shows that such a contact is possible if and only if the
determinant  
$$
g(\der {\bf c},\der{\bf c})=(\der c_1)^2-\der c_0\der c_2$$
is \emph{nonnegative}, since otherwise the quadratic equation for
$x_0$ has no solutions. Unexpectedly, we find that the 
requirement for the two neighbouring solution curves of $y'''=0$ 
to have zero order contact at some point is equivalent to the requirement 
that the corresponding two neighbouring points ${\bf c}$ and ${\bf
  c}+\der{\bf c}$ 
in $\bbR^3$ are
\emph{spacelike separated} in the \emph{Minkowski metric} $g$ in
$\bbR^3$. This is the discovery of W\"unschmann that is quoted in Elie
Cartan's 1941 year's paper\footnote{It is worthwhile to remark, that W\"unschmann thesis is dated 
`1905',  the same year in which Einstein published his famous special 
relativity theory paper \cite{ein}. It was not until three years later 
when Minkowski gave the geometric
interpretation of Einstein's theory in terms of his metric \cite{min}. Perhaps
W\"unschmann was the first who ever wrote such metric in a scientific paper. 
This is a very interesting feature of
W\"unschmann thesis: he calls the expressions like $(\der c_1)^2-\der c_0\der
c_2=0$, a \emph{Mongesche Gleichung} rather than a \emph{cone in the
metric}, because the notion of a metric with signature different than
the Riemannian one was not yet abstracted!} \cite{carspan}.

Now we consider 
the neighbouring solutions ${\bf c}$ and ${\bf c}+\der{\bf c}$ of $y'''=0$ which
are \emph{null} separated in the metric $\der s^2$. What we can say
about the corresponding curves in the plane $(x,y)$? 

To answer this we need the notion of a \emph{first order contact}: Two
neighbouring solution curves $y=c_0+2c_1x+c_2x^2$ and $y+\der y=c_0+2c_1x+c_2x^2+(\der c_0+2x\der
c_1+x^2\der c_2)$ of $y'''=0$, corresponding to ${\bf c}$ and ${\bf
  c}+\der{\bf c}$ 
in $\bbR^3$, have first order contact at $(x_0,y_0)$
iff they have zero order contact at $(x_0,y_0)$ and, in addition,
their curves of \emph{first
derivatives}, $y'=2 c_1+2 c_2 x$ and $y'+\der y'=2 (c_1+\der c_1)+2
(c_2+\der c_2) x$, have zero order contact at $(x_0,y_0)$. Thus
the condition of first order contact at $(x_0,y(x_0))$ is
equivalent to
$\der y(x_0)=0$ and $\der y'(x_0)=0$, 
i.e. to the condition that $x_0$ is a \emph{simultaneous} root for  
\begin{eqnarray}
\der c_0+2x_0\der
c_1+x_0^2\der c_2&=&0\label{3ci}\\
\der c_1+x_0\der c_2&=&0.\nonumber
\end{eqnarray}
Solving the second of these equations for $x_0$, and inserting it into the
first, after an obvious simplification, we conclude that $(\der c_1)^2-\der
c_0\der c_2=0$. Thus we get the interpretation of the
\emph{null separated} neighbouring points in $\bbR^3$ as the solutions
of $y'''=0$ whose curves in the $(x,y)$ plane are neighbouring and have first order
contact at some point.

W\"unschmann notes that the procedure described here for the 
equation $y'''=0$ can be repeated for the equation $y^{(n)}=0$ for
arbitrary $n\geq 3$. In the cases of $n=4$ and $n=5$ he however passes
to the discussion of the solutions that have contact of order
$(n-2)$ rather then one. This is an interesting possibility, complementary in a sense to
the one in which the solutions have first order contact. W\"unschmann
spends rest of the thesis studying it. But we will
not discuss it here. 

Since W\"unschmann does not discuss the first order contact of the
solutions in $n=4,5$, let us look closer onto these two cases:

The general solution to $y^{(4)}=0$ is $y=c_0+3 c_1 x+3c_2 x^2+c_3
x^3$, and the general solution to $y^{(5)}=0$ is $y=c_0+4 c_1 x+6c_2
x^2+4c_3 x^3+c_4 x^4$. Thus now the solutions are points ${\bf c}$ in $\bbR^4$
and $\bbR^5$, respectively. The condition that the neighbouring
solutions ${\bf c}=(c_0,c_1,c_2,c_3)$ and 
${\bf c}+\der{\bf c}=(c_0+\der c_0,c_1+\der c_1,c_2+\der
c_2,c_3+\der c_3)$ of $y^{(4)}=0$ have first order contact at $(x_0,y(x_0))$ is
equivalent to the requirement that the system    
\begin{eqnarray}
\der c_0+3x_0\der c_1+3x_0^2\der c_2+x_0^3\der c_3&=&0\label{4ty}\\
\der c_1+2x_0\der c_2+x_0^2\der c_3&=&0\nonumber
\end{eqnarray}
has a common root $x_0$. Similarly, the condition 
that the neighbouring 
solutions ${\bf c}=(c_0,c_1,c_2,c_3,c_4)$ and 
${\bf c}+\der{\bf c}=(c_0+\der c_0,c_1+\der c_1,c_2+\der
c_2,c_3+\der c_3,c_4+\der c_4)$ of $y^{(5)}=0$ have first order contact at $(x_0,y(x_0))$ is
equivalent to the requirement that the system   
\begin{eqnarray}
\der c_0+4x_0\der c_1+6x_0^2\der c_2+4x_0^3\der c_3+x_0^4\der c_4&=&0\label{5ty}\\
\der c_1+3x_0\der c_2+3x_0^2\der c_3+x_0^3 \der c_4&=&0\nonumber
\end{eqnarray}
has a common root $x_0$.
Calculating the resultants 
for the systems (\ref{3ci}), (\ref{4ty}), and (\ref{5ty}) we get:
\begin{itemize}
\item $R_3=g(\der{\bf c},\der{\bf c})\der c_2$\hspace{2.4cm} if $n=3$,
\item $R_4={^4I}(\der{\bf c},\der{\bf c},\der{\bf
  c},\der{\bf c})\der c_3$ \hspace{1.15cm} if $n=4$,
\item $R_5={^5I}(\der{\bf c},\der{\bf c},\der{\bf
  c},\der{\bf c},\der{\bf c},\der{\bf c})\der c_4$ \hspace{.1cm} if $n=5$,
\end{itemize}
where $g$, ${^4I}$ and ${^5I}$ are as in (\ref{mnmnn}).

This confirms our earlier statement that two neighbouring solutions of
   $y'''=0$ have first order contact iff $g(\der{\bf c},\der{\bf
  c})=0$, since if $\der c_2=0$ the system (\ref{3ci}) collapses to 
$\der c_1=\der c_0=0$. Similarly, one can prove that two
  neighbouring solutions of $y^{(4)}=0$ or $y^{(5)}=0$ have first
  order contact if and only if they are \emph{null separated} in the respective
symmetric multilinear forms ${^4I}$ or ${^5I}$. Our previous discussion of 
the invariant properties of these forms, shows that in the solution
   space of an ODE $y^{(n)}=0$, for $n\geq 4$, there is a naturally
   defined action of the the $\Gl$ group. This group is the stabilizer
   of the invariant polynomial $I(\der{\bf c})$ which distinguishes neighbouring
   solutions having first order contact.  

The main question one can ask in this context is if one can retain
this $\Gl$ structure in the solution space for more complicated
ODEs. In other words, one may asks the following: What does one
have to assume about the function $F$, defining an ODE 
$$y^{(n)}=F(x,y,y',\dots,y^{(n-1)}),$$ 
in order \emph{to have} a well defined \emph{conformal} tensor
$g$, 
$^4I$ or $^5I$, in the
respective cases $n=3,4,5$, on the solution space of the ODE? The same
question can be repeated for any $n>5$ and the invariant $I$.

The answer to this question in the $n=4$ case was found by Robert Bryant
in \cite{bryantsp4}. Later, the answer for $n>4$ case was given by
Boris Doubrov \cite{dub} who established a connection between 
the \emph{Wilczynski invariants} \cite{wilk} for a linear ODE, and 
certain \emph{contact invariant
conditions} for a \emph{nonlinear} ODE associated with it. 
For any $n\geq 3$, given $F$, Doubrov conditions are built from 
the Wilczynski invariants
calculated for the linearization of $y^{(n)}=F$ about one of its
solutions (see \cite{dub} for details). In a quite different
perspective, these conditions, 
were also discovered by Maciej Dunajski and Paul K Tod \cite{dun}.

Doubrov-Wilczynski conditions differ from Bryant ones for $n=4$. They
also differ from the conditions we are going to discuss in the present paper for $n\geq 4$. Doubrov, Dunajski and Tod have $(n-2)$ nonlinear PDEs for $F$ of ODE $y^{(n)}=F$. Although this
number, $(n-2)$,  agrees with the number of conditions we present here, there is an important difference: each of the 
$(n-2)$ conditions for $F$, defined by the above authors, 
has a \emph{different} 
order in the derivatives of $F$. When we arrange Doubrov-Wilczynski 
conditions according to
the order of the corresponding PDEs for $F$, 
we find that the first
condition is of order $3$, the second is of order $4$, and so on, up to
the order $n$ of the $(n-2)$-th condition. On the contrary \emph{each} of \emph{our}
$(n-2)$ conditions is of the \emph{third} order in the derivatives of $F$. 
The simple explanation of this
discrepancy is as follows: We obtain our conditions, 
by applying a variant of \emph{Cartan's equivalence method}; in the process of extracting
them we obtain the first condition to be of the third order as everybody
does. But the second condition which, if we were not applying Cartan's
method, would be of order four, actually collapses in our derivation 
to a condition of order \emph{three}. This is because Cartan's
method automatically utilises the first condition of order three 
by differentiating it, and then eliminating the
fourth derivative from the second condition by means of the fourth
derivative from the differentiated condition of order three. The same
situation is automatically accomplished for the condition of order
\emph{five} and so on. As a result we have $(n-2)$ conditions of order
\emph{three}. They are different from the Doubrov-Wilczynski conditions
already for the ODE of order \emph{four}. In the $n=3$ case all the
conditions, namely those of W\"unschmann, Doubrov, Dunajski and Tod,
and ours are the same. In dimension $n=4$ our conditions agree with
the Bryant ones. Since W\"unschman was the first who obtained these type of conditions in $n=3$ 
and found method of their constructing for arbitrary $n$ we call the
conditions discussed in this paper {\em generalized W\"unschmann's
  conditions}, or \emph{W\"unschmann's conditions}, for short.

Finding the W\"unschmann conditions for $F$ of order $n\geq4$,
although important, is only a byproduct of our analysis. The present
paper is devoted to a thorough study of the \emph{irreducible} $\Gl$
\emph{geometry in dimension five}. This is done from two points of
view: first as a study of an \emph{abstract} geometry on a manifold
and, second, as a study of a \emph{contact geometry of fifth order ODEs}.

We define an abstract 5-dimensional $\Gl$ geometry in two steps.

First, in section \ref{s1}, we show how to construct the algebraic model for
the $\Gl$ geometry in dimension five utilising properties of a \emph{rational normal
  curve}. Second, instead of obtaining the reduction from $\glg(5,\bbR)$
to $\Gl$ by stabilizing the 6-tensor $^5I$, we get the desired
reduction by stabilizing 
a conformal metric $g_{ij}\to e^{2\phi}g_{ij}$ of signature $(3,2)$
and a conformal totally symmetric 3-tensor $\ten_{ijk}\to e^{3\phi}
\ten_{ijk}$. These tensors are supposed to be 
related by the following algebraic constraint:
\be\label{intro1}g^{lm}(\ten_{ijl}\ten_{kmp}+\ten_{kil}\ten_{jmp}+\ten_{jkl}\ten_{imp})=g_{ij}g_{kp}+g_{kl}g_{jp}+g_{jk}g_{ip}.\ee
It is worthwile to note that condition (\ref{intro1}) is a non-Riemannian 
counterpart of the condition considered by Elie 
Cartan in the context of isoparametric surfaces \cite{Cartan1}, \cite{Cartan}.

Our main object of study is then defined as follows:
\begin{defintro}
An irreducible $\glg(2,\bbR)$ geometry in dimension five is a 5-dimensional
manifold $M^5$ equipped with a class of triples $[g,\ten, A]$ such
that on $M^5$:
\begin{itemize}
\item[$(a)$] $g$ is a
metric of signature $(3,2)$,
\item[$(b)$] $\ten$ is a traceless symmetric 3rd rank tensor,
\item[$(c)$] $A$ is a 1-form,
\item[$(d)$] the metric $g$ and the tensor $\ten$ satisfy the identity \eqref{intro1},
\item[$(e)$] two triples $(g,\ten,A)$ and $(g',\ten',A')$ are in the same class
  $[g,\ten,A]$ if and only if 
  there exists a function $\phi:M^5\to \bbR$ such that
$$g'={\rm e}^{2\phi}g,\quad\quad\quad \ten'={\rm
    e}^{3\phi}\ten,\quad\quad\quad A'=A-2\der\phi.$$
\end{itemize}
\end{defintro}

This definition places $\Gl$ geometries in dimension five among the \emph{Weyl
geometries} $[g,A]$. They are special Weyl geometries i.e. such for which the structure
group is reduced from $\cog(3,2)$ to $\Gl$. A natural description of
such geometries should be then obtained in terms of a certain
$\gla(2,\bbR)$-valued connection. However, unlike the usual Weyl
case, the
choice of such a connection is ambiguous, due to the fact that
such a connection has non-vanishing torsion in general, and one must
find admissible conditions for the torsion that specify connection
uniquely. Pursuing this problem in section \ref{trzy} we find an interesting 
subclass of $\Gl$ geometries.
\begin{propintro}
A $\Gl$ geometry $[g,\ten,A]$ is called nearly integrable 
if the Weyl connection $\stackrel{W}{\nabla}$ of $[g,A]$ satisfies
$$(\stackrel{W}{\nabla}_X\ten)(X,X,X)=-\tfrac12 A(X)\ten(X,X,X).$$
It turns out, see section \ref{trzy}, that the nearly integrable $\Gl$
geometries \emph{uniquely} define a $\gl$ connection $D$. This is
characterized by the following requirements: 
\begin{itemize}
\item it preserves the structural tensors:
\begin{eqnarray*}
&&Dg_{ij}=-A g_{ij},\\
&&D\ten_{ijk}=-\tfrac32 A\ten_{ijk},
\end{eqnarray*}
 \item and it {\em has totally skew symmetric torsion}.
\end{itemize}
We call this unique connection the {\em characteristic} connection for
the nearly integrable structure $\Gl$.
\end{propintro}

The rest of section \ref{trzy} is devoted to study the 
algebraic structure of the torsion and the 
curvature of the characteristic connection of a nearly integrable
structure. 
Since the tensor products of tangent spaces are reducible under the 
action of $\Gl$, we decompose the torsion and the curvature tensors
into components belonging to the irreducible representations. 
In particular, the skew symmetric torsion $T$ has two components,
$T^{(3)}$ and $T^{(7)}$, lying in the three-dimensional and the
seven-dimensional irreducible representations respectively. Likewise
the Maxwell 2-form $\der A$ and the Ricci tensor $R$ decompose
according to $\der A=\der A^{(3)}+\der A^{(7)}$ and 
$R=R^{(1)}+R^{(3)}+R^{(5)}+R^{(7)}+R^{(9)}$. The last problem we
adress in section \ref{trzy} concerns with the properties of geometries whose
characteristic connections have `the smallest possible' torsion,
that is the torsion of the pure three-dimensional type. In proposition
\ref{p.purtor} we prove that Ricci tensor for such structures
satisfies remarkable identities:
$$ R^{(3)}=\tfrac14\der A^{(3)},\qquad\qquad R^{(7)}=\tfrac32\der A^{(7)}, \qquad\qquad R^{(9)}=0.$$
The third equation is equivalent to 
$$ R_{(ij)}=\tfrac{1}{5}R g_{ij}+\tfrac{2}{7}R_{kl}\ten^{klm}\ten_{ijm}.$$

In section \ref{nabu} we briefly describe $\Gl$ geometry in the language of the bundle $\Gl\to P\to M^5$. 
We also show how an appriopriate coframe defined on a nine-dimensional
manifold $P$ turns this manifold into a bundle $\Gl\to P\to M^5$ and
generates the $\Gl$ geometry on $M^5$. This construction is the core of the proof of the main theorem in section \ref{secsec}. This closes the part of the paper that is devoted to abstract $\Gl$ geometries.

Section \ref{secsec} contains the main result of this paper, theorem \ref{maint}, 
which links $\Gl$ geometry with the realm of ordinary differential equations. 
It can be encapsulated as follows.
\begin{thintro}
A 5th order ODE $y^{(5)}=F(x,y,y',y'',y''',y^{(4)})$ 
that satisfies three W\"unschmann conditions defines a nearly integrable irreducible $\glg(2,\bbR)$ geometry
$(M^5,$ $[g,\ten,A])$ on the space $M^5$ of its solutions. 
This geometry has the characteristic connection with torsion 
of the `pure' type in the 3-dimensional irreducible representation of $\Gl$.
Two 5th order ODEs that are equivalent under contact transformation of variables define equivalent $\Gl$ geometries. 
\end{thintro}
The theorem has numerous applications. For example, we use it to 
characterise various classes of W\"unschmann 5th order ODEs, 
by means of the algebraic type of the tensors associated with the
corresponding characteristic connection. For example 
iff $F_{y^{(4)}y^{(4)}}=0$, the torsion of the characterstic connection
vanishes, and iff $F_{y^{(4)}y^{(4)}y^{(4)}}=0$, then we have $\der A^{(7)}=0$.

The proof of the theorem consists of an application of the Cartan method of equivalence. 
We write an ODE, considered modulo contact transformation of
variables, as a $G$-structure on the four-order jet space. Starting
from this $G$-structure we explicitly construct a 9-dimensional
manifold $P$, which is a $\Gl$ bundle over the solution space and
carries certain distinguished coframe. This construction is only
possible provided that the ODE satisfies the W\"unschmann conditions,
which we write down explicitly. By means
of proposition \ref{p.constr} the coframe on $P$ defines the nearly
integrable geometry on the solution space of the ODE. It has the 
characteristic connection with torsion in the 3-dimensional representation. 

Section \ref{s.examples} includes examples of 5th order equations in
the W\"unschmann class. We find equations generating all the
structures with vanishing torsion, equations possessing at least
6-dimensional group of contact symmetries and yielding geometries with
$\der A=0$. We also give highly nontrivial 
examples of equations for which $\der A\neq 0$, including a 
family of examples with function $F$ being a solution of a 
certain second order ODE.

Finally, in section \ref{sue} we consider ODEs of order $n>5$. We
apply results of the Hilbert theory of algebraic invariants, to define
the tensors responsible for the reduction
$\glg(n,\bbR)\to\glg(2,\bbR)$. We also 
give the explicit formulae for the $(n-2)$ third order 
W\"unschmann conditions for $n=6$ and $n=7$.

\section{A peculiar third rank symmetric tensor}
\label{s1}
Consider $\bbR^n$ equipped with a \emph{Riemannian} metric $g$ and a
3rd rank tracefree symmetric tensor $\ten\in S^3_0\bbR^n$ satisfying:
\begin{itemize}
\item[(i)] $\ten_{ijk}=\ten_{(ijk)}$ - (symmetry)
\item[(ii)] $g^{ij}\ten_{ijk}=0$ - (tracefree)
\item[(iii)] $g^{lm}(\ten_{ijl}\ten_{kmp}+\ten_{kil}\ten_{jmp}+\ten_{jkl}\ten_{imp})=g_{ij}g_{kp}+g_{kl}g_{jp}+g_{jk}g_{ip}$.        
\end{itemize} 
It turns out that the third condition is very restrictive. In
particular Cartan has shown \cite{Cartan1,Cartan} that for (iii)
to be satisfied the dimension $n$ must be one of the following: 
$n=5,8,14,26$. Moreover Cartan constructed $\ten$ in each of these
dimensions and has shown that it is unique up to an $\og(n)$ transformation. 
Restricting to $n=5,8,14,26$, we
consider the stabilizer $H_n$ of $\ten$ under the action of $\glg(n,\bbR)$:
$$H_n=\{\glg(n,\bbR)\ni a: \ten(aX,aY,aZ)=\ten(X,Y,Z), ~~
\forall X,Y,Z\in\bbR^n\}.$$
Then, one finds that: 
\begin{itemize}
\item $H_5=\sog(3)\subset\sog(5)$ in the 5-dimensional irreducible
representation,
\item $H_8=\sug(3)\subset\sog(8)$ in the 8-dimensional irreducible
representation,
\item $H_{14}=\spg(3)\subset\sog(14)$ in the 14-dimensional irreducible
representation,
\item $H_{26}={\bf F}_4\subset\sog(26)$ in the 26-dimensional irreducible
representation.
\end{itemize}
The relevance of conditions (i)-(iii) is that they are invariant under
the $\og(n)$ action on the space of tracefree symmetric tensors
$S^3_0\bbR^n$. Moreover they totally characterize the orbit 
$\og(n)/H_n\subset S^3_0$ of the tensor $\ten$ under this action \cite{bobi,iso}.\\

The question arises if one can construct tensors satisfying (i)-(iii)
for metrics having non-Riemannian signatures. Below we show how to do
it if $n=5$ and the metric $g$ has signature $(3,2)$. This construction
described to us by E. Ferapontov \cite{fera,fer1} is as follows. 

Consider $\bbR^5$ with coordinates $(\theta^0,\theta^1,\theta^2,\theta^3,\theta^4)$, and a curve
$$\gamma(x)=(1,x,x^2,x^3,x^4)\subset\bbR^5.$$
Associated to the curve $\gamma$ there are \emph{two} algebraic varieties in
$\bbR^5$:
\begin{itemize}
\item \emph{The bisecant variety}. This is defined to be a set consisting
  of all the points on all straight lines \emph{crossing the curve} $\gamma$
	in precisely \emph{two} points. It is given parametrically as
$$B(x,s,u)=(1,x,x^2,x^3,x^4)+u(0,x-s,x^2-s^2,x^3-s^3,x^4-s^4),$$
where $x,s,u$ are \emph{three} real parameters.
\item \emph{The tangent variety}. This is defined to be a set consisting
  of all the points on all straight lines \emph{tangent} to 
the curve $\gamma$. It is given parametrically as
$$T(x,s)=(1,x,x^2,x^3,x^4)+s(0,1,2x,3x^2,4x^3).$$
\end{itemize} 

One of many interesting features of these two varieties is that 
they define (up to a scale) a tri-linear symmetric form
\be\label{ten1} \ten(\theta)=3\sqrt{3}(\theta^0\theta^2\theta^4+2\theta^1\theta^2\theta^3-(\theta^2)^3-\theta^0(\theta^3)^2-\theta^4(\theta^1)^2) \ee 
and a bi-linear symmetric form
\be \label{met1} g(\theta)=\theta^0\theta^4-4\theta^1\theta^3+3(\theta^2)^2.\ee
These forms are distinguished by the fact that the bisecant and tangent varieties are contained in their null cones.
In the homogeneous coordinates 
$(\theta^0,\theta^1,\theta^2,\theta^3,\theta^4)$ in $\bbR^5$ all the points $\theta$ of $B(x,s,u)$ satisfy
$$\ten(\theta)=0,$$
whereas all the points $\theta$ of $T(x,s)$ satisfy
$$\ten(\theta)=0\quad{\rm and}\quad g(\theta)=0.$$

Writing the forms as $\ten(\theta)=\ten_{ijk}\theta^i\theta^j\theta^k$,
$g(\theta)=g_{ij}\theta^i\theta^j$, $i,j,k=0,1,2,3,4$ one can check that so
defined $g_{ij}$ and $\ten_{ijk}$ satisfy relations (i)-(iii) of the
previous section.

Although it is obvious we remark that the above defined metric
$g_{ij}$ has signature $(3,2)$. 

As we have already noted the forms $\ten(\theta)$ and $g(\theta)$ are defined only \emph{up to a scale}. 
We were also able to find a factor, the $3\sqrt{3}$ in
expression (\ref{ten1}), that makes the corresponding $g_{ij}$ and
$\ten_{ijk}$ satisfy (i)-(iii). Note that these conditions are 
\emph{conformal} under the simultaneous change:
$$g_{ij}\to {\rm e}^{2\phi}g_{ij},\quad\quad\ten_{ijk}\to {\rm
  e}^{3\phi}\ten_{ijk}.$$

Thus it is interesting to consider in $\bbR^5$ a 
\emph{class} of pairs $[g,\ten ]$, such that: 
\begin{itemize}
\item in each pair $(g,\ten)$
\begin{itemize}
\item $g$ is a
metric of signature $(3,2)$,
\item $\ten$ is a traceless symmetric 3rd rank tensor, 
\item the metric $g$ and the tensor $\ten$ satisfy the identity
$$g^{lm}(\ten_{ijl}\ten_{kmp}+\ten_{kil}\ten_{jmp}+\ten_{jkl}\ten_{imp})=g_{ij}g_{kp}+g_{kl}g_{jp}+g_{jk}g_{ip},$$
\end{itemize}
\item two pairs $(g,\ten)$ and $(g',\ten')$ are in the same class
  $[g,\ten]$ if and only if 
  there exists $\phi\in\bbR$ such that
\be
g'={\rm e}^{2\phi}g,\quad\quad\quad \ten'={\rm e}^{3\phi}\ten.\label{cotr}
\ee
\end{itemize}

Given a structure $(\bbR^5,[g,\ten])$ we define a group $CH$ to be a
subgroup of the general linear group
$\glg(5,\bbR)$ preserving $[\ten]$. This means that, choosing a
representative $\ten$ of the class $[\ten]$, we define $CH$ to be:
$$CH=\{\glg(5,\bbR)\ni a: \ten(ax,ax,ax)=(\det
a)^{(3/5)}\ten(x,x,x)\}.$$
Note that the exponent $\tfrac35$ in the above expression is caused by
the fact that the r.h.s. of the equation defining the group elements must be
\emph{homogeneous of degree} 3 in $a$, similarly as the l.h.s. is. 

This definition does not depend on the
choice of a representative $\ten\in [\ten]$. We have the following
\begin{proposition}
The set $CH$ of $5\times 5$ real matrices $a\in\glg(5,\bbR)$
preserving $[\ten]$ is the $\glg(2,\bbR)$ group
in its 5-dimensional irreducible representation. Moreover, we have
natural inclusions
$$CH=\glg(2,\bbR)\subset\cog(3,2)\subset\glg(5,\bbR),$$ 
where 
$\cog(3,2)$ is the $11$-dimensional group of homotheties 
associated with the conformal class $[g]$.
\end{proposition}
\begin{remark}
According to our Introduction, there is another $\Gl$ invariant
symmetric conformal tensor that stabilizes $\glg(5,\bbR)$ to the
irreducible $\Gl$. This is the tensor ${^5I}_{ijklpq}$ defined via 
${^5I}(\theta)=\tfrac{1}{120}{^5I}_{ijklpq}\theta^i\theta^j\theta^k\theta^l\theta^p\theta^q$
with ${^5I}$ as in (\ref{mnmnn}). We prefer however to work with a pair
$(g_{ij},\ten_{ijk})$ rather then with ${^5I}_{ijklpq}$, because of
the lower rank, and more importantly, because of the evident conformal
\emph{metric} properties of the $(g_{ij},\ten_{ijk})$ approach. Also, it is
worthwhile to note that the invariants $g_{ij}$, $\ten_{ijk}$ and
${^5I}_{ijklpq}$ are \emph{not} independent. Indeed, one can easilly
check that ${^5I}$ of (\ref{mnmnn}), $\ten$ of (\ref{ten1}) and
$g$ of (\ref{met1}) are related by ${^5I}=\ten^2-g^3$. We interpret 
this relation as the definition of $^5I$ in terms of more primitive
quantities $g$ and $\ten$. 
\end{remark}
The isotropy condition for
the group elements $a$ of $CH$ has its obvious counterpart at the
level of the Lie algebra $\gla(2,\bbR)=(\bbR\oplus\sla(2,\bbR))\subset\coa(3,2)\subset\gla(5,\bbR)$ of $CH=\glg(2,\bbR)$. 
Writing $a={\rm exp}(t\Gamma)$ we find that the infinitesimal version of
the isotropy condition, written in terms of the $5\times 5$ matrices 
$\Gamma=(\Gamma^i_{~j})$ is:
\be
\Gamma^l_{~i}\ten_{ljk}+\Gamma^l_{~j}\ten_{ilk}+\Gamma^l_{~k}\ten_{ijl}=\tfrac35 Tr(\Gamma)\ten_{ijk},
\label{ica}
\ee 
where $Tr(\Gamma)=\Gamma^m_{~m}$. Given $\ten_{ijk}$, these 
\emph{linear} equations can be solved for $\Gamma$. Taking the most
general matrix $\Gamma\in\glg(5,\bbR)$ and $\ten_{ijk}$ given by 
$\ten(x,x,x)=\ten_{ijk}x^ix^jx^k$ of (\ref{ten1}) we find the explicit
realization of the 5-dimensional representation of the $\gla(2,\bbR)$
Lie algebra as:
\be
\Gamma=\Om_- E_-+\Om_+E_++\Om_0 E_0+\Om_1 E_1,\label{om1}
\ee
where  $\Om_-,\Om_+,\Om_0,\Om_1$ are free real parameters, 
and $(E_-,E_+,E_0,E_1)$ are $5\times 5$ matrices given by:
\begin{eqnarray}
&&E_+=\bma
0&4&0&0&0\\0&0&3&0&0\\0&0&0&2&0\\0&0&0&0&1\\0&0&0&0&0\ema,\quad
E_-=\bma
0&0&0&0&0\\1&0&0&0&0\\0&2&0&0&0\\0&0&3&0&0\\0&0&0&4&0\ema,\label{mo21}\\
&&E_0=\bma
-4&0&0&0&0\\0&-2&0&0&0\\0&0&0&0&0\\0&0&0&2&0\\0&0&0&0&4\ema,\quad E_1=-4\bma
1&0&0&0&0\\0&1&0&0&0\\0&0&1&0&0\\0&0&0&1&0\\0&0&0&0&1\ema.\nonumber
\end{eqnarray}

The commutator in
$$\gla(2,\bbR)={\rm Span}_\bbR(E_-,E_+,E_0,E_1)$$ is 
the usual commutator of matrices. In 
particular, the non-vanishing commutators are:
$$[E_0,E_+]=-2E_+\quad,\quad [E_0,E_-]=2E_-\quad,\quad
[E_+,E_-]=-E_0.\quad$$
Note that 
$$\sla(2,\bbR)={\rm Span}_\bbR(E_-,E_+,E_0)$$
is a subalgebra of $\gla(2,\bbR)$ isomorphic to $\sla(2,\bbR)$. It
provides the 
$5$-dimensional irreducible representation of $\sla(2,\bbR)$.

\section{Irreducible $\glg(2,\bbR)$ geometries in dimension five}\label{trzy}
In this section we describe 5-dimensional manifolds whose tangent
space at each point is equipped with the structure $[g,\ten]$ of the
previous section. We will analyze such manifolds in terms of an
appropriately chosen connection. We will describe connections on a
manifold $M$ in terms of Lie-algebra-valued 1-forms on $M$. To be more
specific, let ${\rm dim}M=n$ and let $\mathfrak{g}$ denote an 
$n$-dimensional representation of some Lie algebra. The connection
1-forms $\Gamma^i_{~j}$ on $M$ are the matrix entries of an element 
$\Gamma\in {\mathfrak g}\otimes \Lambda^1 M$. They 
define the covariant exterior derivative $D$. This acts on
tensor-valued-forms via the 
extension to the higher order tensors of the formula:
$$Dv^i=\der v^i+\Gamma^i_{~j}\dz v^i.$$

Now suppose that we have a 5-dimensional manifold $M^5$ equipped with
a class of pairs $[g,\ten]$ such that $g$ is a metric, $\ten$ is a 3rd
rank tensor related to the metric via properties (i)-(iii) of the
previous section, and two pairs $(g,\ten)$ and $(g',\ten')$ are in the
same pair iff they are related by (\ref{cotr}), where $\phi$ is now a
\emph{function} on $M^5$. If we want to
associate a connection with such a structure we have to specify how this
connection is related to the pair $[g,\ten]$. A possible approach is
to choose a representative $(g,\ten)$ of $[g,\ten]$ and declare what is 
$Dg$ and $D\ten$. A first possible choice $Dg=0$ or $D\ten=0$ is
definitely not good since, in general, $Dg'$ and $D\ten'$ 
would not be vanishing for another choice of the representative of
$[g,\ten]$. A remedy for this situation comes from the \emph{Weyl
  geometry} where, given a conformal class $(M,[g])$, a 1-form $A$ is
introduced so that the connection satisfies $Dg_{ij}=-Ag_{ij}$. In our
case we introduce a 1-form $A$ on
$M^5$ and require that $Dg_{ij}=-Ag_{ij}$ and 
$D\ten_{ijk}=-\tfrac32 A\ten_{ijk}$. Then, if we transform $(g,\ten)$
according to (\ref{cotr}), the transformed objects will satisfy 
$Dg'_{ij}=-A'g'_{ij}$ and $D\ten'_{ijk}=-\tfrac32 A'\ten'_{ijk}$
provided that $A'=A-2\der\phi$. This motivates the following    
 
\begin{definition}\label{def1}
An irreducible $\glg(2,\bbR)$ structure in dimension five is a 5-dimensional
manifold $M^5$ equipped with a class of triples $[g,\ten, A]$ such
that on $M^5$:
\begin{itemize}
\item[$(a)$] $g$ is a
metric of signature $(3,2)$,
\item[$(b)$] $\ten$ is a traceless symmetric 3rd rank tensor,
\item[$(c)$] $A$ is a 1-form,
\item[$(d)$] the metric $g$ and the tensor $\ten$ satisfy the identity
$$g^{lm}(\ten_{ijl}\ten_{kmp}+\ten_{kil}\ten_{jmp}+\ten_{jkl}\ten_{imp})=g_{ij}g_{kp}+g_{kl}g_{jp}+g_{jk}g_{ip},$$
\item[$(e)$] two triples $(g,\ten,A)$ and $(g',\ten',A')$ are in the same class
  $[g,\ten,A]$ if and only if 
  there exists a function $\phi:M^5\to \bbR$ such that
$$g'={\rm e}^{2\phi}g,\quad\quad\quad \ten'={\rm
    e}^{3\phi}\ten,\quad\quad\quad A'=A-2\der\phi.$$
\end{itemize}
\end{definition}

If $M^5$ was only equipped with a class of \emph{pairs} $[g,A]$ satisfying
conditions $(a)$, $(c)$ and $(e)$ (with $\ten$, $\ten'$ omitted),
then $(M^5,[g,A])$ would define the \emph{Weyl} geometry. 
Such geometry, which has the
structure group $\cog(3,2)$, is usually studied in
terms of the \emph{Weyl connection}. This is the \emph{unique} torsionfree
connection preserving the conformal structure $[g,A]$. It is defined
by the following two equations:
\begin{eqnarray}
&&\stackrel{W}{D}g_{ij}=-Ag_{ij}\quad\quad{\rm (preservation~ of~
    the~ class}~[g,A]),\label{wc02}\\
&&\stackrel{W}{D}\theta^i=0\quad\quad{\rm (no~torsion),}\label{wc01}
\end{eqnarray}
where $\theta^i$ is a coframe related to the representative $g$ of the
class $[g]$ by $g=g_{ij}\theta^i\theta^j$.
We describe the Weyl connection 
in terms of the Weyl connection 1-forms $\we^i_{~j}$,
$i,j=0,1,2,3,4$.  

Take a representative $(g,A)$ of the Weyl structure $[g,A]$ on
$M^5$. Choose a coframe $(\theta^i)$, $i=0,1,2,3,4$, such that
$g=g_{ij}\theta^i\theta^j$, with all the metric coefficients $g_{ij}$
being \emph{constant}. Then the above two equations define 
$\we^i_{~j}$ together with $\we_{ij}=g_{ik}\we^k_{~j}$ 
to be 1-forms on $M^5$ satisfying
\begin{eqnarray}
&&\we_{ij}+\we_{ji}=Ag_{ij}\quad\quad{\rm (preservation~ of~
    the~ class}~[g,A]),\label{wc2}\\
&&\der\theta^i+\we^i_{~j}\dz\theta^j=0\quad\quad{\rm (no~torsion).}\label{wc1}
\end{eqnarray}
It follows that once the representative $(g,A)$ and the coframe
$\theta^i$ is chosen the above equations \emph{uniquely} determine the Weyl
connection 1-forms $\we^i_{~j}$.  
   
We note that, due to condition (\ref{wc2}), matrix
$\we=(\we^i_{~j})$ of the Weyl connection 1-forms belongs to the
5-dimensional defining representation of the Lie algebra
$\coa(3,2)\subset End(5,\bbR)$ of the Lie group 
$\cog(3,2)\subset\glg(5,\bbR)$. Consequently, the Weyl connection
coefficients $\we^i_{jk}$, defined by
$\we^i_{~j}=\we^i_{~jk}\theta^k$ belong to the tensor product
$\coa(3,2)\otimes \bbR^5$, the vector space of dimension (1+10)5=55.

Now we assume that we have an irreducible $\glg(2,\bbR)$ structure
$[g,\ten,A]$ on a 5-manifold $M^5$. Forgetting about $\ten$ gives the Weyl
geometry as before. In particular there is the \emph{unique} Weyl
connection $\we$ associated with $[g,\ten,A]$. But the existence of
a metric compatible class of tensors $\ten$ makes this Weyl geometry 
more special. To analyze it we introduce a \emph{new} connection,
which will be respecting the entire structure $[g,\ten,A]$. 
This is rather a complicated procedure which we describe below. 

Firstly we require that the new connection preserves $[g]$ and
 $[\ten]$:
\begin{eqnarray}
&&Dg_{ij}=-A g_{ij}\label{cc1}\\
&&D\ten_{ijk}=-\tfrac32 A\ten_{ijk}\label{cc2}.
\end{eqnarray}

This does not determine the connection uniquely -- to have the
uniqueness we need to specify what the torsion of $D$ is. 
We need some preparations to discuss it.

\begin{definition}\label{32de}
Let $(g,\ten,A)$ be a representative of an irreducible $\glg(2,\bbR)$
structure on a 5-dimensional manifold $M^5$. A coframe $\theta^i$,
$i=0,1,2,3,4$, on $M^5$ is called \emph{adapted} to the representative 
$(g,\ten,A)$ if 
$$g=g_{ij}\theta^i\theta^j=\theta^0\theta^4-4\theta^1\theta^3+3(\theta^2)^2$$
and
$$\ten=\ten_{ijk}\theta^i\theta^j\theta^k=3\sqrt{3}(\theta^0\theta^2\theta^4+2\theta^1\theta^2\theta^3-(\theta^2)^3-\theta^0(\theta^3)^2-\theta^4(\theta^1)^2).$$
\end{definition}

Locally such a coframe always exists and is given up to a 
$\glg(2,\bbR)$ transformation. 

Let us now choose an adapted coframe $\theta^i$ to a representative
$(g,\ten,A)$ of $[g,\ten,A]$. In this coframe equations
(\ref{cc1})-(\ref{cc2}) can be rewritten in terms of the connection
1-forms $\Gamma^i_{~j}$ as
\begin{eqnarray}
&&\Gamma^l_{~i}g_{lj}+\Gamma^l_{~}jg_{li}=Ag_{ij}\label{cc11}\\
&&\Gamma^l_{~i}\ten_{ljk}+\Gamma^l_{~j}\ten_{ilk}+\Gamma^l_{~k}\ten_{ijl}=\tfrac32A\ten_{ijk}.\label{icacc}
\end{eqnarray}
When we contract the first equation in indices $i$ and $j$ we get
\be
A=\tfrac25 \Gamma^l_{~l}=\tfrac25 Tr(\Gamma).
\ee
Inserting this into (\ref{icacc}) we get
\be
\Gamma^l_{~i}\ten_{ljk}+\Gamma^l_{~j}\ten_{ilk}+\Gamma^l_{~k}\ten_{ijl}=\tfrac35\Gamma^l_{~l}\ten_{ijk}.\label{icac}
\ee 
Comparing this with (\ref{ica}) we see that the \emph{general solution} for
the connection 1-forms $\Gamma^i_{~j}$ are given by (\ref{om1}), i.e. 
$$\Gamma=\Om_- E_-+\Om_+E_++\Om_0 E_0+\Om_1 E_1,$$
where  
$(\Om_-,\Om_+,\Om_0,\Om_1)$ are four 1-forms on $M^5$ such
that
\be
\Om_1=-\tfrac18 A.
\ee

To fix the remaining three 1-forms 
$(\Om_-,\Om_+,\Om_1)$ we introduce an operator
$$\bar{\ten}:\coa(3,2)\otimes\bbR^5\to S^4\bbR^5$$
defined by:
$$\bar{\ten}(\we)_{ijkm}=\ten_{l(ij}\we^l_{~km)}-\tfrac15\we^l_{~l(m}\ten_{ijk)},$$
and analyze its kernel $\ker\bar{\ten}$. 

Writing equation (\ref{icac}) in terms of the coefficients 
$\Gamma^l_{~im}\in\gla(2,\bbR)\otimes\bbR^5$ and symmetrizing it over 
the indices
$\{imjk\}$, we see that the \emph{whole} $\gla(2,\bbR)\otimes\bbR^5$ is included
in $\ker\bar{\ten}$. 

We use the metric to identify $\bbR^5$ with
$(\bbR^5)^*$, and more generally to identify tensor spaces
$\bgt^k(\bbR^5)^*\bgt^l\bbR^5$ with
$\bgt^{(k+l)}(\bbR^5)^*$. This enables us to identify the objects with
upper indices with the corresponding objects with lower 
indices, e.g. $T_{ijk}=g_{il}T^l_{~jk}$. Having in mind these 
identifications we easily see that, due to the antisymmetry in the last
two indices, every 3-form $T_{ijk}=T_{[ijk]}$ is included in
$\ker\bar{\ten}$.

Thus we have:
$$\gla(2,\bbR)\otimes\bbR^5\subset\ker\bar{\ten},$$
$$\bgw^3\bbR^5\subset\ker\bar{\ten}.$$
The following proposition can be checked by a direct calculation
involving the explicit form of the $\gla(2,\bbR)$ representation given in
(\ref{om1}), (\ref{mo21}). 

\begin{proposition}
The vector space $\ker\bar{\ten}$
has the following properties:
$$\ker\bar{\ten}=(\gla(2,\bbR)\otimes\bbR^5)\oplus\bgw^3\bbR^5$$
and 
$$\dim\ker\bar{\ten}=30.$$
\end{proposition} 
Now we interpret the condition
$\we^l_{~im}\in\ker\bar{\ten}$, i. e. the equation
\be
\ten_{l(ij}\we^l_{~km)}=\tfrac15\we^l_{~l(m}\ten_{ijk)},\label{rewc}
\ee
as a \emph{restriction} on possible Weyl connections. Let us assume that
we have a structure $(M^5,[g,\ten,A])$ with the Weyl connection coefficients 
$\we^l_{~im}$ satisfying (\ref{rewc}). The coefficients $\we^l_{~im}$ are
written in a coframe adapted to some choice $(g,\ten,A)$. 
It is easy to see, using (\ref{wc2}) and contracting
(\ref{rewc}) over all the free indices with a vector field
$X^i$, that the restriction on the Weyl connection (\ref{rewc})
in coordinate-free language is equivalent to
\be
(\stackrel{W}{\nabla}_X\ten)(X,X,X)=-\tfrac12 A(X)\ten(X,X,X).
\label{nic}
\ee
Here $\stackrel{W}{\nabla}$ denotes the Weyl connection in the
Koszul notation. 

\begin{definition}
An irreducible $\glg(2,\bbR)$ structure $(M^5,[g,\ten,A])$ is called
\emph{nearly integrable} iff its Weyl connection
$\stackrel{W}{\nabla}$ associated to the class $[g,A]$  
satisfies (\ref{nic}).
\end{definition}

A nice feature of nearly integrable structures $(M^5,[g,\ten,A])$ is
that they define the \emph{unique} $\gla(2,\bbR)$-valued connection $\Gamma$.
This follows from the above discussion about the kernel of
$\bar{\ten}$. Indeed, given a nearly integrable structure
$(M^5,[g,\ten,A])$ it is enough to choose a representative
$(g,\ten,A)$ and to write the equation (\ref{nic}) for the Weyl
connection $\we$ in an adapted coframe $\theta^i$. Then the uniquely
given Weyl connection coefficients $\we_{ijk}$ are by definition in
$\ker\bar{\ten}=(\gla(2,\bbR)\otimes\bbR^5)\oplus\bgw^3\bbR^5$, which 
means that they \emph{uniquely} split onto
$\Gamma_{ijk}\in\gla(2,\bbR)\otimes\bbR^5$ and 
$\tfrac12 T_{ijk}\in\bgw^3\bbR^5$.
Thus, for all nearly integrable structures $(M^5,[g,\ten,A])$, in a
coframe adapted to $(g,\ten,A)$, we have
\be
\we_{ijk}=\Gamma_{ijk} +\tfrac12 T_{ijk},\label{unic}
\ee
and both $\Gamma_{ijk}\in\gla(2,\bbR)\otimes\bbR^5$ and 
$T_{ijk}\in\bgw^3\bbR^5$ are uniquely determined in terms of
$\we_{ijk}$. Now we rewrite the torsionfree
condition (\ref{wc1}) for the \emph{Weyl} connection in the form
\be
\der\theta^i+\Gamma^i_{~j}\dz\theta^j=\tfrac12 
T^i_{jk}\theta^j\dz\theta^k.\label{wc11}
\ee
It can be interpreted as follows: The nearly integrable structure
$(M^5,[g,\ten,A])$, via (\ref{unic}), 
uniquely determines $\gla(2,\bbR)$-valued connection $\Gamma_{ijk}$ which
respects the structure $[g,\ten,A]$ due to (\ref{cc1}), (\ref{cc2}),
and has \emph{totally skew symmetric torsion} $T_{ijk}$ due to
(\ref{wc11}).

We summarize this part of our considerations in the following
\begin{proposition}\label{chaco}
Every nearly integrable $\glg(2,\bbR)$ structure $(M^5,[g,\ten,A])$
defines a unique $\gla(2,\bbR)$-valued connection which has totally
skew symmetric torsion.
\end{proposition}
Also the converse is true:
\begin{proposition}\label{chacon}
Let $(M^5,[g,\ten,A])$ be an irreducible $\glg(2,\bbR)$ structure and
$\we_{ijk}$ be the Weyl connection coefficients associated, in an adapted
coframe $\theta^i$, with the Weyl
structure $[g,A]$. Assume that the Weyl structure $[g,A]$ admits a
split
$$\we_{ijk}=\Gamma_{ijk}+\tfrac12 T_{ijk},$$ 
in which $\Gamma_{ijk}\in\gla(2,\bbR)\otimes\bbR^5$ and
$T_{ijk}\in\bgw^3\bbR^5$. Then $[g,\ten,A]$ is nearly integrable, the split is unique and $\Gamma_{ij}=\Gamma_{ijk}\theta^k$ is a $\gla(2,\bbR)$-valued connection 
with totally skew symmetric torsion
$\Theta_i=\tfrac12 T_{ijk}\theta^j\dz\theta^k$. 
\end{proposition}

\begin{definition}
The unique $\gla(2,\bbR)$-valued 
connection with totally skew symmetric torsion naturally associated with a
nearly integrable structure $(M^5,[g,\ten,A])$ is called the
\emph{characteristic connection}.
\end{definition}
In the next paragraph we analize algebraic structure of torsion and curvature of characteristic connections.

\subsection{Nearly integrable $\Gl$ structures}\label{s.nearly}
Let $(M^5,[g,\ten,A])$ be a nearly integrable $\Gl$ structure and let $\Gamma$ be its characteristic connection. Then the $\Gl$ invariant information about $(M^5,[g,\ten,A])$ is encoded
in its totally skew symmetric torsion $\Theta_i=\tfrac12 T_{ijk}\theta^i\dz\theta^k$ and its curvature
$$\Omega_{ij}=\tfrac12 R_{ijkl}\theta^k\dz\theta^l=\der\Gamma_{ij}+\Gamma_{ik}\dz\Gamma^k_{~j}.$$
The spaces $\bgw^3\bbR^5$
and $\gla(2,\bbR)\otimes\bgw^2\bbR^5$ are \emph{reducible} under the action
of $\glg(2,\bbR)$. Their decompositions into the $\glg(2,\bbR)$
irreducible components may be used to classify the torsion types, in
the case of $\bgw^3\bbR$, and the curvature types, in the case of
$\gla(2,\bbR)\otimes\bgw^2\bbR^5$. In particular, to decompose
$\bgw^3\bbR^5$ we use the Hodge star operation associated with
one of the metrics $g$ from the class $[g,\ten,A]$. This identifies
$\bgw^3\bbR^5$ with $\bgw^2\bbR^5$. The $\glg(2,\bbR)$ invariant
decomposition of $\bgw^3\bbR^5$ is then transformed to the
decomposition of $\bgw^2\bbR^5$. This is achieved in terms of the 
operator
$$Y_{ijkl}=4\ten_{ijm}\ten_{klp}g^{mp}.$$
This, viewed as an endomorphism of $\bgt^2\bbR^5$ given by
$$Y(w)_{ik}=g^{mj}g^{pl}Y_{ijkl}w_{mp},$$
has the following eigenspaces:
\begin{eqnarray*}
  &&\bgs_1       = \{~S\in\bgt^2\bbR^5~|~ Y (S)= 14 \cdot S~\} = \{S=\lambda
  \cdot g,\;\lambda\in\bbR~\}, \\
 && \bgw_3 = \{ ~F\in\bgt^2 \bbR^5~|~  Y (F) = 7 \cdot F~\}~
=~\sla(2,\bbR),\\
 && \bgs_5 = \{~ S\in\bgt^2\bbR^5 ~| ~ Y(S) = -3 \cdot S~\}, \\
 && \bgw_7 = \{ ~F\in\bgt^2\bbR^5 ~| ~  Y(F) = -8 \cdot F~\},\\
 &&\bgs_9 = \{~ S\in\bgt^2\bbR^5 ~| ~ Y(S) = 4\cdot S~\}.
\end{eqnarray*}
Here the index $k$ in $\bgs_k$ or $\bgw_k$ denotes the dimension of
the eigenspace.

The decomposition
\be\label{e.decomp}\bgt^2\bbR^5=\bgs_1\oplus\bgs_5\oplus\bgs_9\oplus\bgw_3\oplus\bgw_7\ee
is $\glg(2,\bbR)$ invariant. All the components in this decomposition are
$\glg(2,\bbR)$-irreducible. We have the
following 
\begin{proposition}
Under the action of $\glg(2,\bbR)$ the irreducible components of $\bgw^3\bbR^5=*\bgw^2\bbR^5$ are
$$\bgw^3\bbR^5=\bgw_3\oplus\bgw_7.$$
\end{proposition}

At this stage an interesting question arises: Can we give examples of
nearly 
integrable $\glg(2,\bbR)$ structures
whose characteristic connection has torsion of a `pure' type
$T_{ijk}\in\bgw_3$? 

In section \ref{secsec} we give an affirmative answer to this
question. Here we only state a useful 

\begin{lemma}\label{lem}
The 3-dimensional vector 
space $\bgw_3$, when expressed in terms of an adapted coframe
$\theta^i$ of
Definition \ref{32de} is
\begin{eqnarray*}
&&\bgw_3={\rm Span}_\bbR\Big\{\theta^0\dz\theta^3-3\theta^1\dz\theta^2,~
\theta^0\dz\theta^4-2\theta^1\dz\theta^3,~\theta^1\dz\theta^4-3\theta^2\dz\theta^3\Big\}.
\end{eqnarray*}
Similarly, in an adapted coframe $\theta^i$, the Hodge dual $*\bgw_3$
of $\bgw_3$ is
\begin{eqnarray*}
*\bgw_3&=&{\rm Span}_\bbR\Big\{-\theta^0\dz\theta^1\dz\theta^4+2\theta^0\dz\theta^2\dz\theta^3,~-\theta^0\dz\theta^2\dz\theta^4+8\theta^1\dz\theta^2\dz\theta^3,\\
&&-\theta^0\dz\theta^3\dz\theta^4+2\theta^1\dz\theta^2\dz\theta^4\Big\}.
\end{eqnarray*}
In particular, torsion $T^i_{~jk}$ of the 
characteristic connection $\Gamma$ in system (\ref{wc11}) is of pure
type in $\bgw_3$ if and only if, in an adapted coframe $\theta^i$, we
have $g_{il}T^l_{~jk}=T_{[ijk]}$, and its
corresponding 3-form $T=\tfrac16
g_{il}T^l_{~jk}\theta^i\dz\theta^j\dz\theta^k\in*\bgw_3$.
\end{lemma}

Now we pass to the analysis of the curvature. 
The curvature tensor $R^i_{~jkl}$ of a $\gl$ connection defines the following objects: 
\begin{align*}
 & R_{ij}=R^k_{~ikj}  & &\text{the Ricci tensor,} \\
 & R=R_{ij}g^{ij}  & &\text{the Ricci scalar,} \\
 & \Rv^{i}=\ten^{ijk}R_{jk}  & &\text{the Ricci \emph{vector},} \\
 & (\der A)_{ij}=\tfrac{2}{5}R^k_{~kij} & &\text{the Maxwell 2-form.} 
\end{align*}
The Ricci tensor belongs to the space $\bgt^2\bbR^5$ and decomposes according to \eqref{e.decomp}. The Ricci symmetric tensor reads
\be\label{e.ric_decomp} R_{(ij)}=\tfrac{1}{5}Rg_{ij}+\tfrac{2}{7}\Rv^k\ten_{ijk}+R^{(9)}_{ij},\ee
where $\tfrac{1}{5}Rg_{ij}$ is its $\bgs_1$ part, $\tfrac{2}{7}\Rv^k\ten_{ijk}$ 
is its $\bgs_5$ part and $R^{(9)}_{ij}$ is its $\bgs_9$ part defined by \eqref{e.ric_decomp}. 
The antisymmetric Ricci tensor decomposes into 
$$R_{[ij]}=R^{(3)}_{ij}+R^{(7)}_{ij} $$
with the respective $\bgw_3$ and $\bgw_7$ components given by
$$ R^{(3)}_{ij}=\frac{8}{15}R_{[ij]}+\tfrac{1}{15} Y(R_{[\,]})_{ij},$$
$$ R^{(7)}_{ij}=\frac{7}{15}R_{[ij]}-\tfrac{1}{15} Y(R_{[\,]})_{ij}.$$
Here $Y(R_{[\,]})$ denotes the value of the operator $Y$ on $R_{[ij]}$.
Likewise, for the Maxwell form we have 
$$(\der A)_{ij}=\der A^{(3)}_{ij}+\der A^{(7)}_{ij} $$
and
$$ \der A^{(3)}_{ij}=\frac{8}{15}(\der A)_{ij}+\tfrac{1}{15} Y(\der A)_{ij},$$
$$ \der A^{(7)}_{ij}=\frac{7}{15}(\der A)_{ij}-\tfrac{1}{15} Y(\der A)_{ij}.$$
The Ricci tensor and and the Maxwell 2-form have $25+10=35$ coefficients out of total number of $40$ coefficients of the curvature. Since, c.f. \cite{lie},
\ben\gl\otimes\bgw^2\bbR^5=\bgs_1\oplus 2\bgw_3 \oplus 2\bgs_5\oplus 2\bgw_7\oplus\bgs_9,\een
the remaining $5$ parameters are related to the coefficients of a
vector field $K^m$, which is independent of the Ricci tensor. It is
defined in terms of the totally skew symmetric part of the curvature. Using the volume form $\eta^{ijklm}$, we have
$$K^m=R_{ijkl}\eta^{ijklm},$$
and the so defined $K^m$ yields the missing five components of the curvature.
Thus we have the following
\begin{proposition}
The irreducible components of the curvature $R_{ijkl}$ of a characteristic connection are given by
$$R,\quad \Rv^i,\quad R^{(9)}_{~ij},\quad R^{(3)}_{~ij},\quad R^{(7)}_{~ij},\quad \der A^{(3)}_{~ij},\quad 
\der A^{(7)}_{~ij}, \quad K^i.$$
\end{proposition}

It is interesting to ask what is the decomposition of the curvature if
the characteristic connection has torsion in three-dimensional
representation $\bgw_3$. It appears tha it has a very special algebraic form. Writing the structural equations for a characteristic connection with the torsion in $\bgw_3$ 
\begin{eqnarray*}
T&=&\tfrac{1}{12}t_1(-\theta^0\dz\theta^1\dz\theta^4+2\theta^0\dz\theta^2\dz\theta^3)+\\                              
 &&+\tfrac{1}{12}t_2(-\theta^0\dz\theta^2\dz\theta^4+8\theta^1\dz\theta^2\dz\theta^3)+\\
 &&+\tfrac14t_3(-\theta^0\dz\theta^3\dz\theta^4+2\theta^1\dz\theta^2\dz\theta^4)
\end{eqnarray*}
and utilising Bianchi identites, we get the following
\begin{proposition}\label{p.purtor}
Let $\Gamma$ be a characteristic connection of a $\Gl$ structure with torsion in $\bgw_3$ given above. Then
\begin{itemize}
\item The Ricci tensor component $R^{(9)}_{ij}=0$, which means that 
$$ R_{(ij)}=\tfrac{1}{5}Rg_{ij}+\tfrac{2}{7}\Rv^k\ten_{ijk}.$$
\item The skew symmetric Ricci tensor and the Maxwell 2-form
  are related by 
$$\der A^{(3)}=4R^{(3)},\qquad\qquad \der A^{(7)}=\tfrac{2}{3}R^{(7)}.$$
\item The Ricci vector $\Rv$ is fully determined by $T$:
\begin{align*}
 \Rv^i=&(40)^2(*T)_{jk}(*T)_{lm}g^{kl}\ten^{jmi} \\
     =&\frac{7}{6}\sqrt{3}\left( t_3^2,\ -\frac{1}{3}t_2t_3,\ \frac{1}{9}t_1t_3+\frac{2}{27}t_2^2,\  -\frac{1}{9}t_1t_2,\  \frac{1}{9}t^2_1 \right).
\end{align*}
\end{itemize}
Thus, the curvature is fully described by $T$, $R$, $\der A^{(3)}$, $\der A^{(7)}_{ij}$ and $K^{i}$.
\end{proposition}

\subsection{Arbitrary $\Gl$ structures}\label{s.gl_ogolnie}
So far we have been able to introduce a unique $\glg(2,\bbR)$-valued
connection for a nearly integrable
$(M^5,[g,\ten,A])$ only. Nevertheless such a 
connection can be always introduced. To see this consider a
$\glg(2,\bbR)$-invariant conformal pairing in $\coa(3,2)\otimes\bbR^5$
given by
$$(\we,\we')=g^{il}g^{jm}g^{kp}\we_{ijk}\we'_{lmp},$$
where $\we,\we'\in \coa(3,2)\otimes\bbR^5$. We use 
the orthogonal complement of
$\ker\bar{\ten}\subset\coa(3,2)\otimes\bbR^5$ 
with respect to this pairing:
$$\ker\bar{\ten}^\perp=\{\we\in  \coa(3,2)\otimes\bbR^5~{\rm s.t}~(\ker\bar{\ten},\we)=0\}.$$
This vector space is 30-dimensional. It contains a 5-dimensional
subspace spanned by $g_{ij}A_m$, which is related to the $\bbR$ factor in
the split $\gla(2,\bbR)=\bbR\oplus\sla(2,\bbR)\subset\coa(3,2)=\bbR\oplus\soa(3,2)$. Thus it is reasonable to
consider the intersection, say $V_{25}$, of this 30-dimensional space with
$\soa(3,2)\otimes\bbR^5$. This 25-dimensional space 
$$V_{25}=\ker\bar{\ten}^\perp\cap(\soa(3,2)\otimes\bbR^5)$$
has, in turn, zero intersection with
$(\gla(2,\bbR)\otimes\bbR^5)\oplus\bgw^3\bbR^5$
and provides the $\glg(2,\bbR)$ invariant decomposition of
$\coa(3,2)\otimes\bbR^5$:
$$\coa(3,2)\otimes\bbR^5=(\gla(2,\bbR)\otimes\bbR^5)\oplus\bgw^3\bbR^5\oplus
V_{25}.$$
Therefore, if we choose a coframe adapted to a representative $(g,\ten,A)$ we can
uniquely decompose the Weyl connection coefficients 
$\we_{ijk}\in\coa(3,2)\otimes\bbR^5$ of our \emph{arbitrary}
$\glg(2,\bbR)$ structure according to 
$$\we_{ijk}=\Gamma_{ijk}+\tfrac12 B_{ijk}.$$
Now $\Gamma_{ijk}\in\gla(2,\bbR)\otimes\bbR^5$, and they are interpreted as
new connection coefficients; the tensor $B_{ijk}$ belongs to $\bgw^3\bbR^5\oplus
V_{25}$ and its antisymmetrization $T_{ijk}=B_{i[jk]}$ is now interpreted 
as the torsion of $\Gamma$.
Thus, \emph{every} $\glg(2,\bbR)$ structure
$(M^5,[g,\ten,A])$ uniquely defines a $\gla(2,\bbR)$-valued connection
with torsion in $\bgw^3\bbR^5\oplus
V_{25}$. The torsion is not totally skew anymore. Space $V_{25}$ further 
decomposes onto the $\glg(2,\bbR)$-irreducible components according
to 
$V_{25}=\bgs_5\oplus\bgs_9\oplus\bgs_{11}$. The $\glg(2,\bbR)$ structures
equipped with the unique $\gla(2,\bbR)$ connection which has torsion
in $V_{25}$ find application in the
theory of integrable equations of hydrodynamic type 
\cite{fer1}.

\section{$\glg(2,\bbR)$ bundle} \label{nabu}
\emph{First}, we describe an irreducible $\glg(2,\bbR)$ structure $[g,\ten,A]$ on $M^5$ in the language of principal bundles. 
 
Every irreducible $\glg(2,\bbR)$
structure $[g,\ten,A]$ on a 5-manifold $M^5$ defines the
9-dimensional bundle $\glg(2,\bbR)\to P\to M^5$, the $\glg(2,\bbR)$
reduction of the bundle of linear frames $\glg(5,\bbR)\to F(M^5)\to M^5$.
If $[g,\ten,A]$ is equipped with a $\gl$ connection $\Gamma$ 
then the structural equations on $M^5$ read
$$\der\omega^i+\Gamma^i_{~j}\dz\omega^j=\tfrac12
T^i_{~jk}\omega^j\dz\omega^k,$$
$$\der\Gamma^i_{~j}+\Gamma^i_{~k}\dz\Gamma^k_{~j}=\tfrac12 R^i_{~jkl}\omega^k\dz\omega^l.$$
Here $(\omega^i)$ is an adapted coframe and $\Gamma=(\Gamma^i_{~j})$ is written in the
representation (\ref{om1}). We lift these structural equations to
$P$ obtaining:
\begin{eqnarray}
&&\der\theta^0=4(\Om_1+\Om_0)\dz\theta^0-4\Om_+\dz\theta^1+\tfrac12 T^0_{~ij}\theta^i\dz\theta^j,\nonumber\\
&&\der\theta^1=-\Om_-\dz\theta^0+(4\Om_1+2\Om_0)\dz\theta^1-3\Om_+\dz\theta^2+\tfrac12 T^1_{~ij}\theta^i\dz\theta^j,\nonumber\\
&&\der\theta^2=-2\Om_-\dz\theta^1+4\Om_1\dz\theta^2-2\Om_+\dz\theta^3+\tfrac12 T^2_{~ij}\theta^i\dz\theta^j,\nonumber\\
&&\der\theta^3=-3\Om_-\dz\theta^2+(4\Om_1-2\Om_0)\dz\theta^3-\Om_+\dz\theta^4+\tfrac12 T^3_{~ij}\theta^i\dz\theta^j,\nonumber\\
&&\der\theta^4=-4\Om_-\dz\theta^3+4(\Om_1-\Om_0)\dz\theta^4+\tfrac12 T^4_{~ij}\theta^i\dz\theta^j,\label{sysc}\\
&&\der\Om_+=2\Om_0\dz\Om_++\tfrac12 R_{+ij}\theta^i\dz\theta^j,\nonumber\\
&&\der\Om_-=-2\Om_0\dz\Om_-+\tfrac12 R_{-ij}\theta^i\dz\theta^j,\nonumber\\
&&\der\Om_0=\Om_+\dz\Om_-+\tfrac12 R_{0ij}\theta^i\dz\theta^j,\nonumber\\
&&\der\Om_1=\tfrac12 R_{1ij}\theta^i\dz\theta^j,\nonumber
\end{eqnarray}
with the forms $\theta^i$ being the components of the canonical $\bbR^5$-valued form $\theta$ on $P$, c.f. \cite{kono}. 
In a coordinate system $(x,a)$ on $P$, $x\in M^5$, $a\in\Gl$, which 
is compatible with the local trivialisation $P\cong M^5\times\Gl$ they are given by
\ben\theta^i(x,a)=(a^{-1})^i_{~j}\omega^j(x).\een
The connection forms $(\Om_-,\Om_+,\Om_0,\Om_1)$ are defined in terms
of (\ref{mo21}) via 
\ben\Om_- (E_-)^i_{~j}+\Om_+(E_+)^i_{~j}+\Om_0 (E_0)^i_{~j}+\Om_1 (E_1)^i_{~j}=(a^{-1})^i_{~k}\Gamma^k_{~l}(x)a^l_{~j}+(a^{-1})^i_{~k}\der a^k_{~j}.\een
Note that $(\theta^1,\theta^1,\theta^2,\theta^3,\theta^4,\Om_-,\Om_+,\Om_0,\Om_1)$ is a coframe on $P$ and the class of 1-forms $[A]$ lifts to a 1-form $\tilde{A}=-8\Om_1$.

\emph{Second}, we change the point of view.
Suppose that we \emph{are given} a nine dimensional manifold 
$P$ equipped with a coframe of nine 1-forms
$(\theta^0,\theta^1,\theta^2,\theta^3,\theta^4,\Om_-,$ $\Om_+,\Om_0,\Om_1)$
on it. Suppose that these linearly indpendent forms, together with
some functions $T^i_{~jk}$, $R^l_{~ijk}$, satisfy the
system (\ref{sysc}) on $P$. What we can say about such a 9-dimensional
manifold $P$?

To answer this question consider a distribution $\mathfrak{h}$ on $P$
which annihilates the forms
$(\theta^0,\theta^1,\theta^2,\theta^3,\theta^4)$:
$$\mathfrak{h}=\{X\in{\rm T}P~{\rm
  s.t.}~X\hook\theta^i=0,~i=0,1,2,3,4\}.$$
Then the first five equations of the system (\ref{sysc}) guarantee
  that the forms $(\theta^0,\theta^1,\theta^2,$ $\theta^3,\theta^4)$
  satisfy the Fr\"obenius condition,
$$\der\theta^i\dz\theta^0\dz\theta^1\dz\theta^2\dz\theta^3\dz\theta^4=0,\quad\quad\forall~i=0,1,2,3,4$$
and that, in turn, the distribution $\mathfrak h$ is integrable. Thus
manifold $P$ is foliated by 4-dimensional leaves tangent to the
distribution $\mathfrak h$. 

Now on $P$ we consider two multilinear symmetric forms. The bilinear
one, defined by
\be\tilde{g}=\theta^0\theta^4-4\theta^1\theta^3+3(\theta^2)^2,
\label{fex1}\ee
and the three-linear one given by
\be
\tilde{\ten}=3\sqrt{3}(\theta^0\theta^2\theta^4+2\theta^1\theta^2\theta^3-(\theta^2)^3-\theta^0(\theta^3)^2-\theta^4(\theta^1)^2).
\label{fex2}\ee

Of course, since the 1-forms$(\Om_-,\Om_+,\Om_0,\Om_1)$ 
are not present in the definitions (\ref{fex1}), (\ref{fex2}), then
$\tilde{g}$ and $\tilde{\ten}$ are \emph{degenerate}. 
For example, the signature of the bilinear form $\tilde{g}$ is
$(+,+,+,-,-,0,0,0,0)$. The degenerate directions for these two forms
are just the directions tangent to the leaves of the foliation
generated by $\mathfrak h$. Let us denote by
$(X_0,X_1,X_2,X_3,X_4,X_5,X_6,X_7,X_8)$ the frame of vector fields on
$P$ 
dual to the 1-forms 
$(\theta^0,\theta^1,\theta^2,\theta^3,\theta^4,\Om_-,\Om_+,\Om_0,\Om_1)$.
In particular $(X_5,X_6,X_7,X_8)$ constitutes a basis for $\mathfrak
h$, and we have $X_\mu\hook\theta^i=0$ for each $\mu=5,6,7,8$
and $i=0,1,2,3,4$. Using this, and the exterior derivatives of $\theta^i$ given
in the first five equations (\ref{sysc}), we easily find the Lie
derivatives of $\tilde{g}$ and $\tilde{\ten}$ along the directions
  tangent to the leaves of $\mathfrak h$. These are:
$${\mathcal
    L}_{X_\mu}\tilde{g}=8(X_\mu\hook\Om_1)\tilde{g},\quad\quad
  {\mathcal
    L}_{X_\mu}\tilde{\ten}=12(X_\mu\hook\Om_1)\tilde{\ten},\quad\quad\forall
    \mu=5,6,7,8.$$
Moreover, if we denote 
\be
\tilde{A}=-8\Om_1,\label{fex3}
\ee
and we use the last of equations (\ref{sysc}), we also find that
$$ {\mathcal
    L}_{X_\mu}\tilde{A}=-8\der(X_\mu\hook\Om_1),\quad\quad\forall\mu=5,6,7,8.$$
This is enough to deduce that the objects
    $(\tilde{g},\tilde{\ten},\tilde{A})$ descend to the 5-dimensional
    leaf space $M^5=P/\mathfrak{h}$. There they define a conformal class of triples 
$(g,\ten,A)$ with the transformation rules $g\to {\rm e}^{2\phi}g$,
    $\ten\to{\rm e}^{3\phi}\ten$, $A\to A-2\der\phi$. Due to the fact
    that, when passing to the quotient $M^5=P/\mathfrak{h}$, we reduced the
    degenerate directions of $\tilde{g}$ and $\tilde{\ten}$ to points
    of $M^5$, the resulting descended triples $(g,\ten,A)$ have
    non-degenerate $g$ of signature
    $(3,2)$ and non-degenerate $\ten$. It is clear that together with
    $A$ they define an irreducible $\glg(2,\bbR)$ structure on
    $M^5$: a section $s\colon M^5\to P$ is an adapted coframe on $M^5$, the triple $(s^*\tilde{g},s^*\tilde{\Upsilon},s^*\tilde{A})$ is a representative of the structure, the forms $s^*\Om_-,s^*\Om_+,s^*\Om_0,s^*\Om_1$ are $\gl$ connection 1-forms on $M^5$ and $s^*T$, $s^*R$ are torsion and curvature of this connection, respectively. We have the following
\begin{proposition}\label{p.constr}
Every 9-dimensional manifold $P$ equipped with nine 1-forms
$(\theta^0,\theta^1,\theta^2,\theta^3,\theta^4,\Om_-,\Om_+,\Om_0,\Om_1)$
which 
\begin{itemize}
\item are linearly independent at every point of $P$, 
\item satisfy system (\ref{sysc}) with some functions $T^i_{~jk}$, $R^i_{~jkl}$ on $P$, 
\end{itemize}
is foliated by 4-dimensional leaves over a 5-dimensional space $M^5$,
which is the base for the fibration $P\to M^5$. The manifold $M^5$ is equipped with a natural irreducible $\glg(2,\bbR)$
structure $[g,\ten,A]$ and a $\gl$ connection compatible with it. The 
torsion and the curvature of this connection is given by $T^i_{~jk}$ and $R^i_{~jkl}$.
\end{proposition}

\section{5th order ODE as nearly integrable $\glg(2,\bbR)$ geometry with `small' torsion. Main theorem}
\label{secsec}
A large number of examples of nearly integrable $\glg(2,\bbR)$
structures in dimension five is related to 5th order ODEs. This is
mainly due to the following, well known,
\begin{proposition}\label{odefl}
An ordinary differential equation
$y^{(5)}=0$
has $\glg(2,\bbR)\times_{\rho_5}\bbR^5$ as its group of contact
symmetries. Here $\rho_5 :\glg(2,\bbR)\to\glg(5,\bbR)$ is the
5-dimensional irreducible representation of $\glg(2,\bbR)$. 
\end{proposition}

To explain the above statement we consider a general 5th order ODE
\be
y^{(5)}=F(x,y,y',y'',y^{(3)},y^{(4)})\label{piec}
\ee
for a real function $\bbR\ni x\mapsto y(x)\in\bbR$. Let us introduce the notation
$y_1=y'$, $y_2=y''$, $y_3=y^{(3)}$, $y_4=y^{(4)}$ and $F_i=\frac{\partial F}{\partial y_i}$, $i=1,2,3,4$, $F_y=\frac{\partial F}{\partial y}$. The functions $(x,y,y_1,y_2,y_3,y_4)$ form a local coodinate system in the 4-order jet space $J$ of curves in $\bbR^2$. 
Define a total derivative, which is a vector field in $J$ 
\be
 \D=\partial_x+y_1\partial_y+y_2\partial_{y_1}+y_3
\partial_{y_2}+y_4\partial_{y_3}+F\partial_{y_4}.\label{D}
\ee
With the help of $\D$ the derivatives are given by formulae $y_1=\D y/\D x$, $y_2=\D y_1/\D x$ and so on, up to $y_5=\D y_4/\D x$.

A contact transformation of variables in a 5-order ODE is a
transformation that mixes the independent variable $x$, the dependent
variable $y$ and the first derivative $y_1$ in such a way that the
meaning of the first derivative is retained:

\begin{definition}
A contact transformation of variables is an invertible, sufficiently smooth transformation of the form
\be\label{contact}
\bma x \\ y \\ y_1 \ema 
 \mapsto  
\bma \bar{x} \\ \bar{y} \\ \bar{y}_1 \ema 
 = 
\bma \bar{x}(x,y,y_1) \\
     \bar{y}(x,y,y_1)  \\
     \bar{y}_1(x,y,y_1)
\ema
\ee
satisfying the condition
$$\bar{y}_1= \frac{\D\bar{y}}{\D\bar{x}}.\quad \text{(preservation of first derivative)}$$
The higher order derivatives are given by the iterative formula
$$y_{n+1}\mapsto \bar{y}_{n+1}=\frac{\D\bar{y}_n}{\D\bar{x}},\qquad i=1,2,3,4. $$
\end{definition}

Let us now consider the equation $y^{(5)}=0$. We show how the flat torsionfree 5-dimensional irreducible $\glg(2,\bbR)$ structure is naturally generated on its space of solutions by means of the symmetry group.
A solution to $y^{(5)}=0$ is of the form
\be\label{solo} y(x)=c_4x^4+4c_3x^3+6c_2x^2+4c_1x+c_0\ee
with five integration constants $c_0,c_1,c_2,c_3,c_4$. Then a solution
of $y^{(5)}=0$ may be ideantified with a point $c=(c_0,c_1,c_2,c_3,c_4)^T$ in $\bbR^5$. A contact symmetry of $y^{(5)}=0$ is a contact transformation of variables that transforms its solutions into solutions. Group of contact symmetries of $y^{(5)}=0$ is generated by the following one-parameter groups of transformations on the $xy$-plane:
\begin{align*}
 &\varphi^0_t(x,\,y)=(x,\,y+t), &  & \varphi^1_t(x,\,y)=(x,\,y+4xt),\\
 &\varphi^2_t(x,\,y)=(x,\,y+6x^2t), &  & \varphi^3_t(x,\,y)=(x,\,y+4x^3t),\\
 &\varphi^4_t(x,\,y)=(x,\,y+x^4t), &  & \varphi^5_t(x,\,y)=(xe^{2t},\,ye^{4t}),\\
 &\varphi^6_t(x,\,y)=(x,\,ye^{4t}), &  & \varphi^7_t(x,\,y)=(x+t,\,y),\\
 & \varphi^8_t(x,y)=\left(\frac{x}{1+xt},\,\frac{y}{(1+xt)^4}\right) &  &
\end{align*}
and the transformation rules for $y_1$ are given by $\varphi^A(y_1)=\D(\varphi^A(y))/\D(\varphi^A(x))$, 
$A=0,\ldots, 8$.

Transforming \eqref{solo} according to the above formulae we find that $\varphi^0_t,\ldots,\varphi^4_t$ are translations in the space of solutions: 
$$\varphi^0_t(c)=(c_0-t,c_1,c_2,c_3,c_4)^T \qquad ,\ldots, \qquad \varphi^4_t(c)=(c_0,c_1,c_2,c_3,c_4-t)^T,$$ while transformations $\varphi^5_t,\ldots,\varphi^8_t$ generate $\Gl$ and act through the 5-dimensional irreducible representation \eqref{mo21}:
\begin{align*}
 &\varphi^5_t(c)=\exp(tE_0)c, &
 &\varphi^6_t(c)=\exp(tE_1)c, \\
 &\varphi^7_t(c)=\exp(tE_+)c, &
 &\varphi^8_t(c)=\exp(tE_-)c.
\end{align*}
Of course, $\Gl$ stabilizes the origin $(0,0,0,0,0)$ in $\bbR^5$, thus the space of solutions is the homogeneous space $\glg(2,\bbR)\to\glg(2,\bbR)\times_{\rho_5}\bbR^5\to\bbR^5$. The total space of this bundle is equipped with the Maurer -- Cartan form $\omega_{MC}$ of $\glg(2,\bbR)\times_{\rho_5}\bbR^5$. Choosing an approriate basis in $\gl$ and writing explicitly the structural equations $\der\omega_{MC}+\omega_{MC}\dz\omega_{MC}=0$ we get
\begin{eqnarray*}
&&\der\theta^0=4(\Om_1+\Om_0)\dz\theta^0-4\Om_+\dz\theta^1,\\
&&\der\theta^1=-\Om_-\dz\theta^0+(4\Om_1+2\Om_0)\dz\theta^1-3\Om_+\dz\theta^2,\\
&&\der\theta^2=-2\Om_-\dz\theta^1+4\Om_1\dz\theta^2-2\Om_+\dz\theta^3,\\
&&\der\theta^3=-3\Om_-\dz\theta^2+(4\Om_1-2\Om_0)\dz\theta^3-\Om_+\dz\theta^4,\\
&&\der\theta^4=-4\Om_-\dz\theta^3+4(\Om_1-\Om_0)\dz\theta^4,\\
&&\der\Om_+=2\Om_0\dz\Om_+,\\
&&\der\Om_-=-2\Om_0\dz\Om_-,\\
&&\der\Om_0=\Om_+\dz\Om_-,\\
&&\der\Om_1=0,
\end{eqnarray*} 
which is the system \eqref{sysc} with all the torsion and curvature coefficients equal to zero. According to proposition \ref{p.constr} it yields a flat and torsionfree irreducible $\glg(2,\bbR)$ structure on the space of solutions of $y^{(5)}=0$. 
Again, as in the case of the algebraic geometric realization of section \ref{s1}, we learned about that 
from E. V. Ferapontow \cite{fera}. 

We now pass to a more general situation, namely to the equation
(\ref{piec}) with a \emph{general} $F$. The following questions are in
order:

What shall one assume about $F$ to be able to construct an 
irreducible $\glg(2,\bbR)$ structure on the solution space of the
corresponding ODE? Is the case $F=0$ very special, or there are other
ODEs, contact nonequivalent to the $F=0$ case, which define a $\glg(2,\bbR)$
geometry on the solution space? If the answer is affirmative, how do
we find such $F$s and  what can we say about the corresponding
$\glg(2,\bbR)$ structures?

Answer to these questions is given by the following
\begin{theorem}[Main theorem]\label{maint}
Every contact equivalence class of 5th order ODEs satisfying the W\"unschmann conditions 
\begin{align}
&50\D^2F_4 - 75 \D F_3 + 50 F_2 - 60 F_4\D F_4 + 30 F_3F_4 + 8 F_4^3=0, \notag \\ \notag \\
&375 \D^2F_3 - 1000 \D F_2 + 350 \D F_4^2 + 1250 F_1 - 650F_3 \D F_4  + 
    200 F_3^2 - \notag \\
& 150 F_4 \D F_3 + 200 F_2 F_4 - 140 F_4^2 \D F_4 + 
    130 F_3 F_4^2 + 14 F_4^4=0,\notag  \\ \label{wilki} \\
&1250 \D^2F_2 - 6250 \D F_1 + 1750 \D F_3\D F_4 - 2750 F_2\D F_4  - \notag \\
& 875  F_3\D F_3+1250 F_2 F_3 - 500F_4 \D F_2  + 700 (\D F_4)^2 F_4 + \notag \\
& 1250 F_1 F_4 -1050 F_3 F_4\D F_4  + 350 F_3^2 F_4 - 350 F_4^2\D F_3  + \notag \\
& 550 F_2 F_4^2 - 280 F_4^3 \D F_4 + 210 F_3 F_4^3 + 28 F_4^5 + 18750 F_y=0 \notag 
\end{align}   
defines a nearly integrable irreducible $\glg(2,\bbR)$ geometry
$(M^5,[g,\ten,A])$ on the space $M^5$ of its solutions. This geometry
has the characteristic connection with torsion $T$ of the `pure'
type in the 3-dimensional irreducible representation $\bgw_3$. The first structural equation for this connection are the following:
\begin{eqnarray}
\der\theta^0&=&4(\Om_1+\Om_0)\dz\theta^0-4\Om_+\dz\theta^1+\nonumber\\ &&-\tfrac13t_1\theta^0\dz\theta^1-\tfrac13t_2\theta^0\dz\theta^2-t_3\theta^0\dz\theta^3+2t_3\theta^1\dz\theta^2, \nonumber\\
\der\theta^1&=&-\Om_-\dz\theta^0+(4\Om_1+2\Om_0)\dz\theta^1-3\Om_+\dz\theta^2+\nonumber\\
 &&-\tfrac16t_1\theta^0\dz\theta^2-\tfrac14t_3\theta^0\dz\theta^4-\tfrac23t_2\theta^1\dz\theta^2, \nonumber\\
\der\theta^2&=&-2\Om_-\dz\theta^1+4\Om_1\dz\theta^2-2\Om_+\dz\theta^3+\nonumber\\
&&-\tfrac19t_1\theta^0\dz\theta^3+\tfrac{1}{18}t_2\theta^0\dz\theta^4-\tfrac49t_2\theta^1\dz\theta^3-\tfrac13t_3\theta^1\dz\theta^4,\nonumber\\
\der\theta^3&=&-3\Om_-\dz\theta^2+(4\Om_1-2\Om_0)\dz\theta^3-\Om_+\dz\theta^4+\nonumber\\
&&+\tfrac{1}{12}t_1\theta^0\dz\theta^4-\tfrac23t_2\theta^2\dz\theta^3-\tfrac12t_3\theta^2\dz\theta^4,\nonumber\\
\der\theta^4&=&-4\Om_-\dz\theta^3+4(\Om_1-\Om_0)\dz\theta^4+\nonumber\\
&&-\tfrac13t_1\theta^1\dz\theta^4+\tfrac23t_1\theta^2\dz\theta^3-\tfrac13t_2\theta^2\dz\theta^4-t_3\theta^3\dz\theta^4\nonumber 
\end{eqnarray}
with the torsion coefficients
\begin{eqnarray}
t_3&=&\frac{6 (\al^5_{~5})^2}{5\al^1_{~1}}F_{44},\nonumber\\
t_2&=&\frac{9\al^5_{~5}}{50(\al^1_{~1})^2}~[\al^1_{~1}(10\D F_{44}+3F_4F_{44})+5\al^1_{~0}F_{44}],\nonumber\\
&&\nonumber\\
t_1&=&[1000(\al^1_{~1})^3]^{-1}\times\nonumber\\
&&\Big(225(\al^1_{~0})^2F_{44}+90\al^1_{~0}\al^1_{~1}(10\D F_{44}+3F_4F_{44})+\nonumber \\
&&-9(\al^1_{~1})^2[20(5\D F_{34}+20F_{24}-15F_{33}+3F_4\D F_{44}-11F_4F_{34})+\nonumber\\
&&+F_{44}(-120\D F_4+340F_3+51F_4^2)]\Big),\nonumber
\end{eqnarray}
where $(y,y_1,y_2,y_3,y_4,x,\al^1_{~1},\al^1_{~0},\al^5_{~5})$ is a local coordinate system on $\Gl\to P\to M^5$.
The second structural equations are the following:
\begin{eqnarray}
\der\Om_+&=&2\Om_0\dz\Om_++\left(\tfrac{1}{6}b_2-\tfrac {1}{81}t_1^{2}+\tfrac{5}{3}c_5\right)
 \theta_{{0}},\theta_{{1}}
 +\left(-\tfrac {2}{81}t_1t_2-\tfrac{10}{3}c_4+\tfrac {5}{12}b_3 \right)\theta_{{0}}\dz\theta_{{2}}+\nonumber\\
&&+\left(-\tfrac {1}{243}t_2^{2}-\tfrac {1}{162}t_1t_3+\tfrac{10}{3}c_3-\tfrac{1}{30}R+b_4
 -\tfrac{1}{4}a_2 \right)\theta_{{0}}\dz\theta_{{3}}+ \nonumber\\
&&+ \left(\tfrac {1}{54}t_2t_3-\tfrac{1}{8}a_3-\tfrac{5}{3}c_2+\tfrac{1}{12}b_5 \right)
 \theta_{{0}}\dz\theta_{{4}}+ \nonumber\\
&&+\left( -\tfrac{1}{27}t_2^{2}-\tfrac{1}{18}t_1t_3+\tfrac{1}{10}R+2b_4+\tfrac{3}{4}a_2 \right)
 \theta_{{1}}\dz\theta_{{2}}+\nonumber\\
&&+\left(-\tfrac{1}{9}t_2t_3+\tfrac{1}{4}a_3+\tfrac{2}{3}b_5\right)\theta_{{1}}\dz\theta_{{3}}
 + \left(\tfrac{1}{18}t_3^{2}+\tfrac{5}{3}c_1+\tfrac{1}{6}b_6\right) \theta_{{1}}\dz\theta_{{4}}+\notag \\
&&+ \left( -\tfrac {5}{18}t_3^{2}-\tfrac{10}{3}c_1+\tfrac{1}{3}b_6 \right)\theta_{{2}}\dz\theta_{{3}} 
 +\tfrac{1}{4}b_7\theta_{{2}}\dz\theta_{{4}},\nonumber\\
\der\Om_-&=&-2\Om_0\dz\Om_-+\tfrac{1}{4}b_1\theta_{{0}}\dz\theta_{{2}}+ 
 \left( \tfrac{1}{6} b_2-\tfrac {1}{162}t_1^{2}-\tfrac{5}{3}c_5 \right) \theta_{{0}}\dz\theta_{{3}}+\nonumber\\
&&+ \left( -\tfrac {1}{162}t_1t_2+\tfrac{5}{3} c_4+\tfrac{1}{12}b_3-\tfrac{1}{8}a_1\right)
  \theta_{{0}}\dz\theta_{{4}}+\nonumber\\
&&+ \left(\tfrac {5}{162}t_1^{2}+\tfrac{1}{3}b_2+\tfrac{10}{3}c_5 \right) \theta_{{1}}\dz\theta_{{2}}
 +\left( \tfrac{1}{27}t_1t_2+\tfrac{2}{3}b_3+\tfrac{1}{4}a_1 \right) \theta_{{1}}\dz\theta_{{3}}+  \notag\\
&& + \left(b_4-\tfrac{1}{4}a_2+\tfrac {1}{162}t_1t_3+\tfrac {1}{243}t_2^{2}-\tfrac{10}{3}c_3
 +\tfrac{1}{30}R \right)\theta_{{1}}\dz\theta_{{4}} +\nonumber\\
&& +\left( \tfrac{1}{27}t_2^2+\tfrac{1}{18}t_1t_3-\tfrac{1}{10}R+2b_4+\tfrac{3}{4}a_2\right)
 \theta_{{2}}\dz\theta_{{3}} +\nonumber\\
&& + \left(\tfrac {2}{27}t_2t_3+\tfrac{10}{3}c_2+\tfrac {5}{12}b_5 \right)\theta_{{2}}\dz\theta_{{4}} 
 +\left(\tfrac{1}{9}t_3^{2}-\tfrac{5}{3}c_1+\tfrac{1}{6}b_6\right)\theta_{{3}}\dz\theta_{{4}}+ \nonumber \\
&&\tfrac{1}{12}(2t_{14}+t_{23}-6t_2t_3-3t_{32})\theta^2\dz\theta^4+\tfrac{1}{12}(4t_{24}-9t_3^2-3t_{33})
\theta^3\dz\theta^4,\nonumber \\
\der\Om_0&=&\Om_+\dz\Om_--\tfrac{1}{4}b_1\theta_{{0}}\dz\theta_{{1}}
 +\left( -\tfrac{1}{6}b_2-\tfrac {1}{162}t_1^{2}+\tfrac{5}{6}c_5 \right)\theta_{{0}}\dz\theta_{{2}}+\nonumber \\
&&+\left(-\tfrac {1}{54}t_1t_2-\tfrac{1}{12}b_3+\tfrac{1}{8}a_1\right) \theta_{{0}}\dz\theta_{{3}}+\nonumber\\
&&-\left(\tfrac {1}{81}t_1t_3+\tfrac {2}{243}t_2^{2}+\tfrac{5}{6}c_3+\tfrac {1}{60}R \right)
 \theta_{{0}}\dz\theta_{{4}}+\nonumber \\
&&+ \left(\tfrac {1}{162}t_1t_2-\tfrac {20}{3}c_4-\tfrac{1}{6}b_3-\tfrac{3}{8}a_1\right)
 \theta_{{1}}\dz\theta_{{2}}+ \nonumber \\
&&+ \left(-\tfrac {1}{81}t_1 t_3-\tfrac {2}{243}t_2^{2}+\tfrac{20}{3}c_3
 +\tfrac{1}{30}R \right) \theta_{{1}}\dz\theta_{{3}}+ \nonumber \\
&& + \left(-\tfrac{1}{18}t_2t_3-\tfrac{1}{8}a_3+\tfrac{1}{12}b_5 \right) 
 \theta_{{1}}\dz\theta_{{4}}+\nonumber \\
&& +\left(\tfrac {1}{54}t_2t_3+\tfrac{3}{8}a_3-\tfrac {20}{3}c_2+\tfrac{1}{6}b_5 \right)
 \theta_{{2}}\dz\theta_{{3}}+\nonumber \\
&&+ \left(-\tfrac{1}{18}t_3^{2}+\tfrac{5}{6}c_1+\tfrac{1}{6}b_6 \right) \theta_{{2}}\dz\theta_{{4}} 
+\tfrac{1}{4}b_7\theta_{{3}}\dz\theta_{{4}},\nonumber\\
 \der\Om_1&=&-\tfrac{1}{8}b_1\theta_{{0}}\dz\theta_{{1}}-\tfrac{1}{8}b_2\theta_{{0}}\dz\theta_{{2}}-\tfrac{1}{8} \left(b_3+a_1 \right)\theta_{{0}}\dz\theta_{{3}} 
 -\tfrac{1}{8} \left(b_4+a_2\right)\theta_{{0}}\dz\theta_{{4}}+\notag \\
&&+ \left( \tfrac{3}{8}a_1-\tfrac{1}{4}b_3\right) \theta_{{1}}\dz\theta_{{2}}
 + \left( \tfrac{1}{4}a_2 -b_4\right) \theta_{{1}}\dz\theta_{{3}}
 -\tfrac{1}{8}\left(a_3+b_5 \right)\theta_{{1}}\dz\theta_{{4}} + \notag \\
 &&+\left( \tfrac{3}{8}a_3-\tfrac{1}{4}b_5\right) \theta_{{2}}\dz\theta_{{3}} 
 -\tfrac{1}{8}b_6\theta_{{2}}\dz\theta_{{4}}-\tfrac{1}{8}b_7\theta_{{3}}\dz\theta_{{4}}, \nonumber
\end{eqnarray}
 where $a_1$, $a_2$, $a_3$, $b_1$, $b_2$, $b_3$, $b_4$, $b_5$, $b_6$, $b_7$, $c_1$, $c_2$, $c_3$, $c_4$,  $c_5$ and $R$ are functions. All of these functions but $R$ are determined by the differentials of torsions:
\begin{eqnarray}
\der t_1&=& 2t_2\Om_--2t_1\Om_0-4t_1\Om_1+\tfrac{3}{2}b_1\theta_{{0}}
  + \left( 2b_2-\tfrac{4}{27}t_1^{2}+20 c_5 \right) \theta_{{1}}+\notag \\
 &&+\left( -\tfrac49 t_1t_2-60c_4+3b_3-\tfrac{9}{2}a_1 \right) \theta_{{2}}+ \notag \\
 &&+ \left( -\tfrac49 t_1t_3-\tfrac{8}{27}t_2^{2}+60c_3+6b_4-9a_2\right) \theta_{{3}}+\notag \\
 &&+\left(-\tfrac49t_2t_3-\tfrac92a_3-20c_2+\tfrac12 b_5 \right) \theta_{{4}},\notag \\
\der t_2 &=& 3t_3\Om_-+t_1\Om_+-4t_2\Om_1+ \left( \tfrac12 b_2+\tfrac{2}{27}t_1^2-10c_5 \right) 
 \theta_{{0}}+\notag \\ 
 &&+\left(\tfrac{4}{27}t_1t_2+20c_4+2b_3+\tfrac92 a_1\right)\theta_1
 + 9\left(a_2+b_4\right) \theta_{{2}}+\notag \\
 &&+\left(-\tfrac49 t_2t_3+\tfrac92 a_3-20 c_2+2 b_5 \right) \theta_{{3}}+
 \left(-\tfrac23 t_3^2+10c_1+\tfrac12 b_6 \right)\theta_{{4}}, \notag \\
\der t_3 &=& 2t_3\Om_{{0}}+\tfrac23 t_2\Om_+-4t_3\Om_1
+\left(\tfrac{4}{81}t_1t_2+\tfrac{20}{3}c_4+\tfrac16 b_3-\tfrac32 a_1 \right) \theta_{{0}}+\notag \\ 
 &&+\left(\tfrac{4}{27}t_1t_3+\tfrac{8}{81}t2^2-20c_3+2b_4-3a_2 \right) \theta_1+ \notag \\
 &&+\left(\tfrac49 t_2 t_3-\tfrac32 a_3+20 c_2+b_5 \right) \theta_2
  +\left(\tfrac49 t_3^2-\tfrac{20}{3}c_1+\tfrac23 b_6\right) \theta_3+\tfrac12 b_7\theta_{{4}}. \notag
\end{eqnarray}
The function $R$ is the Ricci scalar for the connection.
\end{theorem}

Before presenting the proof let us notice several facts and corollaries.

The theorem guarantees that every equivalence class of ODEs 
satisfying conditions \eqref{wilki} has its
corresponding nearly integrable $\glg(2,\bbR)$ geometry
$(M^5,[g,\ten,A])$ with torsion in $\bgw_3$.  It may happen, however,
that there are contact \emph{non-equivalent} classes of ODEs defining \emph{the same} $\Gl$ geometries. (See also remark \ref{r.loss}).

The W\"unschmann conditions, although very complicated, possess nontrivial solutions. For example the equation
$$y^{(5)}=c\Big(\frac{5y^{(3)3}( 5 - 27cy''^{2})}{9( 1+c y''^{2})^2} +
10\frac{y''y^{(3)}y^{(4)}}{ 1+c y''^2}\Big),$$
where $c=\pm 1$ satisfies the W\"unschmann conditions and is not contact equivalent to $F=0$. Other examples are considered in section \ref{s.examples}.

The connection of theorem \ref{maint} is a characteristic connection
with torsion in $\bgw_3$. If the W\"unschmann ODE is general enough,
the torsion may be quite arbitrary within $\bgw_3$. From 
proposition \ref{p.purtor} we know that the independent components of 
the curvature of  a characteristic connection with $T\in\bgw_3$ are
$R$, $\der A^{(3)}$,  $\der A^{(7)}$ and $K^i$. In the notation of
theorem \ref{maint} they read:
\ben\der A^{(3)}=
  \bma 0 & 0 & 0 & a_1 & a_2 \\
       0 & 0 & -3a_1 & -2a_2 & a_3 \\
       0 & 3a_1 & 0 & -3a_3 & 0 \\
      -a_1 & 2a_2 & 3a_3 & 0 & 0 \\
      -a_2 & -a_3 & 0 & 0 & 0 \\ 
  \ema, 
\een
\ben\der A^{(7)}=
  \bma 0 & b_1 & b_2 & b_3 & b_4 \\
       -b_1 & 0 & 2b_3 & 8b_4 & b_5 \\
       -b_2 & -2b_3 & 0 & 2b_5 & b_6 \\
      -b_3 & -8b_4 & -2b_5 & 0 & b_7 \\
      -b_4 & -b_5 & -b_6 & -b_7 & 0 \\ 
  \ema, 
\een
\ben
K=\frac{\sqrt{3}}{3} \bma c_1 & c_2 & c_3 & c_4 & c_5 \ema^T,
\een
and, as we said above, the Ricci scalar is given by the function $R$.
The Ricci vector $\Rv^i=\Upsilon^{ijk}g_{jk}$ is as follows 
\ben
 \Rv^i=\frac{7}{6}\sqrt{3}\left( t_3^2,\ -\frac{1}{3}t_2t_3,\ \frac{1}{9}t_1t_3+\frac{2}{27}t_2^2,\  -\frac{1}{9}t_1t_2,\  \frac{1}{9}t^2_1 \right).
\een
The Ricci tensor satisfies the following equations
$$ R_{(ij)}=\tfrac{1}{5}Rg_{ij}+\tfrac{2}{7}\Rv^k\ten_{ijk},$$
$$\der A^{(3)}=4R^{(3)},\qquad\qquad \der A^{(7)}=\tfrac{2}{3}R^{(7)}.$$

Using theorem \ref{maint} we can also express the Ricci tensor $(Ric)^i_{~j}=g^{ik}R_{kj}$ in
terms of the endomorphisms $E_-$, $E_0$, $E_+$, $E_1$ of \eqref{om1}:
\begin{corollary}\label{cor.ricci}
The Ricci tensor of a characteristic connection with torsion in $\bgw_3$ has the following form in any adapted coframe
\begin{align*}
Ric=&\left(\frac {1}{54}t_2^{2}+\frac{1}{36}t_1t_3-\frac{1}{20}R \right)E_1
+\frac{1}{8}b_1E_-^{3}+\frac{1}{108}t_1^{2}E_-^{2}+ \\ \\
& +\left(-\frac {1}{54}t_1t_2+\frac{1}{8}a_1-\frac{1}{2}b_3\right)E_- 
+\frac {5}{16}b_4E_0^{3}+ \left(\frac {1}{108}t_2^{2}+\frac {1}{72}t_1t_3 \right)E_0^{2}+ \\ \\
&+ \left(-\frac {17}{4}b_4+\frac{1}{8}a_2 \right)E_0
-\frac{1}{8}b_7E_+^{3}+\frac{1}{12}t_3^{2}E_+^{2} 
+ \left(-\frac{1}{8}a_3+\frac{1}{2}b_5-\frac{1}{18}t_2t_3 \right)E_++ \\ \\
&-\frac{5}{32}b_5E_0E_+E_0+\frac{1}{8}b_6E_+E_0E_+
+\frac{1}{54}t_1t_2E_0E_-
+\frac {5}{32}b_3E_0E_-E_0+ \\ \\ &
+\frac{1}{8}b_2E_-E_0E_--\frac{1}{18}t_2t_3E_0E_+.
\end{align*}
\end{corollary}

Of course, since the geometry is constructed from an ODE determined by
the choice of $F=F(x,y,y_1,y_2,y_3,y_4)$, the coefficients 
$a_1,\ldots,a_3,b_1,\ldots,b_7,R$ are expressible in terms of $F$ and
its derivatives. 
Given the connection of theorem \ref{maint} we calculated the explicit formulae for these coefficients and obtained the following 
\begin{corollary}\label{cor.vanish}
A $\Gl$ geometry generated by a 5th order ODE satisfying W\"unschmann
conditions (\ref{wilki}) has the following properties.

The torsion $T$ vanishes iff
\ben F_{44}=0.\een

The 2-form $\der A^{(3)}$ vanishes iff
\begin{align*} 
&(\D F_4)_{34}-(\D F_3)_{44} -\tfrac{3}{5}(\D F_4)_4F_{44} -\tfrac{4}{5}\D F_4 F_{444}+\tfrac {6}{25}F^2_{44}F_4
+\tfrac {4}{25}F^2_4F_{444}+ \\
&+\tfrac{3}{10}F_{34}F_{44}-\tfrac{1}{5}F_4F_{344}
+\tfrac{3}{5}F_3F_{444}+F_{244}
-\tfrac{1}{2}F_{433}=0.
\end{align*}

The 2-form $\der A^{(7)}$ vanishes iff
\ben F_{444}=0. \een

The Ricci vector $\Rv$ is aligned with the vector $K$, i.e. $K=u\Rv$, $u\in\bbR$, iff
\ben
(\D F_4)_{44}-\tfrac{1}{2}F_{344}-\tfrac{2}{5}F_{4}F_{444}-\tfrac {8}{15}F^2_{44}+7uF^2_{44}=0.
\een
\end{corollary}
\noindent We skip writing the 
formula for the Ricci scalar since it is very complicated.

We now pass to the proof of theorem \ref{maint}. On doing this we
will apply a variant of the Cartan method of equivalence. This will be a
rather long and complicated procedure. Thus, for the clarity of the
presentation, we will divide the proof into three main steps, each of
which will occuppy its own respective section 
\ref{s.g-structure}, \ref{s.construction} and
\ref{s.purtor}. First, in section \ref{s.g-structure} we will 
prove lemma \ref{l.g-structure}, which assures that a
class of contact equivalent 5th order ODEs is a $G$-structure on a 4-order jet space
$J$. Thus, we will have a bundle $G\to J\times G \to J$, a reduction of the
frame bundle $F(J)$. In the second step, in section
\ref{s.construction}, we will use the Cartan method of
equivalence in order to construct a submanifold $P\subset J \times G$
together with a coframe on $P$ which fulfills the requirements of 
proposition \ref{p.constr}. This coframe, via proposition
\ref{p.constr}, will define an irreducible $\Gl$ structure for us and
simultanously will provide us with  
a $\gl$ connection on the space of solutions of the ODE. 
The obstructions for an ODE to possess this structure,  
W\"unschmann's expressions for $F$, will appear automatically in the
course of the construction. This part of considerations is summarized
in theorem \ref{mt}. The $\Gl$ structure obtained in this way will
turn out to be nearly integrable, but the connection constructed will
differ from the characteristic one. Therefore, in section
\ref{s.purtor}, we will construct the characteristic connection  
associated with the $\Gl$ structure obtained. This will have torsion
in $\bgw_3$. This construction is described by lemma \ref{prol}.

\subsection{5th order ODE modulo contact transformations}\label{s.g-structure}
 
Let us consider a general 5th order ODE \eqref{piec}. We define the following coframe
\begin{eqnarray}
&&\om^0=\der y-y_1\der x,\nonumber\\
&&\om^1=\der y_1-y_2\der x,\nonumber\\
&&\om^2=\der y_2-y_3\der x,\label{formy}\\
&&\om^3=\der y_3-y_4\der x,\nonumber\\
&&\om^4=\der y_4-F(x,y,y_1,y_2,y_3,y_4)\der x,\nonumber \\
&&\om_+=\der x \nonumber
\end{eqnarray}   
on $J$. We see that every solution of (\ref{piec}) is
a curve $c(x)=(x,y(x),y_1(x),y_2(x),y_3(x)$, $y_4(x))\subset J$ and the vector field $\D$ on $J$ has curves $c(x)$ as the integral curves. The 1-forms $(\om^0,\om^1,\om^2,\om^3,\om^4)$ annihilate $\D$ whereas $\D\lrcorner \om_+=1$. The 5-dimensional space $M^5$ of integral curves of $\D$ is clearly the space of solutions of
(\ref{piec}) and we have a fibration $\bbR\to J\to M^5$.  

Suppose now, that  equation \eqref{piec} undergoes a contact transformation \eqref{contact}, which brings it to $\bar{y}_5=\bar{F}(\bar{x},\bar{y},\bar{y}_1,\bar{y}_2,\bar{y}_3,\bar{y}_4)$.
Then the coframe transforms according to

\be\label{g-structure}
\bma \om^0 \\ \om^1 \\ \om^2 \\ \om^3 \\ \om^4 \\ \om_+  \ema 
\mapsto
\bma \bar{\om}^0 \\ \bar{\om}^1 \\ \bar{\om}^2 \\ \bar{\om}^3 \\ \bar{\om}^4 \\ \bar{\om}_+  \ema 
=
\bma \al^0_{~0} & 0 & 0 & 0 & 0 & 0  \\
     \al^1_{~0} & \al^1_{~1} & 0 & 0 & 0 & 0 \\
     \al^2_{~0} & \al^2_{~1} & \al^2_{~2} & 0 & 0 & 0 \\
     \al^3_{~0} & \al^3_{~1} & \al^3_{~2} & \al^3_{~3} & 0 & 0 \\
     \al^4_{~0} & \al^4_{~1} & \al^4_{~2} & \al^4_{~3} & \al^4_{~4} & 0 \\
     \al^5_{~0} & \al^5_{~1} & 0 & 0 & 0 & \al^5_{~5}
\ema
\bma \om^0 \\ \om^1 \\ \om^2 \\ \om^3 \\ \om^4 \\ \om_+  \ema.
\ee 
Here $\al^i_{~j}$, $i,j=0,1,2,3,4,5$, are real functions on $J$
defined by the formulae \eqref{contact}. They satisfy the
nondegeneracy condition $$\al^0_{~0}\al^1_{~1}\al^2_{~2}\al^3_{~3}\al^4_{~4}\al^5_{~5}\neq 0.$$
The transformed coframe encodes all the contact invariant information about the ODE. In particular, it preserves the simple ideal $(\om^0,\ldots,\om^4)$, from which we can recover solutions of the transformed equation. Hence we have

\begin{lemma}\label{l.g-structure}
A 5th order ODE $y_5=F(x,y,y_1,y_2,y_3,y_4)$ considered modulo contact transformations of variables  is a $G$-structure on the 4-jet space $J$, such that the coframe $(\om^0,\om^1,\om^2,\om^3,\om^4,\om_+)$ of \eqref{formy} belongs to it and the group $G$ is given by the matrix in \eqref{g-structure}
\end{lemma}

\subsection{$\glg(2,\bbR)$ bundle over space of solutions}\label{s.construction}
Using the Cartan method we explicitly construct a submanifold $P\subset J\times G$ and a coframe $(\theta^0,\theta^1,\theta^2,\theta^3,\theta^4,\Om_-,\Om_+,\Om_0,\Om_1)$ on $P$ satisfying proposition \ref{p.constr}.
This part of the proof is divided into eight steps.

\subsubsection*{{\bf 1)}} We observe that there is a natural choice for the forms $(\theta^0,\theta^1,\theta^2,\theta^3,\theta^4)$ of the coframe. Since we are going to build a $\Gl$ structure on the space of solutions $P$ must be a bundle over $M^5$ and the forms $(\theta^0,\theta^1,\theta^2,\theta^3,\theta^4)$ must annihilate vectors tangent to leaves of the projection $P\to M^5$. But on $J\times G$ there are six distinguished 1-forms given by
\be\label{canonical}
\left(\begin{aligned} \theta^0 \\ \theta^1 \\ \theta^2 \\ \theta^3 \\ \theta^4 \\ \theta_+ \end{aligned}\right) =
\left(\begin{aligned}
 &\al^0_{~0}\om^0  \\
 &\al^1_{~0}\om^0+\al^1_{~1}\om^1 \\
 &\al^2_{~0}\om^0+\al^2_{~1}\om^1+\al^2_{~2}\om^2 \\
 &\al^3_{~0}\om^0+\al^3_{~1}\om^1+\al^3_{~2}\om^2+\al^3_{~3}\om^3 \\ 
 &\al^4_{~0}\om^0+\al^4_{~1}\om^1+\al^4_{~2}\om^2+\al^4_{~3}\om^3+\al^4_{~4}\om^4 \\
 &\al^5_{~0}\om_0+\al^5_{~1}\om^1+\al^5_{~5}\om_+
\end{aligned}\right).
\ee
These forms are the components of the 
canonical $\bbR^6$ valued 1-form on $J\times G$. 
Five among these forms, $\theta^0,\theta^1,\theta^2,\theta^3,\theta^4$
also annihilate vectors tangent to the projection $J\times G\to
M^5$. We choose them to be the members of the sought coframe  $(\theta^0,\theta^1,\theta^2,\theta^3,\theta^4,\Om_-,\Om_+,\Om_0,\Om_1)$. Now we must construct a 9-dimensional submanifold $P$ on which $\theta^i$ satisfy equations \eqref{sysc} with some linearly independent forms $\Om_-,\Om_+,\Om_0,\Om_1$.

\subsubsection*{{\bf 2)}} We calculate $\der\theta^0$ and get
$$\der\theta^0=\left(\frac{\der\al^0_{~0}}{\al^0_{~0}} -\frac{\al^1_{~0}}{\al^1_{~1}\al^5_{~5}}\theta_+\right)\dz\theta^0+\frac{\al^0_{~0}}{\al^1_{~1}\al^5_{~5}}\theta_+\dz\theta^1 -\frac{\al^5_{~0}}{\al^1_{~1}\al^5_{~5}}\theta^0\dz\theta^1 $$
For this equation to match \eqref{sysc} we define
\begin{align}
&\Om_+=\theta_+ \label{con1} \\
&4(\Om_1+\Om_0)=\frac{\der\al^0_{~0}}{\al^0_{~0}} -\frac{\al^1_{~0}}{\al^1_{~1}\al^5_{~5}}\theta_+ \mod \theta^i, \label{con2}
\end{align}
with yet unspecified $\theta^i$ terms in \eqref{con2}, and set
\be\label{al1}\al^0_{~0}=-4\al^1_{~1}\al^5_{~5}\ee
to get $-4$ coefficient in the $\Om_+\dz\theta^1$ term. Thereby
$$\der\theta^0=4(\Om_1+\Om_0)\dz\theta^0-4\Om_+\dz\theta^1 \mod
\theta^i\dz\theta^j$$ on the 23-dimensional subbundle of $J\times G\to
M^5$ given by \eqref{al1}. We see that the form $\theta_+$
plays naturally the role of the connection 1-form $\Om_+$.

\subsubsection*{\bf 3)} We calculate $\der\theta^1$ on the 23-dimensional bundle. In order to get $$\der\theta^1=-\Om_-\dz\theta^0+(4\Om_1+2\Om_0)\dz\theta^1-3\Om_+\dz\theta^2\mod\theta^i\dz\theta^j$$ 
we set 
\begin{align}
4\Om_1&+2\Om_0=\frac{\der\al^1_{~1}}{\al^1_{~1}}+\frac{\al^1_{~0}\al^2_{~2}-\al^1_{~1}\al^2_{~1}}{\al^1_{~1}\al^2_{~2}\al^5_{~5}}\theta_+ \mod \theta^i, \label{con3} \\
\Om_-&=-\frac{\der\al^1_{~0}}{4\al^1_{~1}\al^5_{~5}}+\frac{\al^1_{~0}\der\al^1_{~1}}{4(\al^1_{~1})^2\al^5_{~5}}
\label{con4}\\
&+\frac{(\al^1_{~0})^2\al^2_{~2}+(\al^1_{~1})^2\al^2_{~0}-\al^1_{~1}\al^2_{~1}\al^1_{~0})^2}{4(\al^1_{~1})^2\al^2_{~2}(\al^5_{~5})^2}\theta_+\mod\theta^i,  \notag
\end{align}
and
\be\label{al2} \al^2_{~2}=-\frac{\al^1_{~1}}{3\al^5_{~5}} \ee
obtaining a 22-dimensional subbundle  of $J\times G\to M^5$ on which
$\der\theta^0$ and $\der\theta^1$ are in the desired form.

\subsubsection*{\bf 4)} At this point all four connection 1-forms
$\Om_-,\Om_+,\Om_0,\Om_1$ are fixed up to the $\theta^i$ terms. They
are determined by the equations \eqref{con1}, \eqref{con2},
\eqref{con3}, \eqref{con4}. Thus we can not introduce any new 1-forms
to bring $\der\theta^2$ into the desired form. Now to get
$\der\theta^2$ in the form as in theorem \ref{maint}, we may only use the
yet unspecified coefficients $\alpha$s. That is why $\der\theta^2$ imposes more conditions on $\alpha$s. It follows that for
$\der\theta^0,\,\der\theta^1$ and $\der\theta^2$ to be of the form \eqref{sysc} the subbundle $P$ must satisfy
\begin{eqnarray}
&&\al^0_{~0}=-4\al^1_{~1}\al^5_{~5},\nonumber\\
&&\al^2_{~0}=\frac{-75 (\al^1_{~0})^2+(\al^1_{~1})^2(-20\D F_4+20F_3+7F_4^2)}{300\al^1_{~1}\al^5_{~5}},\nonumber\\
&&\al^2_{~1}=\frac{-15\al^1_{~0}+\al^1_{~1}F_4}{30\al^5_{~5}},\nonumber\\
&&\nonumber\\
&&\al^2_{~2}=-\frac{\al^1_{~1}}{3\al^5_{~5}},\nonumber\\
&&\al^3_{~0}=[1800(\al^1_{~1}\al^5_{~5})^2]^{-1}\times\nonumber\\
&&[1125(\al^1_{~0})^3+45\al^1_{~0}(\al^1_{~1})^2(20\D F_4-20F_3-7F_4^2)+\nonumber\\
&&2(\al^1_{~1})^3(100\D^2F_4-200F_2-30F_4 \D F_4-60F_3F_4-11F_4^3)],\nonumber\\
&&\nonumber\\
&&\al^3_{~1}=\frac{225 (\al^1_{~0})^2-30\al^1_{~0}\al^1_{~1}F_4+(\al^1_{~1})^2(80\D F_4-100F_3-31F_4^2)}{1200\al^1_{~1}(\al^5_{~5})^2}\label{alfas},\\
&&\al^3_{~2}=\frac{5\al^1_{~0}-\al^1_{~1}F_4}{20 (\al^5_{~5})^2},\nonumber\\
&&\al^3_{~3}=\frac{\al^1_{~1}}{6(\al^5_{~5})^2},\nonumber\\
&&\nonumber\\
&&\al^4_{~1}=[18000(\al^1_{~1})^2(\al^5_{~5})^3]^{-1}\times\nonumber\\
&&[-1125(\al^1_{~0})^3+225(\al^1_{~0})^2\al^1_{~1}F_4-15\al^1_{~0}
(\al^1_{~1})^2(80\D F_4-100F_3-31F_4^2)+\nonumber\\
&&(\al^1_{~1})^3(-400\D^2F_4+1400F_2+240F_4 \D F_4+180F_3F_4+11F_4^3)],\nonumber\\
&&\nonumber\\
&&\al^4_{~2}=\frac{-75(\al^1_{~0})^2+30\al^1_{~0}\al^1_{~1}F_4+(\al^1_{~1})^2(-40\D F_4+80F_3+17F_4^2)}{600\al^1_{~1}(\al^5_{~5})^3},\nonumber\\
&&\al^4_{~3}=\frac{-5\al^1_{~0}+3\al^1_{~1}F_4}{30(\al^5_{~5})^3},\nonumber\\
&&\al^4_{~4}=-\frac{\al^1_{~1}}{6(\al^5_{~5})^3}.\nonumber
\end{eqnarray}
 
The necessity of these conditions can be checked by a direct, quite 
lengthy calculations. We performed these calculations using the
symbolic computation programs Maple and Mathematica. 

We stress that conditions (\ref{alfas}) are only \emph{necessary}
for $\der\theta^2$ to satisfy (\ref{sysc}). It is because certain
unwanted terms cannot be removed by \emph{any} choice of subbundle
$P$. Vanishing of these unwanted terms is a property of the ODE
itself, and this is the reason for the W\"unschmann conditions to
appear. 

More specifically, to achieve 
$$ \der\theta^2=-2\Om_-\dz\theta^1+4\Om_1\dz\theta^2-2\Om_+\dz\theta^3\mod\theta^i\dz\theta^j$$ 
on the bundle defined by \eqref{al1}, \eqref{al2} and \eqref{alfas}
an ODE must satisfy 
\be
50\D^2F_4 - 75 \D F_3 + 50 F_2 - 60 F_4\D F_4 + 30 F_3F_4 + 8 F_4^3=0.\label{wilk1}
\ee  
It follows from the construction that this condition, the first 
of \eqref{wilki}, is invariant under the contact transformation of variables. 

From now on we restrict our considerations only to contact equivalence
class of ODEs satisfying (\ref{wilk1}). If (\ref{alfas}) and (\ref{wilk1}) are satisfied then the three 
differentials $\der\theta^0$, $\der\theta^1$ and $\der\theta^2$ are
precisely in the form (\ref{sysc}).

\subsubsection*{{\bf 5)}} The requirement that also 
$\der\theta^3$ is in the form (\ref{sysc}) is equivalent to the
following equation for $\al^4_{~0}$:
\begin{eqnarray}
&&\al^4_{~0}=[120000(\al^1_{~1}\al^5_{~5})^3]^{-1}\times\nonumber\\
&&[-1875(\al^1_{~0})^4-150(\al^1_{~0}\al^1_{~1})^2(20\D F_4-20F_3-7F_4^2)-\nonumber\\
&&40\al^1_{~0}(\al^1_{~1})^3(50\D F_3-100F_2+30F_4
    \D F_4-40F_3F_4-9F_4^3)+\label{al4}\\
&&(\al^1_{~1})^4\Big(400(-5\D^2F_3+10\D F_2-6(\D F_4)^2+10F_3\D F_4-3F_3^2+F_4\D F_3)+
\nonumber\\
&&120F_4^2(7\D F_4-5F_3)-63F_4^4\Big)].\nonumber
\end{eqnarray}   

\subsubsection*{{\bf 6)}} If condition \eqref{al4} is also imposed we have
$$(\der\theta^4+4\Om_-\dz\theta^3-4(\Om-\Om_0)\dz\theta^4)\dz\theta^0\dz\theta^1=0 \mod \theta^i.$$
However, 
$$\der\theta^4\dz\theta^0\dz\theta^2\dz\theta^3\dz\theta^4=0$$ 
if and only if second condition of \eqref{wilki} is satisfied:
\begin{eqnarray}
&&375 \D^2F_3 - 1000 \D F_2 + 350 \D F_4^2 + 1250 F_1 - 650F_3 \D F_4  + 
    200 F_3^2 -\nonumber\\&& 150 F_4 \D F_3 + 200 F_2 F_4 - 140 F_4^2 \D F_4 + 
    130 F_3 F_4^2 + 14 F_4^4=0.
\label{wilk2}
\end{eqnarray}
Again it follows from the construction that condition (\ref{wilk2}),
considered simultaneously with (\ref{wilk1}), is invariant under
contact transformations of the variables. From now on, we assume that
all our 5th order ODEs (\ref{piec}) satisfy both conditions
(\ref{wilk1}), (\ref{wilk2}). It follows that it is still not
sufficient to force $\der\theta^4$ to satisfy the system
(\ref{sysc}), since without further assumptions on $F$, we do \emph{not}
have $\der\theta^4\dz\theta^1\dz\theta^2\dz\theta^3\dz\theta^4=0$.  
To achieve this it is necessary and sufficient to impose the last
restriction on $F$:
\begin{eqnarray}
&&1250 \D^2F_2 - 6250 \D F_1 + 1750 \D F_3\D F_4 - 2750 F_2\D F_4  -\nonumber\\&& 875  F_3\D F_3
  +1250 F_2 F_3 - 500F_4 \D F_2  + 700 (\D F_4)^2 F_4 +\label{wilk3}\\
&& 1250 F_1 F_4
  - 
    1050 F_3 F_4\D F_4  + 350 F_3^2 F_4 - 350 F_4^2\D F_3  + \nonumber\\
&& 550 F_2 F_4^2
  -
    280 F_4^3 \D F_4 + 210 F_3 F_4^3 + 28 F_4^5 + 18750 F_y=0.\nonumber
\end{eqnarray}

\subsubsection*{{\bf 7)}} Assuming that $F$ satisfies conditions \eqref{wilki} and fixing coefficients $\al^i_{~j}$
according to (\ref{alfas}), (\ref{al4}) we are remained with a 11-dimensional subbundle of $J\times G\to M^5$
parametrized by $(x,y,y_1,y_2,y_3,y_4,\al^1_{~0},\al^1_{~1},\al^5_{~5},\al^5_{~0},\al^5_{~1})$. It follows that the forms $\Om_0,\Om_1,\Om_-,\Om_+$ on this bundle are
\begin{eqnarray}
&&\Om_+=\theta_+,\nonumber\\
&&\nonumber\\
&&\Om_0=\frac{\der
    \al^5_{~5}}{2\al^5_{~5}}-\frac{5\al^1_{0}+\al^1_{~1}F_4}{20\al^1_{~1}\al^5_{~5}}\theta_+ \mod\theta^i,\nonumber\\
&&\nonumber\\
&&\Om_1=\frac{\der\al^1_{~1}}{4\al^1_{~1}}-\frac{\der\al^5_{~5}}{4\al^5_{~5}}+\frac{F_4}{20\al^5_{~5}}\theta_+\mod\theta^i,\label{dt}\\
&&\nonumber\\
&&\Om_-=\frac{\der\al^1_{~0}}{4\al^1_{~1}\al^5_{~5}}-\frac{\al^1_{~0}\der\al^1_{~1}}{4(\al^1_{~1})^2\al^5_{~5}}-\nonumber\\
&&\frac{25(\al^1_{~0})^2+10\al^1_{~0}\al^1_{~1}F_4+(\al^1_{~1})^2(20\D F_4-20F_3-7F_4^2)}{400(\al^1_{~1}\al^5_{~5})^2}\theta_+\mod\theta^i.\nonumber
\end{eqnarray}

\subsubsection*{{\bf 8)}} In order to construct a 9-dimensional bundle and find the $\theta^i$ terms
in (\ref{dt}) we need to consider the $\der\Om_A$ part of equations (\ref{sysc}). Forcing
$\der\Om_A$ not to have $\Om_A\dz\theta^i$ terms we \emph{uniquely} specify the
$\theta^i$ terms in (\ref{dt}). This requirement, in particular, fixes the coefficients
$\al^5_{~1}$ and $\al^5_{~0}$ to be:
\begin{eqnarray}
&&\al^5_{~1} = \frac{\al^5_{~5}(10\D F_{44} + 5 F_{34} + 6 F_4 F_{44})}{50},\nonumber\\
&&\al^5_{~0}=\frac{\al^5_{~5}}{250}~[50(\D F_{34}+7F_{24}-5F_{33})+\label{al5}\\
&&5F_4(6\D F_{44}-37F_{34})+2F_{44}(-60\D F_4+145F_3+21F_4^2)].\nonumber
\end{eqnarray}  
Now all the forms
$(\theta^0,\theta^1,\theta^2,\theta^3,\theta^4,\Om_+,\Om_-,\Om_0,\Om_1)$
are well defined and independent on a 9-dimensional manifold $P$ 
parametrized by
$(y,y_1,y_2,y_3,y_4,x,\al^1_{~0},\al^1_{~1},\al^5_{~5})$. We calculate
structural equations \eqref{sysc} for these forms and have the
following 

\begin{theorem}\label{mt}
A 5th order ODE $y^{(5)}=F(x,y,y',y'',y^{(3)},y^{(4)})$ considered modulo contact transformation of variables has an irreducible $\glg(2,\bbR)$ structure on the space of its solution $M^5$ together with a $\gla(2,\bbR)$ connection $\Gamma$ if and only if its defining function $F=F(x,y,y_1,y_2,y_3,y_4)$ satisfies the contact invariant W\"unschmann conditions (\ref{wilki}).
The bundle $\glg(2,\bbR)\to P\to M^5$ is given by the equations (\ref{alfas}), (\ref{al4}) and (\ref{al5}). 
The first structural equations for the connection $\Gamma=(\Om_+, \Om_-, \Om_0, \Om_1)$ on $P$ read
\begin{eqnarray}
\der\theta^0&=&4(\Om_1+\Om_0)\dz\theta^0-4\Om_+\dz\theta^1+\nonumber\\
&&t_1\theta^0\dz\theta^1+t_2\theta^0\dz\theta^2+t_3\theta^0\dz\theta^3,\nonumber\\
\der\theta^1&=&-\Om_-\dz\theta^0+(4\Om_1+2\Om_0)\dz\theta^1-3\Om_+\dz\theta^2+\nonumber\\
&&\tfrac12t_1\theta^0\dz\theta^2+\tfrac13t_2\theta^0\dz\theta^3+\tfrac14t_3\theta^0\dz\theta^4+\nonumber\\
&&t_2\theta^1\dz\theta^2+t_3\theta^1\dz\theta^3,\nonumber\\
\der\theta^2&=&-2\Om_-\dz\theta^1+4\Om_1\dz\theta^2-2\Om_+\dz\theta^3+\label{mt1}\\
&&\tfrac29t_1\theta^0\dz\theta^3+\tfrac{1}{18}t_2\theta^0\dz\theta^4+\tfrac13t_1\theta^1\dz\theta^2+\nonumber\\
&&\tfrac89t_2\theta^1\dz\theta^3+\tfrac23t_3\theta^1\dz\theta^4+t_3\theta^2\dz\theta^3,\nonumber\\
\der\theta^3&=&-3\Om_-\dz\theta^2+(4\Om_1-2\Om_0)\dz\theta^3-\Om_+\dz\theta^4+\nonumber\\
&&\tfrac{1}{12}t_1\theta^0\dz\theta^4+\tfrac{1}{3}t_1\theta^1\dz\theta^3+\tfrac13t_2\theta^1\dz\theta^4+\nonumber\\
&&t_2\theta^2\dz\theta^3+\tfrac32t_3\theta^2\dz\theta^4,\nonumber\\
\der\theta^4&=&-4\Om_-\dz\theta^3+4(\Om_1-\Om_0)\dz\theta^4+\nonumber\\
&&\tfrac{1}{3}t_1\theta^1\dz\theta^4+t_2\theta^2\dz\theta^4+3t_3\theta^3\dz\theta^4,\nonumber
\end{eqnarray}
with the torsion coefficients
\begin{eqnarray}
t_3&=&\frac{6 (\al^5_{~5})^2}{5\al^1_{~1}}F_{44},\nonumber\\
t_2&=&\frac{9\al^5_{~5}}{50(\al^1_{~1})^2}~[\al^1_{~1}(10\D F_{44}+3F_4F_{44})+5\al^1_{~0}F_{44}],\nonumber\\
&&\nonumber\\
t_1&=&[1000(\al^1_{~1})^3]^{-1}\times\nonumber\\
&&\Big(225(\al^1_{~0})^2F_{44}+90\al^1_{~0}\al^1_{~1}(10\D F_{44}+3F_4F_{44})+\nonumber\\ 
&&-9(\al^1_{~1})^2[20(5\D F_{34}+20F_{24}-15F_{33}+3F_4\D F_{44}-11F_4F_{34})+\nonumber\\
&&+F_{44}(-120\D F_4+340F_3+51F_4^2)]\Big).\nonumber
\end{eqnarray}
\end{theorem}
Also the second structural equations are easily calculable but we skip them due to their complexity.

It is remarkable that the above $\gla(2,\bbR)$ connection has torsion with not more than \emph{three}
functionally independent coefficients $t_1, t_2, t_3$. This
suggests that the $\glg(2,\bbR)$ geometry on the 5-dimensional
solution space $M^5$ of the ODE is \emph{nearly integrable} with
torsion in the irreducible part $\bgw_3$ only. 
That it is really the case will be shown below.

\subsection{Characteristic connection with torsion in $\bgw_3$}\label{s.purtor}
As we know from section \ref{trzy}, given an irreducible
$\glg(2,\bbR)$-structure $(M^5,[g,\ten,A])$, we can ask if such a
structure is nearly integrable. According to propositions
\ref{chaco} and \ref{chacon}, the necessary and sufficient condition for nearly integrability
is that the structure admits a $\gla(2,\bbR)$-valued connection with
\emph{totally skew symmetric} torsion.

In our case of ODEs satisfying W\"unschmann conditions we have a $\gla(2,\bbR)$-valued connection of theorem \ref{mt}, whose torsion is expressible in terms of three independent functions. This torsion, however, has quite complicated 
algebraic structure, in particular it is \emph{not} totally skew symmetric. 

It appears that an irreducible $\glg(2,\bbR)$ structure $(M^5,[g,\ten,A])$ 
associated with any 5th order ODE satisfying conditions \eqref{wilki} 
admits another $\gla(2,\bbR)$-valued connection that
\emph{has} totally skew symmetric torsion. Thus all structures 
$(M^5,[g,\ten,A])$ originating from W\"unschmann 5th order ODEs are nearly
integrable; the new connection is their characteristic connection. Even more
interesting is the fact that its torsion is still more
special: it is always in $\bgw_3$. 

One way of seeing this is to calculate the Weyl connection $\we$ for
the corresponding
$(M^5,[g,\ten,A])$ and to decompose it according to (\ref{unic}). Here
we prefer another method -- the analysis in terms of the Cartan bundle $P$ of
theorem \ref{mt}.

\begin{lemma}\label{prol}
Consider a contact equivalence class
of 5th order ODEs satisfying conditions \eqref{wilki}. 
Let
$\theta^0,\theta^1,\theta^2,\theta^3,\theta^4,\Om_+,\Om_-,\Om_0,\Om_1$
and $t_1,t_2,t_3$ be the objects of theorem \ref{mt}. Then there is a
$\gla(2,\bbR)$ connection $\tilde{\Gamma}=(\Pi_+,\Pi_-,\Pi_0,\Pi_1)$
whose torsion $\tilde{T}^i_{~jk}$ is totally skew symmetric and has
its associated 3-form in $\tilde{T}\in*\bgw_3$. Explicitly:
\begin{eqnarray*}
\tilde{T}&=&\tfrac{1}{12}t_1(-\theta^0\dz\theta^1\dz\theta^4+2\theta^0\dz\theta^2\dz\theta^3)+\\
 &&\tfrac{1}{12}t_2(-\theta^0\dz\theta^2\dz\theta^4+8\theta^1\dz\theta^2\dz\theta^3)+\\
 &&\tfrac14t_3(-\theta^0\dz\theta^3\dz\theta^4+2\theta^1\dz\theta^2\dz\theta^4).
\end{eqnarray*} 
\begin{proof}
Any $\gla(2,\bbR)$ connection $\tilde{\Gamma}=(\Pi_+,\Pi_-,\Pi_0,\Pi_1)$ compatible with the $\glg(2,\bbR)$ structure of theorem \ref{mt} is given by 
\begin{align}
 &\Pi_A=\Om_A+\sum_i\gamma_{Ai} \theta^i,  \quad A\in\{+,0,-\},\quad i=0,\ldots,4, \label{wisc} \\
 &\Pi_1=\Om_1 \notag
\end{align}
with arbitrary functions $\gamma_{Ai}$. We calculate structural equations $\der\theta+\tilde{\Gamma}\dz\theta=\tilde{T}$ for $\tilde{\Gamma}$
utilising equations (\ref{mt1}), and ask if there
exists a choice of $\gamma_{Ai}$ such that the new torsion 
$\tilde{T}^i_{~jk}$ satisfies $g_{il}\tilde{T}^l_{~jk}=\tilde{T}_{[ijk]}$
and $\tilde{T}=\tfrac16
g_{il}\tilde{T}^l_{~jk}\theta^i\dz\theta^j\dz\theta^k\in*\bgw_3$.
Using lemma \ref{lem} we easily find that the unique solution is given by 
\begin{eqnarray*}
 &&\Pi_+ = \Om_+ - \tfrac16t_1\theta^0 - \tfrac13t_2\theta^1 - \tfrac12t_3\theta^2, \\
 &&\Pi_- = \Om_- + \tfrac16t_1\theta^2 + \tfrac13t_2\theta^3 +\tfrac12t_3\theta^4, \\
 &&\Pi_0 = \Om_0 - \tfrac16t_1\theta^1 - \tfrac13t_2\theta^2 -\tfrac12t_3\theta^3, \\
 &&\Pi_1=\Om_1, 
 \end{eqnarray*}
\end{proof}
\end{lemma}

Lemma \ref{prol} together with the results of section \ref{nabu} prove theorem \ref{maint}.

\begin{remark}\label{r.loss}
Note that a passage from $\Om_+$ to $$\Pi_+ = \Om_+ - \tfrac16t_1\theta^0 - \tfrac13t_2\theta^1 - \tfrac12t_3\theta^2$$ belongs to a \emph{larger} class of transformations than the contact transformations (\ref{g-structure}), (\ref{canonical}); it involves a forbidden $\theta^2$ term. Thus it may happen that there are nonequivalent classes of ODEs which define the same $(M^5,[g,\ten,A])$. To
distinguish between nonequivalent ODEs one has to use the connection of theorem \ref{mt}.
\end{remark}

\section{Examples of nearly integrable $\glg(2,\bbR)$ structures from 5th order ODEs}\label{s.examples}

In this section we provide examples of W\"unschmann ODEs and nearly
integrable $\glg(2,\bbR)$ structures related to them. Since such
structures have the torsions of their characteristic connections in
$\bgw_3$, then via proposition \ref{p.purtor}, they are characterized
by the torsion $T$, the Ricci scalar $R$, the components of Maxwell
2-forms $\der A^{(3)}$, $\der A^{(7)}$, and the vector $K$; all
these objects being associated to the
characteristic connection $\Gamma$. There is also the unique Weyl
connection $\we$ associated with these structures.

\subsection{Torsionfree structures}\label{s.notorsion}
We see from corollary \ref{cor.vanish} that $$ T\equiv 0 \quad \iff \quad  F_{44}\equiv0.$$  Then $\we=\Gamma$ and all the curvature components but the Ricci scalar necessarily vanish. The following proposition can be checked by direct calculation.

\begin{proposition}\label{p.notorsion}
The three nonequivalent differential equations  
$$y^{(5)}=c\Big(\frac{5y^{(3)3}( 5 - 27cy''^{2})}{9( 1+c y''^{2})^2} +
10\frac{y''y^{(3)}y^{(4)}}{ 1+c y''^2}\Big),$$
with $c=+1, 0, -1,$ represent the only three contact
nonequivalent classes of 5th order ODEs having the
corresponding nearly integrable $\glg(2,\bbR)$ structures 
$(M^5,[g,\ten,$ $A])$ with the characteristic
connection with vanishing torsion. In all three cases the holonomy of the 
Weyl connection $\we$ of structures $(M^5,[g,\ten,A])$ is 
reduced to the $\glg(2,\bbR)$. For all the three cases the Maxwell
2-form $\der A\equiv 0$. The corresponding Weyl structure is flat
for $c=0$. If $c=\pm1$, then in the conformal class $[g]$ there is an
Einstein metric of positive ($c=+1$) or negative ($c=-1$) Ricci
scalar. In case $c=1$ the manifold $M^5$ can be identified with
the homogeneous space $\sug(1,2)/\slg(2,\bbR)$ with an Einstein 
$g$ descending from
the Killing form on $\sug(1,2)$. Similarly in $c=-1$ case the 
manifold $M^5$ can be identified with
the homogeneous space $\slg(3,\bbR)/\slg(2,\bbR)$ with an Einstein 
$g$ descending from
the Killing form on $\slg(3,\bbR)$. In both cases with $c\neq 0$ the metric
$g$ is not conformally flat.
\end{proposition}

\subsection{Structures with vanishing Maxwell form}
From now on we assume that 
$$F_{44}\neq 0$$
and consider structures with vanishing Maxwell 2-form $\der A=0$. For such structures both torsion and curvature have at most 9 independent coefficients contained in $T$, $K$ and the scalar $R$. The simplest geometries in this class are those satisfying the additional equality 
$$ K^i=u\Rv^i,\qquad u\in\bbR.$$

Putting $\der A=0$ and $ K^i=u\Rv^i$ into structural equations of
theorem \ref{maint} and using Bianchi identities we find that either
$$ u=-\frac{1}{420},\qquad R=\frac{35}{54}(t_2^2-3t_1t_3)$$
or
$$u=\frac{2}{105},\qquad R=\frac{10}{27}(t_2^2-3t_1t_3).$$
Thus in these cases $R$ is funtionally dependent on $t_1,t_2,t_3$ and the only invariants for such $\Gl$ structures are $u$ and the sign of $R$. For each possible values of $u$ and $\sgn R$ we found a generating ODE.

\begin{proposition}
Consider the equations
\be F=\frac{5y_4^2}{3y_3}+\epsilon y_3^{5/3},\quad \epsilon=-1,0,1, \label{ex53} \ee
\be F=\frac{5y_4^2}{4y_3}, \label{ex54} \ee
and
\be F=\frac{5(8y_3^3-12y_2y_3y_4+3y_1y_4^2)}{6(2y_1y_3-3y_2^2)}, \label{exfrac} \ee
where the sign of expression $(2y_1y_3-3y_2^2)$ is an invariant, and the singular locus $2y_1y_3-3y_2^2=0$ separates nonequivalent equations with $\pm$ signs. The equations generate all the six $\Gl$ strucutres satisfying $\der A=0$ and $K^i=u\Rv^i$, $u\in\bbR$. 
\begin{align*}
&\text{For (\ref{ex53})} \qquad u=-\tfrac{1}{420}\quad \text{and}\quad  \sgn R=\epsilon, \\
&\text{for (\ref{ex54})} \qquad  u=\tfrac{2}{105} \quad \text{and}\quad  R=0, \\
&\text{for (\ref{exfrac})}\qquad u=\tfrac{2}{105} \quad \text{and} \quad \sgn R=\sgn (3y_2^2-2y_1y_3).
\end{align*}
\end{proposition}

Morover, the above ODEs can be also described in a geometric way by means of the symmetry group. 
\begin{proposition}
The equations (\ref{ex53}), (\ref{ex54}) and (\ref{exfrac}) are the
only 5th order W\"unschmann ODEs satisfying $F_{44}\neq 0$, $F_{444}=0$
and possessing the maximal group of transitive contact symmetries of dimension grater than five.
Equations $F=\frac{5y_4^2}{3y_3}$ and $F=\frac{5y_4^2}{4y_3}$ have 7-dimensional groups of symmetries, all the remaining have 6-dimensional ones.
\begin{proof}
The proof is based on further application of the Cartan method of equivalence. Let us return to the coframe of theorem \ref{mt}, which encodes all the contact invariant information about the ODE. If there are any nonconstant coefficients in the structural equations for this coframe we can use them for further reduction of the group $\Gl$ and of the bundle $P$.
For an ODE satisfying $F_{44}\neq 0$ we normalize $t_3=1$, $t_2=0$, which implies
$$\al^1_{~1}=\frac{6}{5}(\al^5_{~5})^2F_{44},\qquad \al^1_{~0}=-\frac{6}{25}(\al^5_{~5})^2(10\D F_{44}+3F_4F_{44}).$$

Now the coframe of theorem \ref{mt} is reduced to a 7-dimensional manifold $P_7$ parameterized by $(x,y,y_1,y_2,y_3,y_4,\al^5_{~5})$, three 1-forms $(\Om_0,\Om_-,\Om_1)$ become dependent on each other 
and we can use only one of them, our choice is $\Om_0$, to supplement $(\theta^0,\theta^1,\theta^2,\theta^3,\theta^4,\Om_+)$ to an invariant coframe on $P_7$.
Next we calculate structural equations for the new coframe. The coefficients in these equations are built from $\al^5_{~5}$ and 16 functions $f_1,\ldots,f_{16}$ of $x$, $y$, $y_1$, $\ldots$, $y_4$. In particular
\begin{align*}
\der\theta^0=\,&6\Om_0\dz\theta^0-4\Om_+\dz\theta^1+\frac{f_1}{(\al^5_{~5})^2}\theta^0\dz\theta^1+ \\
&+\frac{f_2}{\al^5_{~5}}\theta^0\dz\theta^2+f_3\theta^0\dz\theta^3+f_4\al^5_{~5}\theta^0\dz\theta^4,
\end{align*}
where for example
\ben
 f_3=-\frac{5F_{344}F_{44}+10 DF_{44}F_{444}+6F_4F_{44}F_{444}}{F_{44}^3}, \quad\quad f_4=5\frac{F_{444}}{F_{44}^2}.
\een

Let us assume $F_{444}=0$ and consider two possibilities: $f_3\neq const$ and $f_3=const$. If $f_3\neq const$ then it follows from the equations $\der^2\theta^i=0$, $\der^2\Om_A=0$ that $f_2$ may not be a constant. Thus $f_2/\al^5_{~5}$ and $f_3$ are \emph{two} functionally independent coefficients in structural equations for the 7-dimensional coframe $(\theta^0,\theta^1,\theta^2,\theta^3,\theta^4,\Om_+,\Om_0)$. According to the procedure of finding symmetries of ODEs, which is described in \cite{olver}, the dimension of the group of contact symmetries of a corresponding 5-order ODE is not larger than the dimension of the coframe minus the number of the independent coefficients in the structural equations, that is $7-2=5$. It follows that ODEs possessing contact symmetry group greater than 5-dimensional necessarily satisfy $f_3=const$. Let us assume $f_3=const$ then and we get from identities $\der^2\theta^i=0$, $\der^2\Om_A=0$ that (i) either $f_3=2$ or $f_3=\tfrac{3}{2}$ and (ii) for both admissible values of $f_3$ all the remaining nonvanishing functions $f_j$ are expressible by $f_1$. For example, the system corresponding to $f_3=\tfrac{3}{2}$ is the following
\begin{eqnarray*}
\der\theta^0&=&6\Om_0\dz\theta^0-4\Om_+\dz\theta^1+\tfrac{f_1}{(\al^5_{~5})^2}\theta^0\dz\theta^1+\tfrac32\theta^0\dz\theta^3\nonumber\\
\der\theta^1&=&4\Om_0\dz\theta^1+\tfrac{2f_1}{7(\al^5_{~5})^2}\Om_+\dz\theta^0-3\Om_+\dz\theta^2+\nonumber\\
&&\tfrac{3f_1}{7(\al^5_{~5})^2}\theta^0\dz\theta^2+\tfrac32\theta^1\dz\theta^3\nonumber\\
\der\theta^2&=&2\Om_0\dz\theta^2+\tfrac{4f_1}{7(\al^5_{~5})^2}\Om_+\dz\theta^1-2\Om_+\dz\theta^3-\nonumber\\
&&\tfrac{2f_1^2}{49(\al^5_{~5})^4}\theta^0\dz\theta^1+\tfrac{4f_1}{21(\al^5_{~5})^2}\theta^0\dz\theta^3+\nonumber\\&&\tfrac{f_1}{7(\al^5_{~5})^2}\theta^1\dz\theta^2+\tfrac16\theta^1\dz\theta^4+\tfrac32\theta^2\dz\theta^3\\
\der\theta^3&=&\tfrac{6f_1}{7(\al^5_{~5})^2}\Om_+\dz\theta^2-\Om_+\dz\theta^4-\tfrac{3f_1^2}{49(\al^5_{~5})^4}\theta^0\dz\theta^2+\nonumber\\&&\tfrac{f_1}{14(\al^5_{~5})^2}\theta^0\dz\theta^4+\tfrac{f_1}{7(\al^5_{~5})^2}\theta^1\dz\theta^3+\tfrac34\theta^2\dz\theta^4\nonumber\\
\der\theta^4&=&-2\Om_0\dz\theta^4+\tfrac{8f_1}{7(\al^5_{~5})^2}\Om_+\dz\theta^3-\tfrac{4f_1^2}{49(\al^5_{~5})^4}\theta^0\dz\theta^3+\nonumber\\&&\tfrac{f_1}{7(\al^5_{~5})^2}\theta^1\dz\theta^4+\tfrac32\theta^3\dz\theta^4\nonumber\\
\der\Om_+&=&2\Om_0\dz\Om_++\tfrac{3f_1^2}{98(\al^5_{~5})^4}\theta^0\dz\theta^1+\tfrac{f_1}{14(\al^5_{~5})^2}\theta^0\dz\theta^3+\tfrac18\theta^1\dz\theta^4\nonumber\\
\der\Om_0&=&\tfrac{f_1^2}{49 (\al^5_{~5})^4}\Om_+\dz\theta^0-\tfrac14\Om_+\dz\theta^4+\tfrac{3f_1^2}{196(\al^5_{~5})^4}\theta^0\dz\theta^2+\nonumber\\&&\tfrac{f_1}{56(\al^5_{~5})^2}\theta^0\dz\theta^4+\tfrac{f_1}{14(\al^5_{~5})^2}\theta^1\dz\theta^3+\tfrac{3}{16}\theta^2\dz\theta^4.\nonumber
\end{eqnarray*} 
If $f_1=0$ then to this system there corresponds a unique equivalence class of ODEs
satisfying W\"unschmann conditions and having 7-dimensional transitive
contact symmetry group. The class is represented by 
\ben
F=\frac{5y_4^2}{3y_3}.
\een

In the  case $f_1\neq 0$ we have next \emph{two} nonequivalent classes of ODEs enumerated by the sign of $f_1$ and possessing 6-dimensional transitive contact symmetry groups. Representatives of these classes are
$$F=\frac{5y_4^2}{3y_3}\pm y_3^{5/3},$$
where $\pm1=\sgn f_1$.

In the similar vein we find that the only ODEs related to the case $f_1=2$ are \eqref{ex54} and \eqref{exfrac}.
\end{proof}
\end{proposition}

\subsection{Simple structures with nonvanishing Maxwell form}
All the previous examples satisfy the contact invariant condition
$F_{444}=0$. In this paragraph we give examples of W\"unschmann ODEs
with $F_{444}\neq 0$. As such they will lead to the $\Gl$ structures
with the Maxwell form having a nonzero $\der A^{(7)}$ part. 
First and the simplest example of such equations is
\be\label{exy4} F=(y_4)^{(5/4)}. \ee  
The $\Gl$ structure associated with this ODE has the following properties
$$\der A^{(3)}=0,\qquad R=0, \qquad K=\frac{2}{105}\Rv.$$
It is then an example of a structure with nonvanishing $\der A$ belonging to the 7-dimensional irreducible representation.

Next example is the ODE given by the formula
\begin{eqnarray}
F&=&\frac{1}{9(y_1^2+y_2)^2} \nonumber \\
&&\Big( 5w\big(y_1^6+3y_1^4y_2+9y_1^2y_2^2-9y_2^3-4y_1^3y_3+12y_1y_2y_3+4y_3^2-3y_4(y_1^2+y_2)\big)+\nonumber\\
&&45y_4(y_1^2+y_2)(2y_1y_2+y_3)-4y_1^9-18y_1^7y_2-54y_1^5y_2^2-90y_1^3y_2^3+270y_1y_2^4+ \label{exw} \\
&&15y_1^6y_3+45y_1^4y_2y_3-405y_1^2y_2^2y_3+45y_2^3y_3+60y_1^3y_3^2-180y_1y_2y_3^2-40y_3^3\Big),\nonumber
\end{eqnarray}
where\footnote{Note that $w= 0$ also gives rise to $F$ satisfying
  conditions (\ref{wilki}). But since such $F$ has only
  quadratic $y_4$-dependence it is equivalent to one of proposition
  \ref{p.notorsion}. Actually the one with $c<0$.}
$$w^2=y_1^6 + 3 y_1^4 y_2 + 9 y_1^2 y_2^2 - 9 y_2^3 - 4 y_1^3 y_3 + 
      12y_1 y_2 y_3 + 4 y_3^2 - 3 y_1^2 y_4 - 3 y_2 y_4.$$
Torsion and curvature for the corresponding $\Gl$ structure are complicated and are of general algebraic form.
Both these examples have 6-dimensional transitive group of contact symmetries.

\subsection{A remarkable nonhomogeneous example}
Finally, we present an example of 5th order ODEs satisfying W\"unschmann conditions (\ref{wilki}), which
are generic, in a sense that the function $F$ representing it satisfies 
$F_{444}\neq 0$, but which have the corresponding group of transitive
symmetries of dimension $D<6$. We consider an ansatz in which function
$F$ depends in a special way on only two coordinates $y_3$ and $y_4$. Explicitly:
\be
F=~(y_3)^{5/3}~q\Big(\frac{y_4^3}{y_3^4}\Big),\label{gor}
\ee
where $q=q(z)$ is a sufficiently differentiable real function of its
argument $$z=\frac{y_4^3}{y_3^4}.$$ 
It is remarkable that the above $F$ satisfies \emph{all} W\"unschmann conditions provided that 
\begin{itemize}
\item either $q(z)=\frac53 z^{2/3}$
\item or function $q(z)$ satisfies the following second order ODE:
\be\label{ree}
90z^{4/3}(3q-4z^{2/3})q''-54z^{4/3}{q'}^2+30z^{1/3}(6q-5z^{2/3})q'-25q=0.
\ee
\end{itemize}
In the first case $F=\frac53\frac{y_4^2}{y_3}$, and we recover 
function (\ref{ex53}) with 7-dimensional group of symmetries. Note that
one of the solutions of equation (\ref{ree}) is $q=\frac54
z^{2/3}$, which corresponds to $F=\frac54\frac{y_4^2}{y_3}$. Thus also
the other solution with seven symmetries, the solution
(\ref{ex54}), is covered by this ansatz.  

We observe that if function $q(z)$ satisfies 
\be
25q-60z q'+27 z^{4/3}{q'}^2=0,\label{e2s}
\ee
then it \emph{also} satisfies the reduction (\ref{ree}) of
conditions (\ref{wilki}). Equation (\ref{e2s}) can be solved by 
first putting it in the form
$$q'=\frac{5(2z^{1/3}\pm\sqrt{(4z^{2/3}-3q)})}{9z^{2/3}}$$
and then by integrating, according to the sign $\pm1$. 
In the upper sign case the integration gives $q$ in an implicit form:  
$$
\frac{(2z^{1/3}+\sqrt{(4z^{2/3}-3q)})^{24}(2\sqrt{(4z^{2/3}-3q)}- z^{1/3})^3}{(2\sqrt{(4z^{2/3}-3q)}+  z^{1/3})^3(5z^{2/3}-4q)^3}={\rm const}.
$$
In the lower sign case the implicit equation for $q$ is:
$$
\frac{(2z^{1/3}+\sqrt{(4z^{2/3}-3q)})^{24}(2\sqrt{(4z^{2/3}-3q)}- z^{1/3})^3(5z^{2/3}-4q)^3}{(2\sqrt{(4z^{2/3}-3q)}+  z^{1/3})^3q^{24}}={\rm const}.
$$
Inserting these $q$s into (\ref{gor}) we have a quite nontrivial
W\"unschmann ODE $F=F_\pm$. We close this section with a remark that 
other solutions to the
second order ODE (\ref{ree}) also provide examples of 
5th order W\"unschmann ODEs.

\section{Higher order ODEs}\label{sue}
All our considerations about $\glg(2,\bbR)$ structures associated with 
ODEs of 5th order can be repeated for other orders. This is due to the
following well known fact generalizing proposition \ref{odefl}:
\begin{proposition}\label{odef2}
For every $n\geq 4$, the ordinary differential equation
$$y^{(n)}=0$$
has $\glg(2,\bbR)\times_{\rho_n}\bbR^n$ as its group of contact
symmetries. Here $\rho_n :\glg(2,\bbR)\to\glg(n,\bbR)$ is the
$n$-dimensional irreducible representation of $\glg(2,\bbR)$. 
\end{proposition}

The representation $\rho_n$, at the level of Lie algebra
$\gla(2,\bbR)$, is given in terms of the Lie algebra generators
$$E_+=\bma
0&n-1&0&...&0&0&0\\0&0&n-2&...&0&0&0\\&&&...&&&\\0&0&0&...&3&0&0\\0&0&0&...&0&2&0\\0&0&0&...&0&0&1\\0&0&0&...&0&0&0\ema,\quad
E_-=\bma
0&0&0&...&0&0&0\\1&0&0&...&0&0&0\\0&2&0&...&0&0&0\\0&0&3&...&0&0&0\\&&&...&&&\\
0&0&0&...&n-2&0&0\\0&0&0&...&0&n-1&0\ema,$$ 
$$E_0=\bma
1-n&0&0&...&0&0&0\\0&3-n&0&...&0&0&0\\0&0&5-n&...&0&0&0\\&&&...&&&\\
0&0&0&...&n-5&0&0\\&&&...&0&n-3&0\\0&0&0&...&0&0&n-1\ema,\quad E_1=(1-n)
{\mathbf{1}},$$
where  ${\mathbf{1}}$ is the
$n\times n$ identity matrix. In case of dimension $n=5$ these matrices
coincide with (\ref{mo21}). They also satisfy the same commutation
relations
$$[E_0,E_+]=-2E_+\quad,\quad [E_0,E_-]=2E_-\quad,\quad
[E_+,E_-]=-E_0,$$
where the commutator in the $\gla(2,\bbR)={\rm
  Span}_\bbR(E_-,E_+,E_0,E_1)\subset End(\bbR^n)$ is 
the usual commutator of matrices. 
   
Now, we consider a general $n$-th order
ODE
\be
y^{(n)}=F(x,y,y',y'',y^{(3)},...,y^{(n-1)}),
\label{oden}
\ee
and as before, to simplify the notation, we introduce coordinates 
$x,y,y_1=y',y_2=y'',y_3=y^{(3)},...,y_{n-1}=y^{(n-1)}$ on the 
$(n+1)$-dimensional jet space $J$. Introducing the $n$ \emph{contact forms} 
\begin{eqnarray}
 \om^0&=&\der y-y_1\der x,\nonumber\\
 \om^1&=&\der y_1-y_2\der x,\nonumber\\
 &\vdots&\nonumber\\
 \om^i&=&\der y_i-y_{i+1}\der x,\label{faj}\\
 &\vdots&\nonumber\\
 \om^{n-2}&=&\der y_{n-2}-y_{n-1}\der x,\nonumber\\
 \om^{n-1}&=&\der y_{n-1}-F(x,y,y_1,y_2,...,y_{n-1})\der x\nonumber
\end{eqnarray} 
and the additional 1-form 
$$w_+=\der x,$$
we define a \emph{contact transformation} to be a diffeomorphism $\phi:J\to
J$ which transforms the above $n+1$ one-forms via:
\begin{eqnarray*}
&&\phi^*\om^i=\sum_{k=0}^i\al^i_{~k}\om^k, \quad\quad i=0,1,...n-1,\\
&&\phi^*w_+=\al^n_{~0}\om^0+\al^n_{~1}\om^1+\al^n_{~n}w_+.
\end{eqnarray*}
Here $\al^i_{~j}$ are functions on $J$ such that
$\displaystyle \prod_{i=0}^n\al^i_{~i}\neq 0$ at
each point of $J$. 

Therefore, as in the case of $n=5$, the contact equivalence problem 
for the $n$th order ODEs (\ref{oden}) 
can be studied in terms of the invariant forms 
$(\theta^0,\theta^1,....,\theta^{n-1},\Om_+)$ defined by
\begin{eqnarray}
&&\theta^i=\sum_{k=0}^i\al^i_{~k}\om^k, \quad\quad i=0,1,...n-1,\label{fa}\\
&&\Om_+=\al^n_{~0}\om^0+\al^n_{~1}\om^1+\al^n_{~n}w_+.\nonumber
\end{eqnarray}
These forms initially live on an $\frac{n^2+3n+8}{2}$-dimensional 
manifold $G\to J\times G\to J$, with 
the $G$-factor parametrized by
$\al^i_{~j}$, such that $\displaystyle\prod_{i=0}^n\al^i_{~i}\neq 0$.

Introducing $\gla(2,\bbR)$-valued forms 
\be
\Gamma=\Om_- E_-+\Om_+E_++\Om_0 E_0+\Om_1 E_1,
\label{gma}
\ee
where $(\Om_+,\Om_-,\Om_0,\Om_1)$ are 1-forms on $J\times G$, we can
specialize to $F\equiv 0$, and reformulate proposition \ref{odef2} 
to 
\begin{proposition}\label{pe}
If $F\equiv 0$ then 
one can chose $\frac{n(n+1)}{2}$ parameters $\al^i_{~j}$, as functions
of $x$, $y$, $y_1$, ..., $y_{n-1}$ and the remaining \emph{three}
$\al$s, say $\al^{i_1}_{~j_1}$, $\al^{i_2}_{~j_2}$,
$\al^{i_3}_{~j_3}$, so that the $(n+4)$-dimensional manifold $P$ 
parametrized by  $(x,y,y_1, ..., y_{n-1},\al^{i_1}_{~j_1},
\al^{i_2}_{~j_2}, \al^{i_3}_{~j_3})$ is locally the 
contact symmetry group, $P\cong\glg(2,\bbR)\times_{\rho_n}\bbR^n$, of 
equation $y^{(n)}=0$. Forms (\ref{fa}), after restriction to $P$, can be
supplemented by three additional 1-forms $(\Om_-,\Om_0,\Om_1)$,  
so that  
$(\theta^0,\theta^1,...,\theta^{n-1},\Om_+,\Om_-,\Om_0,\Om_1)$
constitute a basis of the left invariant forms on the Lie 
group $P$. The choice of $\al$s and $\Omega$s is determined by the
requirement that basis $(\theta^0,\theta^1,...,\theta^{n-1},\Om_+,\Om_-,\Om_0,\Om_1)$ satisfies 
\begin{eqnarray}
&&\der\theta+\Gamma\dz\theta=0,\label{zeron}\\
&&\der\Gamma+\Gamma\dz\Gamma=0,\nonumber
\end{eqnarray}
where $\theta=(\theta^0,\theta^1,...,\theta^{n-1})^T$ is a column 
$n$-vector, and $\Gamma$ is given by (\ref{gma}).
\end{proposition}
The defining equations (\ref{zeron}) of the left invariant basis, when
written explicitly in terms of $\theta^i$s and $\Om$s, read
\begin{eqnarray}
\der\theta^0&=&(n-1)(\Om_1+\Om_0)\dz\theta^0+(1-n)\Om_+\dz\theta^1,\nonumber\\
\der\theta^1&=&-\Om_-\dz\theta^0+[(n-1)\Om_1+(n-3)\Om_0]\dz\theta^1+(2-n)\Om_+\dz\theta^2,\nonumber\\
 &\vdots & \nonumber\\
\der\theta^k&=&-k\Om_-\dz\theta^{k-1}+[(n-1)\Om_1+(n-2k-1)\Om_0]\dz\theta^k+\nonumber\\
&&+(1+k-n)\Om_+\dz\theta^{k+1},\label{ixy}\\
 &\vdots & \nonumber\\
\der\theta^{n-1}&=&(1-n)\Om_-\dz\theta^{n-2}+(n-1)(\Om_1-\Om_0)\dz\theta^{n-1},\nonumber\\
\der\Om_+&=&2\Om_0\dz\Om_+,\nonumber\\
\der\Om_-&=&-2\Om_0\dz\Om_-,\nonumber\\
\der\Om_0&=&\Om_+\dz\Om_-,\nonumber\\
\der\Om_1&=&0.\nonumber
\end{eqnarray}
This system can be analyzed in the same spirit as system (\ref{sysc})
of section \ref{nabu}. Thus, we first consider the distribution 
$$\mathfrak{h}=\{X\in{\rm T}P~{\rm
  s.t.}~X\hook\theta^i=0,~i=0,1,2,..,n-1\}$$
annihilating $\theta$.

Then the first $n$ equations of the system (\ref{ixy}) guarantee
  that forms $(\theta^0,\theta^1,\theta^2,,...,\theta^{n-1})$
  satisfy the Fr\"obenius condition,
$$\der\theta^i\dz\theta^0\dz\theta^1\dz\theta^2\dz...\dz\theta^{n-1}=0,\quad\quad\forall~i=0,1,2,...n-1$$
and that, in turn, the distribution $\mathfrak h$ is integrable. Thus
manifold $P$ is foliated by 4-dimensional leaves tangent to the
distribution $\mathfrak h$. The space of leaves of this distribution
$P/{\mathfrak h}$ can be identified with the solution space
$M^n=P/{\mathfrak h}$ of equation $y^{(n)}=0$. This in particular means, that
all equations (\ref{zeron}) can be interpreted respectively as the first 
and the second structure equations for a 
$\gla(2,\bbR)$-valued connection $\Gamma$ having vanishing torsion and
and vanishing curvature. This $\gla(2,\bbR)$-valued connection
originates from a certain $\glg(2,\bbR)$
(conformal) structure on the solution space $M^n$.

To make this last statement more precise we have to invoke a few results
from Hilbert's theory of algebraic invariants \cite{hilb} adapted to
our situation of ODEs. 

\subsection{Results from Hilbert's theory of algebraic invariants}
First we ask if for a given order $n\geq 4$ of an ODE (\ref{oden})
with $F=0$ there exists a bilinear form $\tilde{g}$ on $P$ of 
proposition \ref{pe} such that it projects to a nondegenerate
conformal metric on $M^n$. This is answered, in a bit more general
form, by applying the
\emph{reciprocity law of Hermite} (see \cite{hilb}, p. 60), and its corollaries, due
to Hilbert (see \cite{hilb}, p. 60). 

To adapt Hilbert's results
to our paper we introduce a definition of an \emph{invariant of
  degree} $q$. 
Let $\tilde{t}$ be a \emph{totally symmetric} covariant tensor field
of rank $q$
defined on the group manifold $P$ of proposition \ref{pe}.
\begin{definition}
The tensor field $\tilde{t}$ is called 
a $\glg(2,\bbR)$-\emph{invariant} of degree $q$, if and only if, it 
is degenerate on $\mathfrak h$ and if for every $X\in\mathfrak h$, there
exists a function $c(X)$ on $P$ such that $${\mathcal
  L}_X\tilde{t}=c(X)\tilde{t}.$$ The degeneracy condition means that 
$\tilde{t}(X,...)=0$, for all $X\in\mathfrak h$.
\end{definition}

In the following we will usually abbreviate the term `a $\glg(2,\bbR)$-invariant' to:
`an invariant'. 

The first result from Hilbert's theory, adapted to our situation, is
given by the following  
\begin{proposition}\label{baba}
For every $n=2m+1$, $m=2,3,...$ there exists a \emph{unique}, up to a scale,
invariant $\tilde{g}$ of second degree on $P$. This invariant, a
degenerate symmetric conformal  bilinear form $\tilde{g}$ of signature
$(m+1,m,0,0,0,0)$ on $P$, satisfies 
$${\mathcal L}_X\tilde{g}=2(n-1)(X\hook\Om_1)\tilde{g},$$ for
all $X\in{\mathfrak h}$. 
\end{proposition}
In case of \emph{even} orders $n=2m$, Hilbert's theory gives the following
\begin{proposition} \label{nnn}
For $n=2m$ every $\glg(2,\bbR)$-invariant has degree $q\geq 4$.
\end{proposition}
Thus, if $n=2m$, we do \emph{not} have a conformal metric on the solution
space $M^n$.

Returning to \emph{odd} orders, we present 
the quadratic invariants $\tilde{g}$, of proposition \ref{baba}, for
$n<10$:
\begin{eqnarray}
&{^5\tilde{g}}=3(\theta^2)^2-4\theta^1\theta^3+\theta^0\theta^4,\quad\quad&{\rm
    if}\quad\quad n=5,\nonumber\\
&{^7\tilde{g}}=-10(\theta^3)^2+15\theta^2 \theta^4-6 \theta^1 \theta^5+\theta^0 \theta^6\quad\quad&{\rm
    if}\quad\quad n=7,\label{corr3}\\
&{^9\tilde{g}}=35 (\theta^4)^2-56\theta^3 \theta^5+28 \theta^2\theta^6-8\theta^1\theta^7+\theta^0 \theta^8\quad\quad&{\rm
    if}\quad\quad n=9.\nonumber
\end{eqnarray}
These expressions can be generalized to higher (odd) $n$s. We have the 
following
\begin{proposition}
If $n=2m+1$ and $m\geq 2$, the invariant $\tilde{g}$ of proposition
\ref{baba} is given by: 
$$\tilde{g}=\sum_{j=0}^{m-1} (-1)^j \binom{2m}{j} \theta^j
\theta^{2m-j}+\tfrac12 (-1)^m\binom{2m}{m}(\theta^m)^2.$$ 
\end{proposition}
\begin{remark}\label{ics}
This proposition is also valid for $m=1$. For such $m$, the value of
$n$ is $n=3$, and we are in the regime of \emph{third} order
ODEs. Such ODEs were considered by W\"unschmann \cite{wun}. Since $^3\tilde{g}=\theta^0\theta^2-(\theta^1)^2$ is the only invariant in this case, the counterpart of the bundle $P$ of proposition \ref{pe} is a 10-dimensional bundle $P\cong \og(2,3)$, the full conformal group in Lorentzian signature $(1,2)$. The counterpart of system (\ref{zeron})/(\ref{ixy}) is given by Maurer-Cartan equations for $\og(3,2)$:  
\begin{eqnarray*}
\der \theta^0&=&2(\Om_1+\Om_0)\dz\theta^0-2\Om_+\dz\theta^1,\\
\der \theta^1&=&-\Om_-\dz\theta^0+2\Om_1\dz\theta^1-\Om_+\dz\theta^2,\\
\der \theta^2&=&-2\Om_1\theta^1+(2\Om_1-2\Om_0)\dz\theta^2, \\
\der \Om_+&=&2\Om_0\dz\Om_++\tfrac{1}{2}\Om_3\dz\theta^0+\Om_4\dz\theta^1, \\
\der \Om_-&=&-2\Om_0\dz\Om_-+\Om_2\dz\theta^1+\tfrac{1}{2}\Om_3\dz\theta^2, \\
\der \Om_0&=&\Om_+\dz\Om_--\tfrac{1}{2}\Om_2\dz\theta^0+\tfrac{1}{2}\Om_4\dz\theta^2, \\
\der \Om_1&=&-\tfrac{1}{2}\Om_2\dz\theta^0-\tfrac{1}{2}\Om_3\dz\theta^1-\tfrac{1}{2}\Om_4\dz\theta^2, \\
\der \Om_2&=&-\Om_3\dz\Om_-+2\Om_2\dz\Om_0+2\Om_2\dz\Om_1, \\
\der \Om_3&=&-2\Om_2\dz\Om_+-2\Om_4\dz\Om_-+2\Om_3\dz\Om_1, \\
\der \Om_4&=&-2\Om_4\dz\Om_0-\Om_3\dz\Om_++2\Om_4\dz\Om_1.
\end{eqnarray*}
Here, apart from $\theta^0,\theta^1,\theta^2$ and $\Om_+,\Om_-,\Om_0,\Om_1$ we have also left invariant forms $\Om_2,\Om_3,\Om_4$.
\end{remark}
Now we pass to the invariants of degree $q=3$. The question of
their existence was again 
determined by Hilbert (see \cite{hilb}, p. 60), in terms of the 
\emph{reciprocity law of Hermite}. In the language of our paper we
have the following
\begin{proposition}
An invariant of third degree $\tilde{\ten}$ 
exists on $P$ if and only if
$$n=4\mu+1,\quad\quad \mu\in\bbN.$$
\end{proposition}
Hilbert's theory,
\cite{hilb}, p. 60, implies also the following:
\begin{proposition}\label{cor8}
In low dimensions $n=4\mu+1$, the unique up to a scale
\emph{cubic} invariant is given by
\begin{itemize}
\item $n=5$:
$${^5\tilde{\ten}}=(\theta^2)^3-2\theta^1\theta^2\theta^3+\theta^0(\theta^3)^2-\theta^0\theta^2\theta^4+(\theta^1)^2\theta^4$$
\item $n=9$:
\begin{eqnarray*}
&&{^9\tilde{\ten}}=15 (\theta^4)^3-36 \theta^3 \theta^4 \theta^5+24 \theta^2 (\theta^5)^2+24
  (\theta^3)^2\theta^6-22 \theta^2 \theta^4 \theta^6-\\
&&8 \theta^1 \theta^5 \theta^6+3 \theta^0(\theta^6)^2-8\theta^2 \theta^3 \theta^7+12 \theta^1
  \theta^4 \theta^7-4 \theta^0 \theta^5 \theta^7+\\
&&3 (\theta^2)^2 \theta^8-4 \theta^1 \theta^3 \theta^8+\theta^0 \theta^4 \theta^8.
\end{eqnarray*}
\end{itemize} 
\end{proposition}
The rough statement about the even orders, $n=2m$, described in
proposition \ref{nnn}, can be again refined in terms of the 
reciprocity law of Hermite. Following Hilbert we have 
\begin{proposition}
If $4\leq n=2m$ the lowest order invariant 
tensor $\tilde{\ten}$ on $P$ has degree \emph{four}. This 
is unique (up to a scale) only if $n=4,6,8,12$. 
If $n=10$ or $n=14$ we have \emph{two}
  independent quartic invariants $\tilde{\ten}$; if
  $n=16,18,20$ we have 
\emph{three} independent quartic invariants; and so on.
\end{proposition}

\begin{proposition}\label{cor9}
In low dimensions $n=2m$, the \emph{quartic} invariant 
tensor $\tilde{\ten}$ on $P$ is given by
\begin{itemize}
\item $n=4$:
$${^4\tilde{\ten}}=-3(\theta^1)^2 (\theta^2)^2 + 4\theta^0 (\theta^2)^3 + 4 (\theta^1)^3\theta^3 -6 \theta^0 \theta^1\theta^2 \theta^3 + (\theta^0)^2 (\theta^3)^2$$
\item $n=6$:
\begin{eqnarray*}
&&{^6\tilde{\ten}}=-32(\theta^2)^2(\theta^3)^2 + 48\theta^1 (\theta^3)^3 + 48 (\theta^2)^3\theta^4 - 76\theta^1\theta^2\theta^3\theta^4 -\\&& 12 \theta^0 (\theta^3)^2 \theta^4 + 9(\theta^1)^2 (\theta^4)^2 + 16 \theta^0 \theta^2 (\theta^4)^2 - 12 \theta^1 (\theta^2)^2 \theta^5+\\&& 16(\theta^1)^2 \theta^3 \theta^5 + 4\theta^0\theta^2 \theta^3 \theta^5 - 10 \theta^0 \theta^1 \theta^4 \theta^5 + (\theta^0)^2 (\theta^5)^2.
\end{eqnarray*}
\item $n=8$:
\begin{eqnarray*}
&&{^8\tilde{\ten}}=-375(\theta^3)^2 (\theta^4)^2 + 600 \theta^2 (\theta^4)^3 + 600 (\theta^3)^3\theta^5 - 
    990 \theta^2 \theta^3 \theta^4 \theta^5 -\\&& 240\theta^1 (\theta^4)^2 \theta^5 + 
    81 (\theta^2)^2(\theta^5)^2 + 360 \theta^1\theta^3 (\theta^5)^2 - 
    240 \theta^2 (\theta^3)^2 \theta^6 + 360 (\theta^2)^2 \theta^4 \theta^6 +\\&& 
    50 \theta^1 \theta^3 \theta^4 \theta^6 + 40 \theta^0 (\theta^4)^2 \theta^6 - 
    234 \theta^1 \theta^2 \theta^5 \theta^6 - 60 \theta^0 \theta^3 \theta^5 \theta^6 + 
    25 (\theta^1)^2 (\theta^6)^2 + \\&&24\theta^0 \theta^2 (\theta^6)^2 + 
    40 \theta^1 (\theta^3)^2\theta^7 - 60 \theta^1\theta^2 \theta^4 \theta^7 - 
    10 \theta^0 \theta^3 \theta^4 \theta^7 + 24 (\theta^1)^2 \theta^5 \theta^7 +\\&& 
    18 \theta^0 \theta^2 \theta^5 \theta^7 - 14 \theta^0 \theta^1 \theta^6 \theta^7 + 
    (\theta^0)^2 (\theta^7)^2.
\end{eqnarray*}
\end{itemize} 
\end{proposition}
Among the small dimensions $n=7$ is quite special, since here the next 
invariant linearly and \emph{functionally} independent of 
the metric $\tilde{g}$ has $q=4$. We have the following 
\begin{proposition}\label{cor7}
In dimension $n=7$, the invariant of the lowest degree is the metric
${^7\tilde{g}}$. There are \emph{no} invariants of degree $q=3$ and \emph{only
  two} linearly independent, invariants of
degree $q=4$. The first of them is ${^7\tilde{g}}^2$. The second
can be chosen to be 
\begin{eqnarray*}
&&{^7\tilde{\ten}}=160(\theta^3)^4 - 480\theta^2 (\theta^3)^2 \theta^4 + 1035 (\theta^2)^2 (\theta^4)^2 - 
    1080 \theta^1 \theta^3 (\theta^4)^2 + 540 \theta^0 (\theta^4)^3 -\\&& 
    1080 (\theta^2)^2\theta^3 \theta^5 + 1920 \theta^1 (\theta^3)^2 \theta^5 - 
    180 \theta^1 \theta^2 \theta^4 \theta^5 - 1080 \theta^0 \theta^3 \theta^4 \theta^5 - 
    288 (\theta^1)^2 (\theta^5)^2 +\\&& 540 \theta^0 \theta^2 (\theta^5)^2 + 540 (\theta^2)^3 \theta^6 - 
    1080 \theta^1 \theta^2\theta^3 \theta^6 + 400 \theta^0 (\theta^3)^2 \theta^6 + 
    540 (\theta^1)^2 \theta^4 \theta^6 -\\&& 330 \theta^0 \theta^2 \theta^4 \theta^6 - 
    84 \theta^0\theta^1 \theta^5 \theta^6 + 7 (\theta^0)^2 (\theta^6)^2.
\end{eqnarray*}
\end{proposition}
\subsection{Stabilizers of the irreducible $\glg(2,\bbR)$ in
    dimensions $n<10$}\label{8812}
In dimensions $n\leq 10$ the $\glg(2,\bbR)$ 
invariant tensors of low order $q\leq 4$ turn out to be sufficient to
reduce the $\glg(n,\bbR)$ group to $\glg(2,\bbR)$ in its irreducible
$n$-dimensional representation. 

Given an invariant tensor
$$\tilde{t}=\frac{1}{q!}t_{i_1i_2...i_q}\theta^{i_1}...\theta^{i_q}$$ 
of degree $q$ on $P$ and a $\glg(n,\bbR)$-valued function
$a=(a^i_{~j})$ on $P$, at every point $p\in P$, we have a $\glg(n,\bbR)$-action 
$$(a^i_{~j},\tilde{t}_{i_1i_2...i_q})\mapsto(\rho_n(a)\tilde{t})_{j_1j_2...j_q}=
a^{i_1}_{~j_1}a^{i_2}_{~j_2}...a^{i_q}_{~j_q}\tilde{t}_{i_1i_2...i_q}.$$
A subgroup $G_{\tilde{t}}$ of $\glg(n,\bbR)$ consisting of
$a=(a^i_{~j})$ 
such that $$\rho_n(a)\tilde{t}=(\det{a})^{q/n}\,\tilde{t},$$
is the stabilizer of $\tilde{t}$ at $p\in P$. Since $\tilde{t}$ is an
invariant then, obviously $\glg(2,\bbR)\subset G_{\tilde{t}}$.

This
leads to the following question: how many invariants is needed 
in dimension $n$ so that its common stabilizer is \emph{precisely}
$\glg(2,\bbR)$ in its $n$ dimensional irreducible
representation?

Inspecting Hilbert's results we checked that in dimensions $4\leq
n\leq 9$ we have 
\begin{theorem}\label{stt}
For each $n=4,5,6,7,8,9$, the \emph{full} stabilizer group of the
respective invariant  
tensor ${^n\tilde{\ten}}$ of propositions \ref{cor8}, \ref{cor9}, \ref{cor7}, is the 
group $\glg(2,\bbR)$ in the $n$-dimensional irreducible representation 
$\rho_n$. In
particular, if $n=5,7,9$ these stabilizers are subgroups of the respective
pseudohomothetic groups $\cog(3,2)$, $\cog(4,3)$ and 
$\cog(5,4)$, each in its defining representation.  
\end{theorem}

Thus in each of these dimensions it is the lowest order \emph{nonquadratic}
invariant what is responsible for the full reduction from $\glg(n,\bbR)$ to $\glg(2,\bbR)$.

\begin{remark}
In dimension $n=5$, using (\ref{corr3}) and proposition \ref{cor8} we
  define a \emph{conformal metric} $[^5g_{ij}]$
  represented by 
$${^5g}_{ij}=\frac12\frac{\partial^2}{\partial \theta^i\partial
    \theta^j}\big({^5\tilde{g}}\big),\quad\quad i,j=0,1,2,3,4$$
  and a \emph{conformal symmetric tensor of
  third degree} $[^5\ten_{ijk}]$
  represented by 
$${^5\ten}_{ijk}=-\frac{\sqrt{3}}{8}\frac{\partial^3}{\partial \theta^i\partial
    \theta^j\partial \theta^k}\big({^5\tilde{\ten}}\big),\quad\quad
    i,j,k,l=0,1,2,3,4.$$

The convenient factor $-\frac{\sqrt{3}}{8}$ in the expression for
    ${^5\ten}_{ijk}$ was chosen so that the pair
    $({^5g}_{ij},{^5\ten}_{ijk})$ satisfies Cartan's
    identities (i)-(iii) of section \ref{s1}. This leads to the $\glg(2,\bbR)$
    geometries in dimension 5 considered in sections \ref{trzy}-\ref{secsec}. 
\end{remark}
\begin{remark}
In the next odd dimension situation is quite similar, but now we have
a quartic invariant ${^7\tilde{\ten}}$. Thus apart from the 
\emph{conformal metric} $[^7g_{ij}]$
  represented by 
$${^7g}_{ij}=\frac12\frac{\partial^2}{\partial \theta^i\partial
    \theta^j}\big({^7\tilde{g}}\big),\quad\quad i,j=0,1,2,3,4,5,6$$
we have a \emph{conformal symmetric tensor of
  fourth degree} $[^7\ten_{ijkl}]$
  represented by 
\be
{^7\ten}_{ijkl}=\frac{1}{24}\frac{\partial^4}{\partial \theta^i\partial
    \theta^j\partial \theta^k\partial \theta^l}\big({^7\tilde{\ten}}\big),\quad\quad
    i,j,k,l=0,1,2,3,4,5,6.\label{7yp}\ee
Note that ${^7\tilde{\ten}}$ of
proposition \ref{cor7} was chosen in such a way that the fourth order
${^7\ten}_{ijkl}$ satisfied 
$$^7g^{ij}~^7\ten_{ijkl}=0,\quad\quad{\rm where}\quad\quad
 ^7g^{ij}~^7g_{jk}={\delta}^i_{~k}.$$
This choice of the fourth order invariant is nevertheless arbitrary,
 since we can always get another invariant of the fourth order by
 replacing ${^7\ten}$ with 
$${^7\bar{\ten}_{ijkl}}=c_1
 {^7\tilde{\ten}_{ijkl}}+c_2\, {^7\tilde{g}_{(ij}}\,{^7\tilde{g}_{kl)}}.$$
It is interesting to note that the choice
$$c_1=\frac{2\sqrt{5}}{\sqrt{3147}},\quad\quad c_2=\frac{34}{\sqrt{15735}}$$
applied to ${^7\bar{\ten}}$, leads, via formula like (\ref{7yp}), to 
${^7{\bar{\ten}}}_{ijkl}$ satisfying Cartan-like identity:
$${^7g^{ih}}\,{^7g^{ef}}\,{^7\bar{\ten}}_{ie(jk}{^7\bar{\ten}}_{lm)fh}={^7g}_{(jk}{^7g}_{lm)}$$
 and $${^7g^{ij}}\,{^7{\bar{\ten}}}_{ijkl}=\tfrac{3}{2}c_2\,{^7g_{kl}},\quad\quad{\rm where}\quad\quad
 {^7g^{ij}}\,{^7g}_{jk}={\delta}^i_{~k}.$$
Note also that the above Cartan--like identities are
 preserved under the
 conformal transformation 
$$({^7g}_{ij},{^7\bar{\ten}}_{ijkl})\mapsto
 ({^7g'}_{ij},{^5\bar{\ten}'}_{ijkl})=({\rm e}^{2\phi}~{^7g}_{ij},{\rm
 e}^{4\phi}~{^7\bar{\ten}}_{ijkl}),$$
where $\phi\in\bbR$.

Thus the $\glg(2,\bbR)$ geometries in dimension $n=7$ may be defined
by a conformal class of pairs of tensors $[{^7g}_{ij},{^7\bar{\ten}}_{ijkl}]$
 with the properties and transformations as above. 
\end{remark}

\begin{remark}
By analogy, in dimensions $n=4,6,8$, the irreducible $\glg(2,\bbR)$
geometries may be described in terms of a conformal tensor 
$[^n\ten_{ijkl}]$
  represented by 
$${^n\ten}_{ijkl}=\frac{1}{24}\frac{\partial^4}{\partial \theta^i\partial
    \theta^j\partial \theta^k\partial
    \theta^l}\big({^n\tilde{\ten}}\big),\quad\quad
    i,j,k,l=0,1,2,...,n-1,$$
and obtained in terms of the respective \emph{quartic} 
invariants ${^n\tilde{\ten}}$ of proposition \ref{cor9}.
\end{remark}

\begin{remark}
Dimension $n=9$ is similar to dimension $n=5$. A periodicity with
period \emph{four} is a remarkable feature of Hilbert's theory of
algebraic invariants \cite{hilb}, p. 60.
\end{remark}
\subsection{W\"unschmann conditions for the existence of $\glg(2,\bbR)$ 
geometries on the solution space of ODEs}\label{8813}
An invariant tensor $\tilde{t}$, by its very definition, has a property
that it descends to a nondegenerate conformal tensor $[t]$ on the
solutions space $M^n=P/\mathfrak h$ of the equation $y^{(n)}=0$. In particular in
dimensions $4\leq n\leq 9$ the conformal class $[^n\ten]$,
corresponding to invariant tensors ${^n\tilde{\ten}}$ reduces the
structure group of $M^n$ to $\glg(2,\bbR)$ defining an irreducible 
$\glg(2,\bbR)$ geometry 
there. We do not know how many invariant tensors are needed 
to achieve this reduction for $n>9$, but it is obvious that for a
given $n$ this number is finite, say $w_n$. Thus for each $n\geq 3$ we
have a finite number of invariants ${^n\tilde{\ten}_I}$,
$I=1,2,...w_n$, which descend to the solution space $M^n$ of the
equation $y^{(n)}=0$ equipping it with a $\glg(2,\bbR)$
structure. It is important that each of the invariants
${^n\tilde{\ten}_I}$ has only  \emph{constant} coefficients when
expressed in terms of the invariant coframe
$(\theta^0,...,\theta^{n-1})$ on $P$ 
(see, for example, every ${^n\tilde{\ten}}$ of the preceding section). 

Now, we return to a \emph{general} $n$-th order ODE (\ref{oden}).
Thus we now have a general function
$F(x,y,y',y'',y^{(3)},...,y^{(n-1)})$, which determines 
the contact forms $(\om^0,\om^1,...,\om^{n-1},w_+)$ by 
(\ref{faj}). Corresponding to these forms we have the 
invariant forms $(\theta^0,...,\theta^{n-1},\Om_+)$ of 
(\ref{fa}), which live on bundle $J\times G$ over $J$. We can now ask the
following question (this generalizes to arbitrary $n>3$ the similar
question of section \ref{secsec}): What shall we assume about $F$ defining the
contact equivalence class of ODEs (\ref{oden}) that there exists a
$(4+n)$-dimensional subbundle $P$ of $J\times G$ on which the forms
$(\theta^0,...,\theta^{n-1},\Om_+)$ satisfy:
\begin{eqnarray}
\der\theta^0&=&(n-1)(\Om_1+\Om_0)\dz\theta^0+(1-n)\Om_+\dz\theta^1+\frac12T^0_{~ij}\theta^i\dz\theta^j,\nonumber\\
\der\theta^1&=&-\Om_-\dz\theta^0+[(n-1)\Om_1+(n-3)\Om_0]\dz\theta^1+\nonumber\\
&&+(2-n)\Om_+\dz\theta^2+\frac12T^1_{~ij}\theta^i\dz\theta^j,\nonumber\\
&\vdots&\nonumber\\
\der\theta^k&=&-k\Om_-\dz\theta^{k-1}+[(n-1)\Om_1+(n-2k-1)\Om_0]\dz\theta^k+\nonumber\\
&&+(1+k-n)\Om_+\dz\theta^{k+1}+\frac12T^k_{~ij}\theta^i\dz\theta^j,\label{nixy}\\
&\vdots&\nonumber\\
\der\theta^{n-1}&=&(1-n)\Om_-\dz\theta^{n-2}+(n-1)(\Om_1-\Om_0)\dz\theta^{n-1}+\frac12T^{n-1}_{~ij}\theta^i\dz\theta^j,\nonumber\\
\der\Om_+&=&2\Om_0\dz\Om_++\frac12R_{+ij}\theta^i\dz\theta^j,\nonumber\\
\der\Om_-&=&-2\Om_0\dz\Om_-+\frac12R_{-ij}\theta^i\dz\theta^j,\nonumber\\
\der\Om_0&=&\Om_+\dz\Om_-+\frac12R_{0ij}\theta^i\dz\theta^j,\nonumber\\
\der\Om_1&=&\frac12R_{ij}\theta^i\dz\theta^j.\nonumber
\end{eqnarray}
As first observed by W\"unschmann \cite{wun} and then successively used by 
Newman and collaborators \cite{newman} this question can be
reformulated into a nicer one. To make this reformulation we repeat our
arguments from section \ref{odef2}.

Suppose that we are able to satisfy system (\ref{nixy}) by forms
(\ref{fa}). Consider the distribution 
$$\mathfrak{h}=\{X\in{\rm T}P~{\rm
  s.t.}~X\hook\theta^i=0,~i=0,1,2,..,n-1\}$$
annihilating $\theta$s. Despite of the fact that system (\ref{nixy})
  involves \emph{new terms}, when compared with system (\ref{ixy}),
  they do not destroy the integrability of the distribution
  $\mathfrak{h}$; the first $n$ equations (\ref{nixy}) still
  guarantee that $\mathfrak{h}$ is \emph{integrable}. Thus
manifold $P$ is foliated by 4-dimensional leaves tangent to the
distribution $\mathfrak h$. The space of leaves of this distribution
$P/{\mathfrak h}$ can be identified with the solution space
$M^n=P/{\mathfrak h}$ of equation (\ref{oden}). Now, on manifold 
$P$ of system (\ref{nixy}), we define $w_n$ tensors
  ${^n\tilde{\ten}_I}$, 
which formally are given by \emph{the same} formulae that defined the
  $w_n$ invariants ${^n\tilde{\ten}_I}$ of the flat system
  (\ref{ixy}) needed to get the full reduction to $\glg(2,\bbR)$. So,
  when defining the present ${^n\tilde{\ten}_I}$, we use the same
  formulae as for the $y^{(n)}=0$ case, replacing forms $\theta$ of
  the flat case, with forms $\theta$ satisfying system
  (\ref{nixy}). It is now easy to verify that the question about the
  conditions on $F$ to admit $P$ with system (\ref{nixy}) is
  equivalent to the requirement that \emph{all} $w_n$ 
tensors ${^n\tilde{\ten}_I}$
  transform \emph{conformally} when Lie transported along the leaves of
  distribution $\mathfrak h$. Infinitesimally this condition is
  equivalent to the existence of 
functions $c_I(X)$ on $P$ such that  
$${\mathcal L}_X({^n\tilde{\ten}_I})=c_I(X)~{^n\tilde{\ten}_I},$$ $\forall
  X\in{\mathfrak h}$, and $\forall I=1,2,...w_n$.
If this is satisfied then tensors ${^n\tilde{\ten}_I}$ descend to
  a conformal class of tensors
 $[{^n\ten_1},{^n\ten_2},..., {^n\ten_{w_n}}]$ on the solution space 
$M^n$ defining a $\glg(2,\bbR)$ there.  

We know that in dimension $n=5$ the conformal preservation of
${^5\tilde{g}}$ and ${^5\tilde{\ten}}$ is equivalent to the requirement
  on function $F=F(x,y,y_1,y_2,y_3,y_4)$ to satisfy W\"unschmann 
  conditions (\ref{wilki}). The generalization of this fact to other
  low dimensions $4\leq n<10$ is given by the following 
\begin{theorem}
Let $M^n$ be the solution space of $n$th order ODE 
\be
y^{(n)}=F(x,y,y',y'',y^{(3)},...,y^{(n-1)}),\label{rano}
\ee
with $4\leq n<10$, and let
$$\D=\partial_x+y_1\partial_y+y_2\partial_{y_1}+\ldots+y_{n-1}\partial_{y_{n-2}}+F\partial_{y_{n-1}}$$
be the total derivative.
The necessary conditions for a contact equivalence
class of ODEs (\ref{rano}) to define a
principal $\glg(2,\bbR)$-bundle $\glg(2,\bbR)\to P\to M^n$ with 
invariants forms 
$(\theta^0,...,\theta^{n-1},\Om_+,\Om_-,\Om_0,\Om_1)$ satisfying system
(\ref{nixy}) is that the defining function $F$ of (\ref{rano}) satisfies $n-2$
W\"unschmann conditions given below:
\begin{itemize}
\item $n=4$:
\begin{eqnarray}
&&4\D ^2F_3-8\D F_2+8F_1-6\D F_3 F_3+4F_2 F_3+F_3^3=0,\nonumber\\
&&\nonumber\\
&&160\D ^2F_2-640\D F_1+144(\D F_3)^2-352 \D F_3 F_2+144 F_2^2 -\nonumber\\
&&80\D F_2 F_3+160 F_1 F_3-72 \D F_3 F_3^2+88 F_2 F_3^2+9
  F_3^4+16000F_y=0,
\nonumber\end{eqnarray}
\item $n=5$:
\begin{eqnarray}
&&50\D ^2F_4 - 75 \D F_3 + 50 F_2 - 60 F_4\D F_4 + 30 F_3F_4 + 8 F_4^3=0\nonumber\\
&&\nonumber\\
&&375 \D ^2F_3 - 1000 \D F_2 + 350 \D F_4^2 + 1250 F_1 - 650F_3 \D F_4  + 
    200 F_3^2 -\nonumber\\&& 150 F_4 \D F_3 + 200 F_2 F_4 - 140 F_4^2 \D F_4 + 
    130 F_3 F_4^2 + 14 F_4^4=0\nonumber\\
&&\nonumber\\
&&1250 \D ^2F_2 - 6250 \D F_1 + 1750 \D F_3\D F_4 - 2750 F_2\D F_4  -\nonumber
\\&& 875  F_3\D F_3
  +1250 F_2 F_3 - 500F_4 \D F_2  + 700 (\D F_4)^2 F_4 +\nonumber\\
&& 1250 F_1 F_4
  - 
    1050 F_3 F_4\D F_4  + 350 F_3^2 F_4 - 350 F_4^2\D F_3  +\nonumber \\
&& 550 F_2 F_4^2
  -
    280 F_4^3 \D F_4 + 210 F_3 F_4^3 + 28 F_4^5 + 18750 F_y=0.\nonumber
\end{eqnarray}
\item $n=6$:
\begin{eqnarray*}
&&45\D ^2F_5 - 54\D F_4 + 27 F_3 - 45 \D F_5 F_5 + 18 F_4F_5 + 5 F_5^3\\
&&\\
&&945\D ^2F_4 - 1890 \D F_3 + 900 (\D F_5)^2 + 1575 F_2 - 1350 \D F_5 F_4 + 
    333 F_4^2 -\\&& 315\D F_4 F_5 + 315 F_3 F_5 - 300 \D F_5 F_5^2 + 
    225 F_4 F_5^2 + 25 F_5^4=0\\
&&\\
&&2835 \D ^2F_3 - 9450\D F_2 + 4320 \D F_4 \D F_5 + 14175 F_1 - 5130 \D F_5 F_3 -\\&& 
    1728 \D F_4 F_4 + 1863 F_3 F_4 - 945 \D F_3 F_5 + 1800 (\D F_5)^2F_5 + 
    1575F_2F_5 -\\&& 2160 \D F_5 F_4 F_5 + 576 F_4^2 F_5 - 720 \D F_4 F_5^2 + 
    855 F_3 F_5^2 -\\&& 600 \D F_5 F_5^3 + 360 F_4 F_5^3 + 50 F_5^5=0\\
\end{eqnarray*}
\begin{eqnarray*}
&&14175 \D ^2F_2 - 85050 \D F_1 + 6480 (\D F_4)^2 + 16200 \D F_3 \D F_5 -\\&& 
    31050 \D F_5 F_2 - 9720 \D F_4 F_3 + 3645 F_3^2 - 6480 \D F_3 F_4 +\\&& 
    5400 \D F_5^2 F_4 + 11475 F_2 F_4 - 4320 \D F_5F_4^2 + 864 F_4^3 - 
    4725 \D F_2 F_5 +\\&& 10800 \D F_4 \D F_5 F_5 + 14175 F_1 F_5 - 
    10800 \D F_5 F_3 F_5 - 6480 \D F_4F_4 F_5 +\\&& 5940 F_3F_4 F_5 - 
    2700 \D F_3 F_5^2 + 4500 (\D F_5)^2 F_5^2 + 5175 F_2 F_5^2 -\\&& 
    7200 \D F_5 F_4 F_5^2 + 2340 F_4^2 F_5^2 - 1800 \D F_4 F_5^3 + 
    1800 F_3 F_5^3 - 1500 \D F_5 F_5^4 + \\&&1050 F_4 F_5^4 + 125 F_5^6 + 
    297675 F_y=0
\end{eqnarray*}
\item $n=7$:
\begin{eqnarray*}
&&245 \D ^2F_6 - 245 \D F_5 + 98 F_4 - 210 \D F_6 F_6 + 70 F_5 F_6 + 20
  F_6^3=0\\
&&\\
&&6860 \D ^2F_5 - 10976 \D F_4 + 6615 (\D F_6)^2 + 6860 F_3 - 8330 \D F_6 F_5 +\\&& 
    1715 F_5^2 - 1960 \D F_5 F_6 + 1568 F_4 F_6 - 1890 \D F_6 F_6^2 + 
    1190 F_5 F_6^2 + 135 F_6^4=0\\&&\\
&&9604 \D ^2F_4 - 24010 \D F_3 + 15435 \D F_5 \D F_6 + 24010 F_2 - 14749 \D F_6
  F_4 -
\\&& 
    5145 \D F_5 F_5 + 4459 F_4 F_5 - 2744 \D F_4 F_6 + 6615 (\D F_6)^2 F_6 + 
    3430 F_3 F_6 -
\\&& 6615 \D F_6 F_5 F_6 + 1470 F_5^2 F_6 - 2205 \D F_5 F_6^2 + 
    2107 F_4 F_6^2 -\\&& 1890 \D F_6 F_6^3 + 945 F_5 F_6^3 + 135
    F_6^5=0\\&&\\
&&336140 \D ^2F_3 - 1344560 \D F_2 + 180075 (\D F_5)^2 + 432180 \D F_4 \D F_6 +\\&& 
    2352980 F_1 - 624260 \D F_6 F_3 - 216090 \D F_5 F_4 + 64827 F_4^2 - \\&&
    144060 \D F_4 F_5 + 154350 (\D F_6)^2 F_5 + 192080 F_3 F_5 - 
    102900 \D F_6 F_5^2 +\\&& 17150 F_5^3 - 96040 \D F_3F_6 + 
    308700 \D F_5 \D F_6 F_6 + 192080 F_2 F_6 -\\&& 246960 \D F_6 F_4 F_6 - 
    154350 \D F_5 F_5 F_6 + 113190 F_4 F_5 F_6 - 61740 \D F_4 F_6^2 +\\&& 
    132300 (\D F_6)^2 F_6^2 + 89180 F_3 F_6^2 - 176400 \D F_6 F_5 F_6^2 + 
    47775 F_5^2 F_6^2 -\\&& 44100 \D F_5 F_6^3 + 35280 F_4F_6^3 - 
    37800 \D F_6F_6^4 + 22050 F_5 F_6^4 + 2700 F_6^6=0\\&&\\
&&2352980 \D ^2F_2 - 16470860 \D F_1 + 1512630 \D F_4 \D F_5 + 
      2268945 \D F_3 \D F_6 - \\&&5126135 \D F_6 F_2 - 1512630 \D F_5 F_3 - 
      907578 \D F_4 F_4 + 648270 (\D F_6)^2 F_4 +\\&& 907578 F_3 F_4 - 
      756315 \D F_3 F_5 + 1080450 \D F_5 \D F_6 F_5 + 1596665 F_2 F_5 - \\&&
      1080450 \D F_6 F_4F_5 - 360150 \D F_5 F_5^2 + 288120 F_4 F_5^2 - 
      672280 \D F_2 F_6 +\\&& 540225 (\D F_5)^2F_6 + 1296540 \D F_4 \D F_6 F_6 + 
      2352980 F_1 F_6 -\\&& 1620675 \D F_6 F_3 F_6 - 864360 \D F_5 F_4 F_6 + 
      324135 F_4^2 F_6 - 648270 \D F_4 F_5 F_6 +\\&& 926100 (\D F_6)^2 F_5 F_6 + 
      756315 F_3 F_5 F_6 - 771750 \D F_6 F_5^2 F_6 + 154350 F_5^3 F_6 -\\&& 
      324135 \D F_3 F_6^2 + 926100 \D F_5 \D F_6 F_6^2 + 732305 F_2 F_6^2 - 
      926100 \D F_6 F_4 F_6^2 -\\&& 617400 \D F_5 F_5 F_6^2 + 
      524790 F_4 F_5 F_6^2 - 185220 \D F_4 F_6^3 + 396900 (\D F_6)^2 F_6^3
  + 
\\&&
      231525 F_3 F_6^3 - 661500 \D F_6 F_5 F_6^3 + 209475 F_5^2 F_6^3 - 
      132300 \D F_5 F_6^4 +\\&& 119070 F_4 F_6^4 - 113400 \D F_6 F_6^5 + 
      75600 F_5 F_6^5 + 8100 F_6^7 + 65883440 F_y=0.
\end{eqnarray*}
\end{itemize} 
\end{theorem}
\begin{remark} Although we calculated the W\"unschmann conditions
  for $n=8$ and $n=9$, we do not present them here due to their
  length. We remark, however that in any order $n\geq 4$, the $n-2$
  W\"unschmann conditions, which by the very definition are conditions needed
  for an ODE to define a $\Gl$ geometry on its solution space, are
  always of the
  \emph{third} order in the derivatives of the function $F$ which
  defines an ODE. In this sense they differ from the generalizations of
  W\"unschmann conditions obtained by \cite{dub} and \cite{dun}. 
\end{remark}

\begin{remark}
If $n=3$ we have only one 
W\"unschmann condition \cite{chern,wun}:
$$9 \D ^2F_2 - 27 \D F_1 - 18 \D F_2 F_2 + 18 F_1 F_2 + 4 F_2^3 + 54 F_y=0.
$$ and, if it satisfied, a conformal Lorentzian geometry associated with
a metric 
$${^3g}=\theta^0\theta^2-(\theta^1)^2$$ is naturally defined 
on the solution space. 
\end{remark}
\begin{remark}
If $n=4$ the ODEs satisfying the two W\"unschmann conditions lead 
to very {\it nontrivial} geometries on 4-dimensional solution
spaces. These are a sort of conformal Weyl geometries, which instead
of a metric are define in terms of the conformal rank four tensor
${^4\ten}$. These geometries define a characteristic connection, which
is $\gla(2,\bbR)$ valued and has an exotic holonomy \cite{bryantsp4}. 
By this we mean
that the holonomy of this \emph{nonmetric} but \emph{torsionless}
connection does not appear on the Berger's list \cite{bryantsp4}. See also our account on this subject in \cite{nur4}.
\end{remark}
\begin{remark}
Our studies of the ODES with $n=3,4,5$, and the preliminary results
about the cases
with $n\geq 7$, make us to conjecture that if $n\geq 7$ then the 
$n-2$ W\"unschmann conditions are too stringent to admit many
solutions for $F$. Thus, we strongly believe, that if $n\geq 7$ the
corresponding $\Gl$ geometries on the solution spaces of the
W\"unschmann ODEs are very special, such that, for example, their
characteristic connections have identically vanishing curvatures. We
intend to discuss these matters in a subsequent paper.
\end{remark}

\end{document}